\documentclass[
12pt,
a4paper,
]{scrbook}

\usepackage{scrpage2}
\usepackage{graphicx}

\usepackage{mathpazo}
\usepackage[scaled]{berasans}
\usepackage[scaled]{beramono}
\usepackage[euler-digits]{eulervm}
\DeclareTextCommandDefault{\nobreakspace}{\leavevmode\nobreak\ }

\usepackage[verbose=silent]{microtype}

\usepackage[utf8]{inputenc}
\usepackage[T1]{fontenc}
\usepackage[english]{babel}

\usepackage{amssymb,amsmath,amsthm}
\usepackage{enumitem}
\usepackage{mathtools}

\usepackage[rgb, table]{xcolor}

\usepackage{kbordermatrix}

\usepackage{booktabs}

\usepackage{tabularx}

\newcommand{\myarraystretch}{1.05}
\newcommand{\temparraystretch}{1.05}

\newcolumntype{I}{!{\vrule width 1.3pt}}
\newlength\savedwidth
\newcommand\whline{\noalign{\global\savedwidth\arrayrulewidth
                            \global\arrayrulewidth 1.3pt}%
           \hline
           \noalign{\global\arrayrulewidth\savedwidth}}

\newcommand{\myFlaTwoByTwo}[4]{
\renewcommand{\temparraystretch}{\arraystretch}
\renewcommand{\arraystretch}{\myarraystretch}
  \left(
	\begin{array}{c@{\;}|@{\;}c} 
      #1 & #2 \\ \hline
      #3 & #4 
    \end{array}
    \renewcommand{\arraystretch}{\temparraystretch}
  \right)
}

\newcommand{\myFlaTwoByOne}[2]{
\renewcommand{\temparraystretch}{\arraystretch}
\renewcommand{\arraystretch}{\myarraystretch}
  \left(
	\begin{array}{c} 
      #1 \\ \hline
      #2  
    \end{array} 
    \renewcommand{\arraystretch}{\temparraystretch}
  \right)
}

\newcommand{\myFlaOneByTwo}[2]{
\renewcommand{\temparraystretch}{\arraystretch}
\renewcommand{\arraystretch}{\myarraystretch}
\left(
\begin{array}{c | c} 
#1 & #2
\end{array}
\renewcommand{\arraystretch}{\temparraystretch}
\right)
}

\newcommand{\myFlaTwoByTwoI}[4]{
\renewcommand{\temparraystretch}{\arraystretch}
\renewcommand{\arraystretch}{\myarraystretch}
\left( 
\begin{array}{c@{\;\;}I@{\;\;}c} 
#1 & #2 \\ \whline
#3 & #4 
\end{array}
\renewcommand{\arraystretch}{\temparraystretch}
\right)
}

\newcommand{\myFlaThreeByThree}[9]{
\renewcommand{\temparraystretch}{\arraystretch}
\renewcommand{\arraystretch}{\myarraystretch}
\left( 
\begin{array}{c | c | c}
#1 & #2 & #3 \\ \hline
#4 & #5 & #6 \\ \hline
#7 & #8 & #9
\end{array}
\renewcommand{\arraystretch}{\temparraystretch}
\right) 
}

\newcommand{\myFlaThreeByTwo}[6]{
\renewcommand{\temparraystretch}{\arraystretch}
\renewcommand{\arraystretch}{\myarraystretch}
\left( 
\begin{array}{c | c}
#1 & #2 \\ \hline
#3 & #4 \\ \hline
#5 & #6
\end{array}
\renewcommand{\arraystretch}{\temparraystretch}
\right) 
}

\newcommand{\FourColHeadTTT}[0]{
\renewcommand{\temparraystretch}{\arraystretch}
\renewcommand{\arraystretch}{\myarraystretch}
\left( 
\begin{array}{c | c | c | c}
}

\newcommand{\FourColHeadFFF}[0]{
\renewcommand{\temparraystretch}{\arraystretch}
\renewcommand{\arraystretch}{\myarraystretch}
\left( 
\begin{array}{c I c I c I c}
}

\newcommand{\FourColHeadFFT}[0]{
\renewcommand{\temparraystretch}{\arraystretch}
\renewcommand{\arraystretch}{\myarraystretch}
\left( 
\begin{array}{c I c I c | c}
}

\newcommand{\FourColHeadTFF}[0]{
\renewcommand{\temparraystretch}{\arraystretch}
\renewcommand{\arraystretch}{\myarraystretch}
\left( 
\begin{array}{c | c I c I c}
}

\newcommand{\FourColHeadFTF}[0]{
\renewcommand{\temparraystretch}{\arraystretch}
\renewcommand{\arraystretch}{\myarraystretch}
\left( 
\begin{array}{c I c | c I c}
}

\newcommand{\FourColHeadTFT}[0]{
\renewcommand{\temparraystretch}{\arraystretch}
\renewcommand{\arraystretch}{\myarraystretch}
\left( 
\begin{array}{c | c I c | c}
}

\newcommand{\FourColHeadTTF}[0]{
\renewcommand{\temparraystretch}{\arraystretch}
\renewcommand{\arraystretch}{\myarraystretch}
\left( 
\begin{array}{c | c | c I c}
}

\newcommand{\FourColTail}[4]{
#1 & #2 & #3 & #4
\end{array}
\renewcommand{\arraystretch}{\temparraystretch}
\right) 
}

\newcommand{\FiveColHeadTTTT}[0]{
\renewcommand{\temparraystretch}{\arraystretch}
\renewcommand{\arraystretch}{\myarraystretch}
\left( 
\begin{array}{c | c | c | c | c}
}

\newcommand{\FiveColHeadFFFT}[0]{
\renewcommand{\temparraystretch}{\arraystretch}
\renewcommand{\arraystretch}{\myarraystretch}
\left( 
\begin{array}{c I c I c I c | c}
}

\newcommand{\FiveColHeadFTFF}[0]{
\renewcommand{\temparraystretch}{\arraystretch}
\renewcommand{\arraystretch}{\myarraystretch}
\left( 
\begin{array}{c I c | c I c I c}
}

\newcommand{\FiveColTail}[5]{
#1 & #2 & #3 & #4 & #5
\end{array}
\renewcommand{\arraystretch}{\temparraystretch}
\right) 
}

\newcommand{\myFlaFourByFourA}[8]{
\renewcommand{\temparraystretch}{\arraystretch}
\renewcommand{\arraystretch}{\myarraystretch}
\left( 
\begin{array}{c | c | c | c}
#1 & #2 & #3 & #4 \\ \hline
#5 & #6 & #7 & #8 \\ \hline
}

\newcommand{\myFlaFourByFourB}[8]{
#1 & #2 & #3 & #4 \\ \hline
#5 & #6 & #7 & #8
\end{array}
\renewcommand{\arraystretch}{\temparraystretch}
\right) 
}

\newcommand{\myFlaThreeByThreeTLI}[9]{
\renewcommand{\temparraystretch}{\arraystretch}
\renewcommand{\arraystretch}{\myarraystretch}
\left( 
\begin{array}{c@{\;\;}|@{\;\;}c@{\;\;}I@{\;\;}c} 
#1 & #2 & #3 \\ \hline
#4 & #5 & #6 \\ \whline
#7 & #8 & #9
\end{array}
\renewcommand{\arraystretch}{\temparraystretch}
\right) 
}

\newcommand{\myFlaThreeByThreeBRI}[9]{
\renewcommand{\temparraystretch}{\arraystretch}
\renewcommand{\arraystretch}{\myarraystretch}
\left( 
\begin{array}{c@{\;\;}I@{\;\;}c@{\;\;}|@{\;\;}c}
#1 & #2 & #3 \\ \whline
#4 & #5 & #6 \\ \hline
#7 & #8 & #9
\end{array}
\renewcommand{\arraystretch}{\temparraystretch}
\right)
}

\newcommand{\myFlaTwoByOneI}[2]{
\renewcommand{\temparraystretch}{\arraystretch}
\renewcommand{\arraystretch}{\myarraystretch}
\left( 
\begin{array}{@{\hspace{1pt}}c@{\hspace{1pt}}}
#1 \\ \whline
#2 
\end{array}
\renewcommand{\arraystretch}{\temparraystretch}
\right)
}

\newcommand{\myFlaOneByTwoI}[2]{
\renewcommand{\temparraystretch}{\arraystretch}
\renewcommand{\arraystretch}{\myarraystretch}
\left( 
\begin{array}{@{\hspace{1pt}}c@{\hspace{2pt}}I@{\hspace{2pt}}c@{\hspace{1pt}}}
#1 & #2 
\end{array}
\renewcommand{\arraystretch}{\temparraystretch}
\right)
}

\newcommand{\myFlaOneByThree}[3]{
\renewcommand{\temparraystretch}{\arraystretch}
\renewcommand{\arraystretch}{\myarraystretch}
\left( 
\begin{array}{c | c | c}
#1 & #2 & #3 
\end{array}
\renewcommand{\arraystretch}{\temparraystretch}
\right)
}

\newcommand{\myFlaOneByThreeI}[3]{
\renewcommand{\temparraystretch}{\arraystretch}
\renewcommand{\arraystretch}{\myarraystretch}
\left( 
\begin{array}{c I c I c}
#1 & #2 & #3 
\end{array}
\renewcommand{\arraystretch}{\temparraystretch}
\right)
}

\newcommand{\myFlaOneByFour}[4]{
\renewcommand{\temparraystretch}{\arraystretch}
\renewcommand{\arraystretch}{\myarraystretch}
\left( 
\begin{array}{c | c | c | c}
#1 & #2 & #3 & #4
\end{array}
\renewcommand{\arraystretch}{\temparraystretch}
\right)
}

\newcommand{\myFlaOneByFourTFF}[4]{
\renewcommand{\temparraystretch}{\arraystretch}
\renewcommand{\arraystretch}{\myarraystretch}
\left( 
\begin{array}{c | c I c I c}
#1 & #2 & #3 & #4
\end{array}
\renewcommand{\arraystretch}{\temparraystretch}
\right)
}

\newcommand{\myFlaOneByThreeLI}[3]{
\renewcommand{\temparraystretch}{\arraystretch}
\renewcommand{\arraystretch}{\myarraystretch}
\left( 
\begin{array}{c@{\;\;}|@{\;\;}c@{\;\;}I@{\;\;}c} 
#1 & #2 & #3 
\end{array}
\renewcommand{\arraystretch}{\temparraystretch}
\right)
}

\newcommand{\myFlaThreeByOne}[3]{
\renewcommand{\temparraystretch}{\arraystretch}
\renewcommand{\arraystretch}{\myarraystretch}
\left( 
\begin{array}{c}
#1 \\ \hline
#2 \\ \hline
#3
\end{array}
\renewcommand{\arraystretch}{\temparraystretch}
\right)
}

\newcommand{\FlaTwoByOneI}[2]{
\renewcommand{\temparraystretch}{\arraystretch}
\renewcommand{\arraystretch}{\myarraystretch}
\left( 
\begin{array}{c}
#1 \\ \whline
#2 
\end{array}
\renewcommand{\arraystretch}{\temparraystretch}
\right)
}

\newcommand{\FlaThreeByOneTI}[3]{
\renewcommand{\temparraystretch}{\arraystretch}
\renewcommand{\arraystretch}{\myarraystretch}
\left( 
\begin{array}{c}
#1 \\ \hline
#2 \\ \whline
#3 
\end{array}
\renewcommand{\arraystretch}{\temparraystretch}
\right) 
}

\newcommand{\FlaThreeByOneBI}[3]{
\renewcommand{\temparraystretch}{\arraystretch}
\renewcommand{\arraystretch}{\myarraystretch}
\left( 
\begin{array}{c}
#1 \\ \whline
#2 \\ \hline
#3 
\end{array}
\renewcommand{\arraystretch}{\temparraystretch}
\right) 
}

\usepackage{ifthen}
\newcommand{\PBefore}{ P_{\mathrm{before}} }
\newcommand{\PAfter} { P_{\mathrm{after}}  }
\newcommand{\PPost}  { P_{\mathrm{post}}   }
\newcommand{\PPre}   { P_{\mathrm{pre}}    }
\newcommand{\PInv}   { P_{\mathrm{inv}}    }
\newcommand{\WSoperation}{
  $[ D, E, F, \ldots ] = {\mathrm{op}}( A, B, C, D, \ldots )$ 
}

\newcommand{\WSpreprocessing}{
  Preprocessing
}

\newcommand{\WSprecondition}{
  $\PPre$
}

\newcommand{\WSpartition}{
}

\newcommand{\WSpartitionsizes}{
  $\dots$
}

\newcommand{\WSinvariant}{
  $\PInv$
}

\newcommand{\WSguard}{ 
  $ G $
}

\newcommand{\WStopofloop}{
  $ \left\{ \left( \mbox{\WSinvariant} \right) \wedge 
    \left( \mbox{\WSguard} \right) \right\} $
}

\newcommand{\WSrepartition}{
  \begin{minipage}[t]{1in}
    \ \\
  \end{minipage}
}

\newcommand{\WSrepartitionsizes}{
  $\dots$
}

\newcommand{\WSbeforeupdate}{
  $\PBefore$ 
}

\newcommand{\WSupdate}{
  \begin{minipage}[t]{.5in}
    $ S_U $
  \end{minipage}
}

\newcommand{\WSafterupdate}{
  $\PAfter$
}

\newcommand{\WSmoveboundary}{
}

\newcommand{\WSafterloop}{
  $\left\{ \left( \mbox{\WSinvariant} \right) \wedge 
   \neg \left( \mbox{\WSguard} \right) \right\} $
}

\newcommand{\WSpostcondition}{
  $\PPost$
}

\newlength{\tabwidth}
\newlength{\tabmargin}
\newcommand{\ALGroutinename}{}

\newboolean{BlockedAlgQ}
\setboolean{BlockedAlgQ}{true}

\newcommand{\worksheetGrayNoNumbersEmpty}{
\begin{center}
\begin{tabular}{| l |} \hline
    Algorithm: \WSoperation \\ \hline \hline
    \rowcolor[gray]{0.83}
    $\left\{ \mbox{\WSprecondition } \right\}$ \\
    {\textbf{Partition}} \\
    \rowcolor[gray]{0.83}
    $\left\{ \mbox{\WSinvariant} \right\}$ \\
    {\textbf{While}} \WSguard { \textbf{do}} \\
    \rowcolor[gray]{0.83}
    \hspace*{0.18in} \WStopofloop \\
    \hspace*{0.15in} 
    {\textbf{Repartition}} \\
    \rowcolor[gray]{0.83}
    \hspace*{0.15in} 
    \ $ \left\{ \mbox{\WSbeforeupdate} \right\}$ \\
    \hspace*{0.15in} \WSupdate \\
    \rowcolor[gray]{0.83}
    \hspace*{0.15in} 
    \ $ \left\{ \mbox{\WSafterupdate} \right\}$ \\
    \hspace*{0.15in} 
    {\textbf{Continue with}} \\
    \rowcolor[gray]{0.83}
    \hspace*{0.21in} 
    $\left\{ \mbox{\WSinvariant} \right\}$ \\
    {\textbf{endwhile}} \\
    \rowcolor[gray]{0.83}
    \WSafterloop \\
    \rowcolor[gray]{0.83}
    $ \left\{ \mbox{\WSpostcondition } \right\} $ \\ \hline
\end{tabular}
\end{center}
}

\newcommand{\worksheetGrayNoNumbersEmptyP}{
\begin{center}
\begin{tabular}{| l |} \hline
    Algorithm: \WSoperation \\ \hline \hline
    \rowcolor[gray]{0.83}
    $\left\{ \mbox{\WSprecondition } \right\}$ \\
    {\textbf{Partition}} \\
    \WSpreprocessing \\
    \rowcolor[gray]{0.83}
    $\left\{ \mbox{\WSinvariant} \right\}$ \\
    {\textbf{While}} \WSguard { \textbf{do}} \\
    \rowcolor[gray]{0.83}
    \hspace*{0.18in} \WStopofloop \\
    \hspace*{0.15in} 
    {\textbf{Repartition}} \\
    \rowcolor[gray]{0.83}
    \hspace*{0.15in} 
    \ $ \left\{ \mbox{\WSbeforeupdate} \right\}$ \\
    \hspace*{0.15in} \WSupdate \\
    \rowcolor[gray]{0.83}
    \hspace*{0.15in} 
    \ $ \left\{ \mbox{\WSafterupdate} \right\}$ \\
    \hspace*{0.15in} 
    {\textbf{Continue with}} \\
    \rowcolor[gray]{0.83}
    \hspace*{0.21in} 
    $\left\{ \mbox{\WSinvariant} \right\}$ \\
    {\textbf{endwhile}} \\
    \rowcolor[gray]{0.83}
    \WSafterloop \\
    \rowcolor[gray]{0.83}
    $ \left\{ \mbox{\WSpostcondition } \right\} $ \\ \hline
\end{tabular}
\end{center}
}

\newcommand{\resetsteps}{
\renewcommand{\ALGroutinename}{}
\renewcommand{\WSoperation}
  { $[ D, E, F, \ldots ] = {\mathrm{op}}( A, B, C, D, \ldots )$ }
\renewcommand{\WSprecondition}{ $\PPre$ }
\renewcommand{\WSpartition}{ }
\renewcommand{\WSpartitionsizes}{}
\renewcommand{\WSinvariant}{ $\PInv$ }
\renewcommand{\WSguard}{ $ G $ }
\setboolean{BlockedAlgQ}{true}
\renewcommand{\WStopofloop}{
  $ \left\{ \left( \mbox{\WSinvariant} \right) \wedge 
    \left( \mbox{\WSguard} \right) \right\} $ }
\renewcommand{\WSrepartition}{
  \begin{minipage}[t]{1in}
    \ \\
  \end{minipage} }
\renewcommand{\WSrepartitionsizes}{}
\renewcommand{\WSbeforeupdate}{ $\PBefore$ }
\renewcommand{\WSupdate}{
  \begin{minipage}[t]{1in}
    $ S_U $
  \end{minipage} }
\renewcommand{\WSafterupdate}{ $\PAfter$ }
\renewcommand{\WSmoveboundary}{}
\renewcommand{\WSafterloop}{
  $\left\{ \left( \mbox{\WSinvariant} \right) \wedge 
   \neg \left( \mbox{\WSguard} \right) \right\} $ }
\renewcommand{\WSpostcondition}{ $\PPost$ }
}

\newcommand{\FlaAlgorithmIter}{
  \begin{center}
    \begin{tabular}{|l|} \hline
      {\textsc Algorithm:} \ALGroutinename \\ \hline
      \hspace*{-3mm}
      \begin{tabular}{ p{5mm} @{\hspace*{0pt}} l }\\[-3mm]
        \multicolumn{2}{l}{\small \textbf{Partition}} \WSpartition  \\
        & {\small \textbf{where}} \WSpartitionsizes  \\[2pt]
      \end{tabular} \\[1mm]
      \WSpreprocessing \\
      {\small \textbf{While}} \WSguard {\small \textbf{do}} \\[1pt] 
      \hspace*{1.9mm}
      \begin{tabular}{ p{3mm} @{\hspace*{0pt}} l}
        \multicolumn{2}{l}{\small \textbf{Repartition}} \\[-1pt]
        & \WSrepartition  \\
        \multicolumn{2}{l}{\rule{5.5cm}{1.2pt}}\\
        & \WSupdate \\[-1mm] 
        \multicolumn{2}{l}{\rule{5.5cm}{1.2pt}}\\
        \multicolumn{2}{l}{\textbf{\small Continue with}} \\[-1pt]
        & \WSmoveboundary \\
      \end{tabular} \\
      {\textbf{endwhile}} \\ 
      \hline
    \end{tabular}
  \end{center}
}

\newcommand{\FlaAlgorithm}{
  \begin{center}
    \begin{tabular}{|l|} \hline
      {\textsc Algorithm:} \ALGroutinename \\ \hline
      \hspace*{-3mm}
      \begin{tabular}{ p{5mm} @{\hspace*{0pt}} l }\\[-3mm]
        \multicolumn{2}{l}{\small \textbf{Partition}} \WSpartition  \\
        & {\small \textbf{where}} \WSpartitionsizes  \\[2pt]
      \end{tabular} \\[1mm]
      {\small \textbf{While}} \WSguard {\small \textbf{do}} \\[1pt] 
      \hspace*{1.9mm}
      \begin{tabular}{ p{3mm} @{\hspace*{0pt}} l}
        \multicolumn{2}{l}{\small 
            \ifthenelse{\boolean{BlockedAlgQ}}
            {\textbf{Determine block size} $b$}
            {}
          } \\[1pt]
        \multicolumn{2}{l}{\small \textbf{Repartition}} \\[-1pt]
        & \WSrepartition  \\
        & {\small \textbf{where}} \WSrepartitionsizes \\[-2mm]
        \multicolumn{2}{l}{\rule{5.5cm}{1.2pt}}\\
        & \WSupdate \\[-1mm] 
        \multicolumn{2}{l}{\rule{5.5cm}{1.2pt}}\\
        \multicolumn{2}{l}{\textbf{\small Continue with}} \\[-1pt]
        & \WSmoveboundary \\
      \end{tabular} \\
      {\textbf{endwhile}} \\ 
      \hline
    \end{tabular}
  \end{center}
}

\newtheoremstyle{mydefinitionstyle}
{}
{}
{\setlength{\parindent}{0cm}}
{}
{\bfseries}
{\\}
{.5em}
{}
\theoremstyle{mydefinitionstyle}
\newtheorem{definition}{Definition}[section]

\newtheoremstyle{mytheoremstyle}
{}
{}
{\setlength{\parindent}{0cm}}
{}
{\bfseries}
{\\}
{.5em}
{}
\theoremstyle{mytheoremstyle}

\newtheoremstyle{mylemmastyle}
{}
{}
{\setlength{\parindent}{0cm}}
{}
{\bfseries}
{\\}
{.5em}
{}
\theoremstyle{mylemmastyle}

\newtheoremstyle{mycorollarystyle}
{}
{}
{\setlength{\parindent}{0cm}}
{}
{\bfseries}
{\\}
{.5em}
{}
\theoremstyle{mycorollarystyle}

\newtheoremstyle{myexamplestyle}
{}
{}
{}
{}
{\itshape}
{\\}
{.5em}
{}
\theoremstyle{myexamplestyle}

\usepackage{tikz}
\usetikzlibrary{positioning}

\definecolor{updateF}{Hsb}{0, 0.4, 1.0}

\definecolor{linvD}{Hsb}{203, 1.0, 0.8}
\definecolor{linvF}{Hsb}{203, 0.4, 0.8}

\tikzset{
	every picture/.style={
		node distance=1.5cm and 1.5cm,
		on grid=true,
		semithick,
		>=latex,
		auto, 
		default/.style={
			circle,
			very thick,
			draw= none,
			fill= lightgray,
		},
		rect/.style={
			rectangle,
			rounded corners,
			very thick,
			draw= none,
			fill= lightgray,
		},
		linv/.style={
			rectangle,
			rounded corners,
			very thick,
			draw= none,
			fill= linvF,
		},
		update/.style={
			rectangle,
			rounded corners,
			very thick,
			draw= none,
			fill= updateF,
		},
	},
	every matrix/.style={
		row sep=1cm,
		column sep= {0.4cm},
	}
}

\usepackage[
]{hyperref} 

\begin{document}


\begin{titlepage}


\begin{tikzpicture}[remember picture, overlay]
	  \node [anchor=north east, inner sep=0pt, xshift=-5mm] at (current page.north east){\includegraphics[height=32mm]{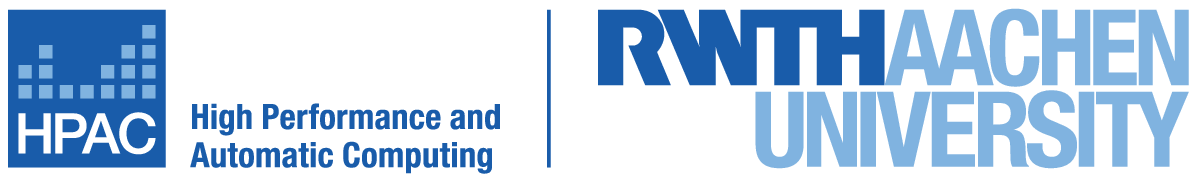}};
\end{tikzpicture}\\[5cm]

\begin{center}
\sffamily
{\Huge \bfseries \sffamily Systematic Generation of Algorithms for Iterative Methods\\[0.5cm]}
{\Large \sffamily Master's Thesis\\[1.5cm]}
{\Large \sffamily Henrik Barthels, B.Sc.\\}
{\vfill}
{Supervised by\\
Prof. Paolo Bientinesi, Ph.D.\\
Prof. Georg May, Ph.D.}
\end{center}
\end{titlepage}

\frontmatter

\cleardoublepage

\subsection*{Version History}

This is a revised edition of the author's thesis.

\subsubsection*{First Revision}

Corrections of typographical errors and clarifications of some passages.

\subsubsection*{Second Revision}

Addition of some references.

\cleardoublepage

\subsection*{Abstract}

The FLAME methodology makes it possible to derive provably correct algorithms from a formal description of a linear algebra problem. So far, the methodology has been successfully used to automate the derivation of direct algorithms such as the Cholesky decomposition and the solution of Sylvester equations. In this thesis, we present an extension of the FLAME methodology to tackle iterative methods such as Conjugate Gradient. As a starting point, we use a formal description of the iterative method in matrix form. The result is a family of provably correct pseudocode algorithms. We argue that all the intermediate steps are sufficiently systematic to be fully automated.

\cleardoublepage
\subsection*{Acknowledgments}
\tikz[remember picture] \node[scale=0.01] (ack1) {};

\noindent
I would like to thank Paolo Bientinesi and Diego Fabregat Traver for their support, guidance, and taking their time for all those lengthy discussions. I want to thank Georg May for agreeing to be the second supervisor of this thesis. Last but not least, I thank Friederike and Robert for proofreading.

\begin{tikzpicture}[remember picture,overlay]

\node (ack11) at (ack1) [xshift=1.845cm, yshift=0.05pt, , yshift=1cm] {};
\node [scale=0.006, rotate=90, color=white] at (ack11) {47847542435635809235762872377575944927660000 83383191170664050450142445419123365066167050916870919608156195489913826093750 6612949327291548436115306813438614139140599264303491211530114314015726475397706073484320072965158098374168300344235229928198676088291682518718750};
\end{tikzpicture}

\tableofcontents



\clearpage
\section*{Glossary}
\markboth{\MakeUppercase{Glossary}}{}

\begin{table}[htp]
\begin{center}
\begin{tabular}{l p{9cm}}
$A$, $B$, $C$, …		&	Matrices. \\
$a$, $b$, $c$, …		&	Row or column vectors. \\
$\alpha$, $\beta$, $\gamma$, …	&	Scalars. \\
BiCG				&	Biconjugate Gradient method.\\
BLAS				& 	Basic Linear Algebra Subprograms. A library for linear algebra operations.\\
CG					& 	Conjugate Gradient method.\\
$e_0$				& 	The unit vector that is one in the first position.\\
$e_r$				& 	The unit vector that is one in the last position.\\
FLAME				& 	Formal Linear Algebra Methods Environment.\\
$G$					&	Loop guard. \\
$I$					&	Identity matrix. \\
$P_\text{after}$			&	After predicate. \\
$P_\text{before}$		&	Before predicate. \\
$P_\text{inv}$			&	Loop invariant. \\
PME					&	Partitioned Matrix Expression. \\
$P_\text{post}$			&	Postcondition. \\
$P_\text{pre}$			&	Precondition. \\
SPD					&	Symmetric Positive Definite. \\
$\| \ldots \|$			&	An arbitrary vector norm.
\end{tabular}

\end{center}
\end{table}

\mainmatter

\chapter{Introduction}

The goal of this thesis is to simplify the development of algorithms for iterative methods. We present a methodology for the systematic derivation of such algorithms and lay the foundations for a system that automates the generation of algorithms and code for iterative solvers for linear systems.


Over the last few decades, iterative methods have become an indispensable tool for solving sparse linear systems. Such systems commonly occur in science and engineering, for instance when discretized partial differential equations have to be solved. While direct methods are a reliable tool to solve linear systems, their ability to use the sparsity of a matrix to their advantage is limited to specific sparsity patterns.
Often enough, those sparse systems are so large that using direct methods is impractical.

Since the introduction of the Conjugate Gradient (CG) method in 1952 \cite{hestenes1952methods}, much progress has been made in the field (for an overview, consider \cite{barrett:templates, saad2000iterative, vanderVorst:book}).
%
%
%
%
However, the way from an expert's idea for an algorithm to a working implementation is still a long one.
At first, a formal description of the algorithm has to be derived on paper. Then, it has to be shown that the algorithm is correct, something that ideally follows from the derivation. Furthermore, to assess how useful the algorithm is in practice, it is desirable to prove that it is numerically stable. Finally, the algorithm has to be translated into code. Usually, that is done multiple times, for different languages, or potentially using different libraries. Every one of those steps takes time and is a potential source of errors.

To speed up this process and eliminate those sources of error, one may try to automate some or all of the steps described above. Automating the ingenuity of the iterative methods expert certainly lies in the distant future. In contrast, automatically translating a sufficiently formal description of an algorithm into code is a lot more feasible.

This thesis covers the systematical derivation of provably correct (pseudocode) algorithms, based on an abstract, formal description of an iterative method.


%
%
\section{Background: FLAME}

The Formal Linear Algebra Methods Environment (FLAME) \cite{Bientinesi2005:504, gunnels2001flame} is a project with the goal to automate the derivation of linear algebra algorithms, as well as their implementation.
In \cite{Bientinesi:thesis}, it was shown that it is possible to systematically derive algorithms for dense linear algebra in a number of well defined steps.

The starting point of this derivation is a formal description of the input, consisting of a precondition ($P_\text{pre}$) and a postcondition ($P_\text{post}$). As an example, the description of a linear system where A is lower triangular is shown below:
\begin{align*}
x:= \Phi ( A, b) \equiv
\left\{
\begin{aligned}
P_\text{pre}: \{ &\text{\ttfamily Input}[A] \land \text{\ttfamily Matrix}[A] \land \text{\ttfamily NonSingular}[A] \land \\
		&\text{\ttfamily LowerTriangular}[A] \land \\
		&\text{\ttfamily Input}[b] \land \text{\ttfamily Vector}[b] \land \\
		&\text{\ttfamily Output}[x] \land \text{\ttfamily Vector}[x] \} \\
P_\text{post}: \{ & Ax = b \}
\end{aligned}
\right. 
\end{align*}
The function $x:= \Phi (A, b)$ is used to abstract from the details of this representation.

At the core of the derivation lies the concept of a Partitioned Matrix Expression (PME). A PME describes all parts of the output operands of a linear algebra operation in terms of parts of the input operands. 
For the lower triangular system, the equation $Ax=b$ is partitioned as
$$\myFlaTwoByTwo{A_{TL}}{ 0 }{A_{BL}}{A_{BR}} \myFlaTwoByOne{x_T}{x_B} = \myFlaTwoByOne{b_T}{b_B}\text{,}$$
where $A_{TL}$ and $A_{BR}$ are square. Then, the PME is
$$\myFlaTwoByOne{x_T := \Phi \left( A_{TL}, b_T \right)}{x_B := \Phi \left( A_{BR}, b_B - A_{BL} x_T \right)} \text{.}$$
%
%
%
%
The approach described in \cite{Bientinesi:thesis} has a number of advantages: 
(1) The resulting algorithms are built around a proof of correctness, so they are correct by construction. (2) It naturally leads to multiple variants of algorithms for the same operation. While they are all correct in exact arithmetic, they potentially behave very different in presence of round-off errors. (3) Furthermore, it was shown that this derivation can be combined with a systematic stability analysis \cite{Bientinesi2011:810}.

The author of \cite{Bientinesi:thesis} presented evidence that this derivation, based on a formal description of the operation and the PME, is systematic enough to be automated. A system that automates this process, including the generation of PMEs, was presented in \cite{Fabregat-Traver:thesis, Fabregat-Traver2011:238}.

All of the efforts described above focused on direct methods. Naturally, the question arises if this methodology extends to iterative methods. In \cite{eijkhout:CGderivation}, Eijkhout et al.\ introduced a matrix representation of the CG method and showed preliminary evidence that a systematic derivation of algorithms with a FLAME-like methodology is possible. However, the approach presented there heavily relied on guidance by a human expert. This thesis can be seen as a continuation of this work.

\section{Challenges}

The systematic derivation of algorithms for iterative methods introduces new challenges, especially  for deriving the PME.
\begin{itemize}
\item[-] For direct methods, the sizes of all operands are constant. In case of iterative methods, a variable number of iterations can be performed. In the matrix representation, this is reflected by the fact that some operands have variable sizes.
\item[-] The matrix representation of iterative methods introduces new types of op\-er\-ands. Usually, an operand is either known or unknown, so it is input or output, respectively. Now, there are operands that are initially partially known and partially unknown.
\item[-] Quite often, equations have to be solved by using properties of certain expressions. Consider the following equation as an example:
$$-P u + p = r$$
$P$, $p$ and $r$ are known, and the goal is to find an assignment for $u$. It can be solved by using the fact that $P^T A P$ is lower triangular and $P^T A p$ is zero ($A$ is known as well). Multiplying $P^T A$ from the left to both sides of the equation results in
$$-P^T A P u = P^T A r \text{,}$$
which is a triangular system that can easily be solved.

This introduces two challenges: The properties of $P^T A P$ and $P^T A p$ are not explicitly part of the matrix representation, so they have to be derived from it by algebraic manipulations and deductive reasoning. Then, to enable a system to solve this equation, it must be capable of identifying that  the properties of those expressions can be used to do so.
\item[-] To ensure that the derived algorithms can be used in practice, it is desirable to derive the exact same algorithms as used today. The difficulty in achieving this goal lies in the fact that some of those algorithms are the result of nontrivial rewritings of easily derivable formulas. Those rewritings are hard to automate.
\end{itemize}

\section{Contributions}

The following contributions are made in this thesis:

\begin{description}
\item[Systematic derivation of algorithms for iterative methods.] We present an approach for sys\-tematically deriving algorithms for iterative methods similar to the one discussed above for direct methods. To achieve this, we present solutions to the following two problems:

\begin{description}

\item[Systematic derivation of matrix properties.] To be able to derive algorithms, properties of matrices or matrix expressions are needed that are not explicitly part of the formal description of the operation. We present a method for deriving those properties from the description, using a set of inference rules that encode linear algebra knowledge.

\item[Systematic derivation of PMEs for iterative methods.] The derivation of PMEs for iterative methods is more complex compared to direct methods. As part of this derivation, equations have to be solved using the derived matrix properties.

\end{description}
\end{description}

The presented methodology is systematic enough to be executed mechanically, that is, without any human intervention, thus setting the ground for a system that indeed automates its application.

Having made clear what is accomplished with this thesis, it should also be pointed out what is beyond its scope: The input to the derivation process is a representation of an iterative method that still has to be derived by an expert in the field. In case of direct methods, the formal description that is input immediately follows from the operation itself, as shown above for the lower triangular system. 
In contrast, the way from a linear system $Ax=b$ to a formal description of an iterative method that is suitable as input for the presented approach is much more complicated.

\section{Outline of the Thesis}

The thesis is structured as follows: Chapter \ref{chap:directMethods} serves as an introduction to familiarize the reader with the FLAME approach for deriving algorithms for direct methods. In Chapter \ref{chap:iterativeMethodsIntro}, we lay the foundations for applying a similar approach to iterative methods. The matrix representation for iterative methods is introduced in Section \ref{sec:matrixRepresentationIntroduction}. In Section \ref{sec:propertyDerivation}, we describe a method for systematically deriving matrix properties from this representation. At the end of this chapter, in Section \ref{sec:MatrixRepresentationPartitionings}, the implications of the matrix representation for the derivation of algorithms are discussed. In Chapter \ref{chap:derivationIM}, we explain how to derive algorithms for iterative methods. Finally, Chapter \ref{chap:conclusion} summarizes the results of this thesis and points out opportunities for future research.

%
%

%
%

\chapter{Derivation of Algorithms for Direct Methods}
\label{chap:directMethods}

The purpose of this chapter is to introduce the reader to the systematic derivation of loop-based, blocked algorithms for direct methods using the FLAME methodology. For the derivation, we will follow the systematic approach and the notation of \cite{Fabregat-Traver:thesis}. We begin with the simple example of a triangular linear system in Section \ref{sec:triLS}. In Section \ref{sec:SPDsystem}, we proceed with a somewhat more elaborate example to go into some details that are not covered in the first one. Next, in Section \ref{sec:generalLS}, a case will be demonstrated where it is not possible to derive algorithms that compute the solution for every input of the operation. Finally, we discuss the equivalence of loop invariants that are obtained with this approach in Section \ref{sec:equivalence}.

%

\subsubsection{Notation}

Throughout this thesis, two different notations are used for the indexing of matrix blocks and elements. The first one is the standard notation that uses numerals. The second one uses capital letters, as shown below:
\begin{align*}
& \bullet \myFlaOneByTwo{A_L}{A_R} 
&& \bullet \myFlaTwoByOne{A_T}{A_B}
&& \bullet \myFlaTwoByTwo{A_{TL}}{A_{TR}}{A_{BL}}{A_{BR}} \\
& \bullet \myFlaOneByThree{A_L}{A_M}{A_R}
&& \bullet \myFlaThreeByOne{A_T}{A_M}{A_B}
&& \bullet \myFlaThreeByThree	{A_{TL}}	{A_{TM}}	{A_{TR}}
					{A_{ML}}	{A_{MM}}	{A_{MR}}
					{A_{BL}}	{A_{BM}}	{A_{BR}}
\end{align*}


%
%
The subscript letters $T$, $B$, $L$, $M$ and $R$ stand for Top, Bottom, Left, Middle and Right, respectively.

\section{Triangular Linear System}
\label{sec:triLS}

We begin with deriving algorithms for the linear system $Ax = b$ where $A$ is triangular.
The starting point of the derivation is a complete description of the problem, consisting of two logical predicates, a precondition ($P_\text{pre}$) and a postcondition ($P_\text{post}$). The precondition lists the properties of all quantities that are part of the operation, while the postcondition consists of the equation (or equations) that constitute the operation. The precondition is true before the execution of the algorithm, and our goal is to find an algorithm that makes the postcondition true upon termination. Additionally, a function is introduced which abstracts from the details of the operation, and instead describes the output as a function of the input. For reasons that will become apparent later, we use $\hat{b}$ to denote the initial content of $b$. The description for a linear system with $A$ lower triangular is:
$$x:= \Phi ( A, \hat{b}) \equiv
\left\{
\begin{aligned}
P_\text{pre}: \{ &\text{\ttfamily Input}[A] \land \text{\ttfamily Matrix}[A] \land \text{\ttfamily NonSingular}[A] \land \\
		&\text{\ttfamily LowerTriangular}[A] \land \\
		&\text{\ttfamily Input}[\hat{b}] \land \text{\ttfamily Vector}[\hat{b}] \land \\
		&\text{\ttfamily Output}[x] \land \text{\ttfamily Vector}[x] \} \\
P_\text{post}: \{ & Ax = \hat{b} \}
\end{aligned}
\right.$$
The actual algorithms are then constructed by filling out a so-called ``worksheet'', a template for a loop-based algorithm \cite{Bientinesi:thesis, Bientinesi2005:460, Fabregat-Traver:thesis}, shown in Figure \ref{fig:ws:empty}.

\begin{figure}
\centering
\begin{minipage}[t]{2.35in}
	\resetsteps
	\setboolean{BlockedAlgQ}{false}
	
	\renewcommand{\WSoperation}{ $\ldots$}
	
	\renewcommand{\WSupdate}{\text{Update}}
	
	{
	\worksheetGrayNoNumbersEmpty
	}
\end{minipage}
\caption{The skeleton of a FLAME worksheet.}
\label{fig:ws:empty}
\end{figure}

\subsubsection{PME Generation}

As a first step towards an algorithm, the PME for the operation is generated. In order to do so, we partition the postcondition. For direct methods, partitioning operands along each dimension in at most two parts proved to be sufficient (this is further discussed in Section \ref{sec:equivalence}). 
Since $A$ is lower triangular, and to preserve the triangularity of the resulting objects, a $2 \times 2$ partitioning is chosen for it where $A_{TL}$ (and $A_{BR}$) is square.
%
$$\myFlaTwoByTwo{A_{TL}}{ 0 }{A_{BL}}{A_{BR}} \myFlaTwoByOne{x_T}{x_B} = \myFlaTwoByOne{\hat{b}_T}{\hat{b}_B}$$
%
Next, in the so called \emph{Matrix Arithmetic} step, we compute the symbolic multiplication and distribute the equality over the partitioned objects, resulting in the following expression:
$$\myFlaTwoByOne{A_{TL} x_T = \hat{b}_T}{A_{BL} x_T + A_{BR} x_B = \hat{b}_B}$$
In the final \emph{Pattern Matching} step, we try to find ways to solve the given equations by writing them as known operations. We begin with the top part. $A_{TL}$ and $b_T$ are known, while $x_T$ is unknown. Since $A_{TL}$ is lower triangular, we recognize that $A_{TL} x_T = \hat{b}_T$ is a lower triangular system. Thus, we can rewrite this expression as  $x_T := \Phi \left( A_{TL}, \hat{b}_T \right)$, utilizing the function introduced earlier. Doing so, we can from now on consider $x_T$ as known, which is necessary to solve the second equation. $A_{BL} x_T$ is a matrix-vector product of two known quantities, and if we subtract it on both sides of the equation, we identify $A_{BR} x_B = \hat{b}_B - A_{BL} x_T$ as another triangular system. Rewriting this as $x_B := \Phi \left( A_{BR}, \hat{b}_B - A_{BL} x_T \right)$ provides us with the PME
$$\myFlaTwoByOne{x_T := \Phi \left( A_{TL}, \hat{b}_T \right)}{x_B := \Phi \left( A_{BR}, \hat{b}_B - A_{BL} x_T \right)} \text{.}$$
\subsubsection{Loop Invariant Identification}

The next step consists of finding loop invariants ($P_{\text{inv}}$). A loop invariant is a logical predicate that is true at certain points of a loop in an algorithm. It is true before the loop is entered and after it is left, as well as at the beginning and the end of the loop body. It allows to formally reason about the correctness of the loop \cite{GrSc:92}.

\begin{figure}[h]
\centering
\begin{tikzpicture}

\node[default]	(A)							{$1$};

\node[default]	(B)	[below=of A, yshift=-0.0cm]	{$2$};

\node[default]	(C)	[right=of B, yshift=-0.0cm]		{$3$};

\path[-]	(-0.75cm, -0.75cm) edge [lightgray]  (2.25cm, -0.75cm);

\path[->]	(A)		edge							(B);
		
\path[->]	(B)		edge							(C);

\end{tikzpicture}
\caption{Dependency graph for a lower triangular system.}
\label{fig:dg:loTriSys}
\end{figure}
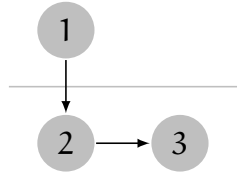

As a first step, the PME is decomposed into basic operations, for example a matrix multiplication, or the function representing the operation we are deriving algorithms for. In general, one can think of those basic operations in terms of BLAS-like functions. For the lower triangular system, the following three operations are obtained:
\begin{enumerate}
\item $x_T := \Phi \left( A_{TL}, \hat{b}_T \right)$
\item $b_B := \hat{b}_B - A_{BL} x_T$
\item $x_B := \Phi \left( A_{BR}, b_B \right)$
\end{enumerate}
Then, from those operations, a dependency graph is constructed. Every node represents one operation, while the edges of this directed graph represent the data dependencies between those operations. If one operation requires the output of another operation, the former depends on the latter. Thus, in the graph, there are edges from the node that computes a quantity to those which need this quantity as input. Consider Figure \ref{fig:dg:loTriSys} for the resulting dependency graph.
%

Subsets of nodes from this graph are now selected as loop invariants. For reasons that will be explained later, those subsets can neither be empty, nor can they contain all nodes. Furthermore, if a node is part of the subset, then all preceding nodes are also in the set. It follows that there are two subsets of nodes of the graph, $\{1\}$ and $\{1, 2\}$, that can be used as loop invariants. 
%
%
\begin{align*}
\{1\}: P_\text{inv}^{1} = &\left\{ x_T := \Phi \left( A_{TL}, \hat{b}_T \right) \right\} \\
\{1, 2\}: P_\text{inv}^{2} = &\left\{ x_T := \Phi \left( A_{TL}, \hat{b}_T \right) \land b_B := \hat{b}_B - A_{BL} x_T \right\}
\end{align*}
%
%
Now, we can also identify the loop guard $G$. The loop is supposed to terminate when all operands are traversed and the complete solution is computed. Thus, the loop guard depends on how the operands are traversed. The first assignment of the dependency graph operates on $x_T$, $A_{TL}$ and $b_T$, so we can conclude that $x$ and $b$ are traversed from the top to the bottom, and $A$ is traversed from the top-left to the bottom-right corner. In the algorithm, this has the result that initially, $A_{TL}$, $x_T$ and $b_T$ are empty. $A_{TL}$ grows with every iteration until it has the same size as $A$. Similarly, $x_T$ and $b_T$ grow to the size of $x$ and $b$. Hence, one possible loop guard is ``$\text{size}(A_{TL}) < \text{size}(A)$''. How $A_{TL}$, $x_T$ and $b_T$ change their sizes is explained in further detail in the next section.

\subsubsection{Algorithm Construction}

Having found a number of loop invariants, we can proceed to the final step, the derivation of the updates, using the worksheet (Figure \ref{fig:ws:empty}). For each loop invariant, a separate worksheet is filled out. Every loop invariant is rewritten in two different ways, resulting in two predicates per loop invariant. The first predicate specifies the state of all operands before the update ($P_\text{before}$). Analogously, the second predicate $P_\text{after}$ represents the situation after the update. The update is then found by identifying the operations that transform the predicate $P_\text{before}$ into $P_\text{after}$. Intuitively speaking, the difference between both predicates is determined.

The rewritings that are applied to the loop invariant are repartitionings of the already partitioned operands. They are chosen in a way that ensures that the operation makes progress and eventually terminates in combination with a suitable loop guard: Those blocks of the output that are already computed grow, and those parts that are not computed yet shrink. This is done by splitting off parts of some quantities by applying the ``Repartition'' rules and merging them with others with the ``Continue with'' repartitioning.
%

We begin with repartitioning the partitioned operands. 
To obtain a blocked algorithm, we use the following ``Repartition'' rules:
\begin{align*}
\myFlaTwoByTwoI{A_{TL}}{ 0 }{A_{BL}}{A_{BR}} &\rightarrow
\myFlaThreeByThreeBRI	{A_{00}}	{0}		{0}
					{A_{10}}	{A_{11}}	{0}
					{A_{20}}	{A_{21}}	{A_{22}} \\
\FlaTwoByOneI{x_T}{x_B} &\rightarrow \FlaThreeByOneBI{x_{0}}{x_1}{x_2} \\
\FlaTwoByOneI{b_T}{b_B} &\rightarrow \FlaThreeByOneBI{b_{0}}{b_1}{b_2}
\end{align*}
Doing so, during every iteration, some parts of $A_{BR}$ are split off. To continue after the update, those parts are merged into $A_{TL}$. Thus, $A_{BR}$ shrinks and $A_{TL}$ grows until the entire matrix, and at the same time both vectors $b$ and $x$, are traversed. This ``Continue with'' partitioning is shown below.
\begin{align*}
\myFlaTwoByTwoI{A_{TL}}{ 0 }{A_{BL}}{A_{BR}} &\leftarrow
\myFlaThreeByThreeTLI	{A_{00}}	{0}		{0}
					{A_{10}}	{A_{11}}	{0}
					{A_{20}}	{A_{21}}	{A_{22}} \\
\FlaTwoByOneI{x_T}{x_B} &\leftarrow \FlaThreeByOneTI{x_{0}}{x_1}{x_2} \\
\FlaTwoByOneI{b_T}{b_B} &\leftarrow \FlaThreeByOneTI{b_{0}}{b_1}{b_2}
\end{align*}
To get to the predicates $P_\text{before}$ and $P_\text{after}$, we have to rewrite the loop invariant using the repartitioned operands. We demonstrate this for the second loop invariant, $P_\text{inv}^{2}$.
%
%
For $P_\text{before}$, we replace all quantities by their counterparts according to the first partitioning (from $2 \times 2$ to $3 \times 3$), resulting in
\begin{align*}
x_0 &:= \Phi \left( A_{00}, \hat{b}_0 \right) \\
\myFlaTwoByOne{b_1}{b_2} &:= \myFlaTwoByOne{\hat{b}_1}{\hat{b}_2} - \myFlaTwoByOne{A_{10}}{A_{20}} x_0 \text{.} 
\end{align*}
Then, we flatten the expression, that is, perform all algebraic operations, distribute the assignments and decompose the result into its parts. This yields the following predicate:
\begin{align*}
P_\text{before} = \{	x_0 &:= \Phi \left( A_{00}, \hat{b}_0 \right) \land \\
				b_1 &:= \hat{b}_1 - A_{10} x_0 \land \\
				b_2 &:= \hat{b}_2 - A_{20} x_0 \}
\end{align*}
The same is done using the ``Continue with'' partitioning. To repartition the function $\Phi$, the PME is used.
\begin{align*}
\myFlaTwoByOne{x_0}{x_1} &:= \Phi \left( \myFlaTwoByTwo{A_{00}}{0}{A_{10}}{A_{11}}, \myFlaTwoByOne{\hat{b}_0}{\hat{b}_1} \right) \\
b_2 &:= \hat{b}_2 - \myFlaOneByTwo{A_{20}}{A_{21}} \myFlaTwoByOne{x_0}{x_1}
\end{align*}
Thus, we obtain the predicate below:
\begin{align*}
P_\text{after} = \{ 	x_0 &:= \Phi \left( A_{00}, \hat{b}_0 \right) \land \\
				x_1 &:= \Phi \left( A_{11}, \hat{b}_1 - A_{10} x_0 \right) \land \\
				b_2 &:= \hat{b}_2 - A_{20} x_0  - A_{21} x_1 \}
\end{align*}
Now, by comparing the predicates $P_\text{before}$ and $P_\text{after}$, we determine the update. Highlighted in red are those parts of $P_\text{after}$ that do not appear in $P_\text{before}$.
\begin{align*}
P_\text{after} = \{ 	&x_0 := \Phi ( A_{00}, \hat{b}_0 ) \land \\
				&\textcolor{red}{x_1 := \Phi} \mathopen{\textcolor{red}{(}}  \textcolor{red}{A_{11}} \mathpunct{\textcolor{red}{,}} \hat{b}_1 - A_{10} x_0 \mathclose{\textcolor{red}{)}} \land \\
				&b_2 := \hat{b}_2 - A_{20} x_0 \mathbin{\textcolor{red}{-}} \textcolor{red}{A_{21} x_1} \}
\end{align*}
This provides us with the following update:
\begin{align*}
x_1 &:= \Phi \left( A_{11}, b_1 \right) \\
b_2 &:= b_2 - A_{21} x_1
\end{align*}
We have now derived a complete algorithm. Consider Figure \ref{fig:ws:triLS} for the filled out worksheet with updates for both loop invariants (omitting some of the repartitionings and all predicates in the interest of visual clarity).

Using the loop guard and the loop invariants, it is easy to see that the derived algorithms are correct and at the end of the computation, the linear system is solved. The algorithms terminate when the loop guard becomes false. The negation of ``$\text{size}(A_{TL}) < \text{size}(A)$'' implies that $A_{TL}$ has the same size as $A$, as it can not be larger. Since $A_{TL}$ is a part of $A$, this means that $A_{TL}$ is equal to $A$. Similarly, $x_T$ equals $x$. $\hat{b}_T$ and $b_T$ are equal to $\hat{b}$ and $b$, respectively. As we updated $b$, we have to distinguish between its initial and its current content to reason about the correctness of the derived algorithm. All other parts have either size $0 \times n$, $0 \times 0$ or $n \times 0$, so they disappear. Plugging that in the loop invariants, in both cases we only get $x := \Phi \left( A, \hat{b} \right)$, which proves that both algorithms compute the solution to a triangular system.

\begin{figure}
\centering
\begin{minipage}[t]{4.5in}
	\resetsteps
	\setboolean{BlockedAlgQ}{false}
	
	\renewcommand{\ALGroutinename}{ $x := \Phi \left( A, b \right)$}
	
	\renewcommand{\WSprecondition}{
	
	}
	
	\renewcommand{\WSpartition}{
        $
        A \rightarrow
		\myFlaTwoByTwo{A_{TL}}{ 0 }{A_{BL}}{A_{BR}}
	$
        }

	\renewcommand{\WSpartitionsizes}{
		$ A_{TL} $ is $ 0 \times 0 $
	}
	
	\renewcommand{\WSguard}{ $\text{size}(A_{TL}) < \text{size}(A)$ }

	\renewcommand{\WSrepartition}{
	$\myFlaTwoByTwoI{A_{TL}}{ 0 }{A_{BL}}{A_{BR}} 
	\rightarrow
	\myFlaThreeByThreeBRI	{A_{00}}	{0}		{0}
						{A_{10}}	{A_{11}}	{0}
						{A_{20}}	{A_{21}}	{A_{22}}
	$
	}
	
	\renewcommand{\WSrepartitionsizes}{
	$A$ is $k \times k$
	}

	\renewcommand{\WSupdate}{
	\begin{tabular}{c p{0.5cm} c}
	Variant 1 & & Variant 2 \\[4pt]
	$\begin{aligned}
	b_1 &:= \hat{b}_1 - A_{10} x_0 \\
	x_1 &:= \Phi \left( A_{11}, b_1 \right)
	\end{aligned}$
	& &
	$\begin{aligned}
	x_1 &:= \Phi \left( A_{11}, b_1 \right) \\
	b_2 &:= b_2 - A_{21} x_1
	\end{aligned}$
	\end{tabular}
	}
	
	\renewcommand{\WSmoveboundary}{
	$\myFlaTwoByTwoI{A_{TL}}{ 0 }{A_{BL}}{A_{BR}} \leftarrow
	\myFlaThreeByThreeTLI	{A_{00}}	{0}		{0}
					{A_{10}}	{A_{11}}	{0}
					{A_{20}}	{A_{21}}	{A_{22}}$
	}
	
	\renewcommand{\WSpostcondition}{
	
	}
	
	{
	\FlaAlgorithm
	}
\end{minipage}
\caption{Worksheet for a lower triangular linear system.}
\label{fig:ws:triLS}
\end{figure}

%

\section{Symmetric Positive Definite Linear System}
\label{sec:SPDsystem}

\begin{figure}[h]
\centering
\begin{tikzpicture}

\node[default]	(A)							{$1$};

\node[default]	(B)	[above=of A, yshift=-0.0cm]	{$2$};

\node[default]	(C)	[right=of A, yshift=-0.0cm]		{$3$};

\node[default]	(D)	[right=of B, yshift=-0.0cm]		{$4$};

\node[default]	(E)	[above=of D, yshift=-0.0cm]	{$5$};

\node[default]	(F)	[above=of E, yshift=-0.0cm]	{$6$};

\path[-]	(-0.75cm, 2.25cm) edge [lightgray]  (2.25cm, 2.25cm);

\path[->]	(A)		edge							(B);
\path[->]	(A)		edge							(C);
\path[->]	(B)		edge							(D);
\path[->]	(C)		edge							(D);
\path[->]	(D)		edge							(E);
\path[->]	(E)		edge							(F);

\end{tikzpicture}
\caption{Dependency graph for a SPD system.}
\label{fig:dg:SPDSys}
\end{figure}
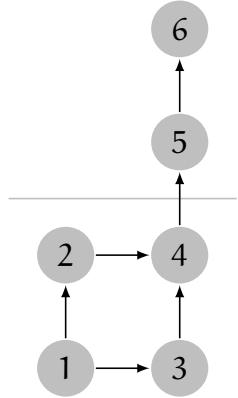

As a second example, we derive algorithms for a symmetric positive definite (SPD) linear system to show in greater detail how the feasibility of loop invariants is checked. Again, the starting point is the formal description of the operation, which is shown below: 
$$x:= \Sigma \left( A, b\right) \equiv
\left\{
\begin{aligned}
P_\text{pre}: \{ &\text{\ttfamily Input}[A] \land \text{\ttfamily Matrix}[A] \land \text{\ttfamily SPD}[A] \land \\
		&\text{\ttfamily Input}[b] \land \text{\ttfamily Vector}[b] \land \\
		&\text{\ttfamily Output}[x] \land \text{\ttfamily Vector}[x] \} \\
P_\text{post}: \{ & \hat{A}x = \hat{b} \}
\end{aligned}
\right.$$
Because of the symmetry of $\hat{A}$, we apply a $2 \times 2$ partitioning where $\hat{A}_{TL}$ (and $\hat{A}_{BR}$) is square. Furthermore, because of the symmetry, it holds that $\hat{A}_{BL} = \hat{A}_{TR}^T$. The following partitioned postcondition is obtained.
$$\myFlaTwoByOne{\hat{A}_{TL} x_T + \hat{A}_{TR} x_B  = \hat{b}_T}{\hat{A}_{BL} x_T + \hat{A}_{BR} x_B = \hat{b}_B}$$
Unfortunately, none of the two equations matches the function $x:= \Sigma \left( A, b\right)$. If we rewrite $\hat{A}_{TL} x_T + \hat{A}_{TR} x_b  = \hat{b}_T$ as $\hat{A}_{TL} x_T = \hat{b}_T - \hat{A}_{TR} x_b $, the right-hand side of this equation is not completely known, as $x_b$ is not known. For the same reason, the second equation does not match the function either. In order to proceed, we need to perform a number of steps that are not part of the systematic approach presented in \cite{Fabregat-Traver:thesis}. As $\hat{A}$ is SPD and $\hat{A}_{TL}$ is square, we know that $\hat{A}_{TL}$ (and $\hat{A}_{BR}$) is also SPD, so we can rewrite the top part of the partitioned postcondition as:
$$x_T = \hat{A}_{TL}^{-1} \left( \hat{b}_T - \hat{A}_{TR} x_B \right)$$
By replacing $x_T$ in the second equation with the right-hand side of the equation above and performing further manipulations, it is possible to obtain an equation that matches the pattern of a SPD linear system. Since $\hat{A}$ is SPD, we can infer that $\hat{A}_{BR} - \hat{A}_{BL} \hat{A}_{TL}^{-1} \hat{A}_{TR}$ is SPD too \cite{stew:98a}.
$$\left( \hat{A}_{BR} - \hat{A}_{BL} \hat{A}_{TL}^{-1} \hat{A}_{TR} \right) x_B = \hat{b}_B - \hat{A}_{BL} \hat{A}_{TL}^{-1} \hat{b}_T$$
Thus, the following PME is derived:
$$\myFlaTwoByOne{x_T := \Sigma \left(\hat{A}_{TL}, \hat{b}_T - \hat{A}_{TR} x_B \right)}{x_B := \Sigma \left( \hat{A}_{BR} - \hat{A}_{BL} \hat{A}_{TL}^{-1} \hat{A}_{TR} , \hat{b}_B - \hat{A}_{BL} \hat{A}_{TL}^{-1} \hat{b}_T \right)}$$
Libraries for linear algebra operations like BLAS usually do not offer operations for products of more than two quantities. Thus, to decompose this PME into its basic buildings blocks, we need to introduce auxiliary variables. In general, there are mutliple ways to compute expressions like $\hat{A}_{BL} \hat{A}_{TL}^{-1} \hat{A}_{TR}$, for example by solving $\hat{A}_{BL} \hat{A}_{TL}^{-1}$ or $\hat{A}_{TL}^{-1} \hat{A}_{TR}$ first. To ensure that the dependency graph is as general as possible and does not impose an ordering on expressions like these, multiple auxiliary variables would have to be introduced. Furthermore, their sizes would have to be left unspecified. To keep this example simple, we compute $\hat{A}_{BL} \hat{A}_{TL}^{-1}$ first, because it appears twice in the PME. While in that case, we only need one auxiliary variable that has the same size as $\hat{A}_{BL}$, it is useful to formally introduce a complete matrix $Z$ and partition it in the same way we partition $\hat{A}$. Doing so, we can treat the auxiliary variable just like all the other operands. Using $Z_{BL}$, we obtain the following decomposition. The corresponding dependency graph is shown in Figure \ref{fig:dg:SPDSys}.
%
%
%
\begin{enumerate}
\item $Z_{BL} := \hat{A}_{BL} \hat{A}_{TL}^{-1}$
\item $A_{BR} := \hat{A}_{BR} - Z_{BL} \hat{A}_{TR}$
\item $b_B := \hat{b}_B - Z_{BL} \hat{b}_T$
\item $x_B := \Sigma \left( A_{BR}, b_B \right)$
\item $b_T := \hat{b}_T - \hat{A}_{TR} x_B$
\item $x_T :=\Sigma \left( \hat{A}_{TL}, b_T \right)$
\end{enumerate}
%
%
%
%
The next step is again to select subsets of the dependency graph to use them as loop invariants and determine the loop guard.

Here, $\hat{A}$ is traversed from the bottom right to the top left, and both $x$ and $\hat{b}$ are traversed from the bottom to the top. Thus, the loop guard $G$ is ``$\text{size}(A_{BR}) < \text{size}(A)$''. If we had solved $\hat{A}_{BL} x_T + \hat{A}_{BR} x_B = \hat{b}_B$ to $x_B$ earlier, as opposed to $x_T$, we would have obtained an algorithm that proceeds in the opposite direction.

For this operation, to select loop invariants, the simplified rule that the subsets can neither be empty nor contain all nodes is not sufficient anymore. We now have to consider the full constraints. In general, it has to be checked whether a loop invariant is feasible, that is, if it leads to an algorithm that actually computes the operation. In \cite{Fabregat-Traver:thesis}, the author lists the following two constraints:
\begin{quote}
\itshape
\begin{enumerate}
\item There must exist a basic initialization of the operands, i.e., an initial partitioning, that renders the predicate $P_\text{\textrm{inv}}$ true:
\begin{align*}
&\{P_\text{\textrm{pre}}\} \\
&\text{\textbf{Partition}} \\
&\{P_\text{\textrm{inv}}\}
\end{align*}
\item $P_\text{\textrm{inv}}$ and the negation of the loop guard, $G$, must imply the postcondition, $P_\text{\textrm{post}}$:
$$P_\text{\textrm{inv}} \land \neg G \Rightarrow P_\text{\textrm{post}}$$
\end{enumerate}
\end{quote}
\label{feasibilityConditions}
The empty subset always fails to satisfy the second constraint, as it translates to an empty predicate. The empty predicate in conjunction with the negation of the loop guard can never imply the postcondition. Similarly, the set that contains all nodes can not satisfy the first condition. In this case, the solution to the operation would already be computed even before the loop is entered. However, merely partitioning the operands does not render such a loop invariant true. In addition to that, predicates $P_\text{before}$ and $P_\text{after}$ would be identical, so no update would be derived.

For SPD linear systems, there are additional subsets that fail to satisfy the second constraint. Let us look at $\{1\}$ for example. The corresponding loop invariant is $P_\text{inv} = \{Z_{BL} := \hat{A}_{BL} \hat{A}_{TL}^{-1}\}$. The negation of the the loop guard ``$\text{size}(A_{BR}) < \text{size}(A)$'' implies that $A_{BR}$ is equal to $A$. Then, $\hat{A}_{BL}$ and $\hat{A}_{TL}$ are empty, and the empty predicate does not imply the postcondition. If we take the subset $\{1, 2\}$, the loop invariant becomes $\{A := \hat{A} \}$ at the end of the operation, which does not imply the postcondition either.

Doing this for all subsets, two feasible loop invariants are obtained:
\begin{enumerate}
\item $\{1, 2, 3, 4\}$ ($P_\text{inv}^1$)
\item $\{1, 2, 3, 4, 5\}$ ($P_\text{inv}^2$)
\end{enumerate}
As mentioned before, the algorithms proceed from the bottom right to the top left, so the repartitionings are different compared to the ones used for the lower triangular system. First, parts of $A_{TL}$ are split off, which are then merged with $A_{BR}$ after the update. Thus, the ``Repartition'' statement for $A$ (and $Z$) is the following:
$$\myFlaTwoByTwoI{A_{TL}}{ A_{TR} }{A_{BL}}{A_{BR}} \rightarrow
	\myFlaThreeByThreeTLI	{A_{00}}	{A_{01}}		{A_{02}}
					{A_{10}}	{A_{11}}	{A_{12}}
					{A_{20}}	{A_{21}}	{A_{22}}$$
This then is the ``Continue with'' partitioning:
$$\myFlaTwoByTwoI{A_{TL}}{ A_{TR} }{A_{BL}}{A_{BR}} 
	\leftarrow
	\myFlaThreeByThreeBRI	{A_{00}}	{A_{01}}		{A_{02}}
						{A_{10}}	{A_{11}}	{A_{12}}
						{A_{20}}	{A_{21}}	{A_{22}}$$
Similarly, the repartitionings for $b$ are:
$$\FlaTwoByOneI{b_T}{b_B} \rightarrow \FlaThreeByOneTI{b_{0}}{b_1}{b_2}$$
$$ \FlaTwoByOneI{b_T}{b_B} \leftarrow \FlaThreeByOneBI{b_{0}}{b_1}{b_2} $$
$x$ is repartitioned in the same way. In the interest of brevity, the derivation of the updates will not be shown here. There are, however, two important points that should be mentioned. When deriving the predicate $P_\text{before}$, the expression
$$\myFlaOneByTwo{Z_{20}}{Z_{21}} := \myFlaOneByTwo{A_{20}}{A_{21}} \myFlaTwoByTwo{A_{00}}{A_{01}}{A_{10}}{A_{11}}^{-1}$$
appears. To flatten it, the PME of a different SPD system, namely $XA = B$, has to be used. Then, in order to find out which expressions appear both in $P_\text{before}$ and $P_\text{after}$ and to identify the differences, it is necessary to rewrite the expressions of both predicates first. Doing so, the auxiliary variables are eliminated.

Finally, the following two updates are found:
\begin{align*}
P_\text{inv}^1: 	b_1 &:= \hat{b}_1 - A_{12} x_2		&	P_\text{inv}^2: 	x_1 &:= \Sigma ( A_{11}, b_1 ) \\
			x_1 &:= \Sigma ( A_{11}, b_1 )	&				b_0 &:= b_0 - A_{01} x_1 
\end{align*}

\section{General Linear System}
\label{sec:generalLS}

In this section, we apply the FLAME methodology to a linear system where $A$ is a general, nonsingular matrix.
$$x:= \Lambda \left( A, b\right) \equiv
\left\{
\begin{aligned}
P_\text{pre}: \{ &\text{\ttfamily Input}[A] \land \text{\ttfamily Matrix}[A] \land \text{\ttfamily NonSingular}[A] \land \\
		&\text{\ttfamily Input}[b] \land \text{\ttfamily Vector}[b] \land \\
		&\text{\ttfamily Output}[x] \land \text{\ttfamily Vector}[x] \} \\
P_\text{post}: \{ & \hat{A}x = \hat{b} \}
\end{aligned}
\right.$$

Since $\hat{A}$ is not lower triangular, there is more than one possibility for the initial partitioning. The different partitionings are shown in Table \ref{tab:LS:part}.
\begin{table}[htp]
\begin{center}
\renewcommand*{\arraystretch}{1.4}
\begin{tabular}{ccc}
\toprule
\#	&	Partitioned postcondition	& Flattened expression	\\ \midrule
1	&	$\myFlaOneByTwo{\hat{A}_L}{\hat{A}_R} \myFlaTwoByOne{x_T}{x_B} = \hat{b}$	&	$\hat{A}_L x_T + \hat{A}_R x_B = \hat{b}$	\\
2	&	$\myFlaTwoByOne{\hat{A}_T}{\hat{A}_B} x = \myFlaTwoByOne{\hat{b}_T}{\hat{b}_B}$	&	$\myFlaTwoByOne{\hat{A}_T x = \hat{b}_T}{\hat{A}_B x = \hat{b}_B}$\\
3	&	$\myFlaTwoByTwo{\hat{A}_{TL}}{ \hat{A}_{TR} }{\hat{A}_{BL}}{\hat{A}_{BR}} \myFlaTwoByOne{x_T}{x_B} = \myFlaTwoByOne{\hat{b}_T}{\hat{b}_B}$	& $\myFlaTwoByOne{\hat{A}_{TL} x_T + \hat{A}_{TR} x_B  = \hat{b}_T}{\hat{A}_{BL} x_T + \hat{A}_{BR} x_B = \hat{b}_B}$ \\ \bottomrule
\end{tabular}
\caption{Possible partitionings for a general linear system.}
\label{tab:LS:part}
\end{center}
\end{table}

At this point, we are again not able to make any further progress with the usual approach. No (sub)expression of the three expressions matches the initial operation. In case of the first two, the parts of $\hat{A}$ are not square. Some or all of the four parts of $\hat{A}$ in the third expression could be square, but since the partitioning does not prescribe any sizes, we can in general not assume that this is the case for any part.

If we use the third partitioning and assume that $\hat{A}_{TL}$ (and $\hat{A}_{BR}$) is square, we are in a similar situation as with the SPD system. Recall that to proceed, we had to rewrite
$$\hat{A}_{TL} x_T + \hat{A}_{TR} x_B  = \hat{b}_T$$
as 
$$x_T = \hat{A}_{TL}^{-1} \left( \hat{b}_T - \hat{A}_{TR} x_B \right) \text{.}$$
This was possible because $\hat{A}$ is SPD, so we were able to infer that $\hat{A}_{TL}$ is SPD, and thus nonsingular, as well. The difference now is that the $\hat{A}_{TL}$ part of a nonsingular matrix $\hat{A}$ is in general not nonsingular. While we could assume that it is, proceed as we did in the previous section and derive similar algorithms, those algorithms would not solve general linear systems.

\section{Equivalence of Loop Invariants}
\label{sec:equivalence}

For direct methods, one observes that finer partitionings lead to additional loop invariants. A finer partitioning here means that we partition along one dimension in more than the usual two parts. Some of those loop invariants lead to truly new algorithms that cannot be derived with a coarser partitioning. Intuitively, one can say that those algorithms expose intermediate steps that are not exposed with a coarser partitioning. Others are in some sense redundant, because they lead to algorithms that are computationally equivalent to algorithms derived from different loop invariants. We then consider those loop invariants to be equivalent. Those loop invariants may or may not have the same granularity. This means that two loop invariants from a $3 \times 3$ partitioning can be equivalent, but it is also possible that a $2 \times 2$ loop invariant is equivalent to a loop invariant obtained from a $3 \times 3$ partitioning.

This is relevant for the derivation of algorithms for iterative methods because we compare different partitionings in Section \ref{sec:MatrixRepresentationPartitionings}. Even more important, the partitioning that is used for iterative methods behaves differently in this regard, as discussed in Section \ref{sec:RemarkEquivalence}.

In this section, we show what equivalence means in this context and under which conditions this behavior occurs. We begin by introducing the notion of equivalence for algorithms and then extend the results to loop invariants.

\subsubsection{Equivalence of Algorithms}

The triangular continuous-time Sylvester equation, which will serve as an example, is defined as $AX+XB=C$ with $A$ and $B$ being triangular. In the following, we assume A and B to be lower and upper triangular, respectively:
%
$$X:= \Psi ( A, B, C) \equiv
\left\{
\begin{aligned}
P_\text{pre}: \{ &\text{\ttfamily Input}[A] \land \text{\ttfamily Matrix}[A] \land \text{\ttfamily LowerTriangular}[A] \land \\
		&\text{\ttfamily Input}[B] \land \text{\ttfamily Matrix}[B] \land \text{\ttfamily UpperTriangular}[B] \land \\
		&\text{\ttfamily Input}[C] \land \text{\ttfamily Matrix}[C] \land \\
		&\text{\ttfamily Output}[X] \land \text{\ttfamily Matrix}[X] \} \\
P_\text{post}: \{ &AX + XB = \hat{C} \}
\end{aligned}
\right.$$
Applying a $1 \times 3$ partitioning to $X$ results in the following PME:
$$
\myFlaOneByThree{
X_L = \Psi(A, B_{TL}, \hat{C}_L)}
{\begin{aligned}
X_M = \Psi(&A, B_{MM}, \\
&\hat{C}_M - X_L B_{TM})
\end{aligned}}
{\begin{aligned}
X_R = \Psi(&A, B_{BL}, \\
&\hat{C}_R - X_L B_{TR} - X_M B_{MR})
\end{aligned}}
$$
At this point, we introduce a new notation for loop invariants taken from \cite{Fabregat-Traver:thesis}. Instead of a logical predicate, the PME is used, leaving out those (sub)expressions that are not part of the loop invariant. If entire blocks are not included, we use the symbol $\neq$ to express that no constraints are imposed on this part. We now choose the following two loop invariants and derive the corresponding updates, skipping the intermediate steps.
\begin{align*}
&P_\text{inv} = \myFlaOneByThree{X_L = \Psi(A, B_{TL}, \hat{C}_L)}	{\neq}	{\neq}\\
&P_\text{inv}' = \myFlaOneByThree{X_L = \Psi(A, B_{TL}, \hat{C}_L)}	{X_M = \Psi(A, B_{MM}, \hat{C}_M - X_L B_{TM})}	{\neq}
\end{align*}
The updates are for $P_\text{inv}$ and $P_\text{inv}'$ are:
\begin{align*}
P_\text{inv}: C_1 &:= \hat{C}_1 - X_0 B_{01}\\
X_1 &:= \Psi(A, B_{11}, C_1) \\
P_\text{inv}': C_2 &:= \hat{C}_2 - X_0 B_{02} - X_1 B_{12} \\
X_2 &:= \Psi(A, B_{22}, C_2)
\end{align*}
Obviously, the updates for $X_1$ and $X_2$ only differ in the indices, which are shifted by one. $B_{22}$ is  the next block following $B_{11}$ on the main diagonal of $B$, and $C_2$ is the next set of columns following $C_1$ in $C$. Considering that those updates happen inside a loop, as the computation unfolds, they will both point to the same parts of $B$ and $C$, respectively.

Written as above, the updates for $C$ do not immediately appear equivalent. Using BLAS operations, $C_1$ would be updated with one call to GEMM, $C_2$ with two. However, when we rewrite the update for $C_2$ as follows, it is easy to see that it is equivalent to just one GEMM operation.
$$C_2 := \hat{C}_2 - \myFlaOneByTwo{X_0}{X_1} \myFlaTwoByOne{B_{02}}{B_{12}}$$
Similarly to $C_1$ and $C_2$, $X_0$ and $\myFlaOneByTwo{X_0}{X_1}$ will, at some point during the computation, contain the same parts of $X$. The same is true for $B_{01}$ and $$\myFlaTwoByOne{B_{02}}{B_{12}}\text{.}$$ Hence, we consider both algorithms to be equivalent. Note that in general, the equivalence can be much less obvious, requiring much more elaborate rewriting.

In the example above, it is easy to see that simply by replacing quantities in one update with the corresponding quantities of the other update, it is possible to transform one update into the other. The only constraint is that if one quantity replaces another, both must have matching sizes. Since we do not specify any concrete dimension and the size of some parts changes during the computation, we have to compare them using a symbolic representation of those sizes. Let us begin with listing the dimensions of all operands of the Sylvester equation:
\begin{itemize}
\item[-] $X$ and $C$ are of size $n \times m$.
\item[-] $A$ is of size $n \times n$.
\item[-] $B$ is of size $m \times m$.
\end{itemize}
After the repartitioning, we specify the sizes as follows.
%
%
%
\renewcommand{\kbldelim}{(}
\renewcommand{\kbrdelim}{)}
$$
\kbordermatrix{
& b i & \omit & b & \omit & b & \omit & b j + c \\
n & X_0 & \omit\vrule & X_1 & \omit\vrule & X_2 & \omit\vrule & X_3
}
$$
$$
\kbordermatrix{
& b i & \omit & b & \omit & b & \omit & b j + c \\
n & C_0 & \omit\vrule & C_1 & \omit\vrule & C_2 & \omit\vrule & C_3
}
$$
$$
\kbordermatrix{
& b i & \omit & b & \omit & b & \omit & b  j + c \\ 
b i & B_{00} & \omit\vrule & B_{01} & \omit\vrule & B_{02} & \omit\vrule & B_{03} \\ \cline{2-8}
b & 0 & \omit\vrule & B_{11} & \omit\vrule & B_{12} & \omit\vrule & B_{13} \\ \cline{2-8}
b & 0 & \omit\vrule & 0 & \omit\vrule & B_{22} & \omit\vrule & B_{23} \\ \cline{2-8}
b j + c & 0 & \omit\vrule & 0 & \omit\vrule & 0 & \omit\vrule & B_{33}
}
$$
$b$ is the block size, $b i$ and $b j$ are unspecified multiples of the block size. $c$ is an additional constant that is nonzero if the block size $b$ does not divide $m$ and/or $n$. Blocks with constant sizes, that is, any combination of $n$, $m$ and $b$, can only be replaced with blocks of the corresponding quantities that have the exact same sizes. For example, it is possible to replace $X_1$ with $X_2$, or $B_{22}$ with $B_{11}$.

Similarly, the dimensions have to match for quantities with variable sizes. If one dimension is a multiple of the block size, any other multiple of the block size matches. As an example, $b i$ matches $b i + b$, and $b j + c + 2b$ matches $b j + c$. It is important to note that $b i$ does not match $b$, because the former is variable, and the latter is constant. Thus, two valid replacements are:
\begin{align*}
B_{00} &\rightarrow \myFlaTwoByTwo{B_{00}}{B_{01}}{0}{B_{11}} \\
\myFlaOneByTwo{C_{2}}{C_{3}} &\rightarrow C_3
\end{align*}
Returning to the initial example, the replacement to transform the update for the loop invariant $P_\text{inv}$ into the one for the loop invariant $P_\text{inv}'$ is the following:
\begin{align*}
B_{01} &\rightarrow \myFlaTwoByOne{B_{02}}{B_{12}} \\
B_{11} &\rightarrow B_{22} \\
C_1 &\rightarrow C_2 \\
X_0 &\rightarrow \myFlaOneByTwo{X_0}{X_1} \\
X_1 &\rightarrow X_2
\end{align*}
%
%

\subsubsection{Equivalence of Loop Invariants}

As mentioned earlier, we consider two loop invariants to be equivalent if the resulting algorithms are equivalent. To determine this equivalence, however, it is not necessary to derive algorithms. The same method of replacing quantities can be performed directly on the loop invariants. If, with such a replacement, one loop invariant can be transformed into another one, both would result in the same update. Consequently, the loop invariants themselves are equivalent.



Let us demonstrate this for the loop invariants $P_\text{inv}$ and $P_\text{inv}'$. We choose the following replacements:
\begin{align*}
X_L &\rightarrow \myFlaOneByTwo{X_{L}}{X_{M}} \\
C_L &\rightarrow \myFlaOneByTwo{C_{L}}{C_{M}} \\
B_{TL} &\rightarrow \myFlaTwoByTwo{B_{TL}}{B_{TM}}{0}{B_{MM}}
\end{align*}
Applying those to the loop invariant
$$P_\text{inv} = \myFlaOneByThree{X_L = \Psi(A, B_{TL}, \hat{C}_L)}	{\neq}	{\neq}\text{,}$$
the equation for computing $X_L$ becomes
$$\myFlaOneByTwo{X_{L}}{X_{M}} = \Psi \left(A,\myFlaTwoByTwo{B_{TL}}{B_{TM}}{0}{B_{MM}}, \myFlaOneByTwo{\hat{C}_{L}}{\hat{C}_{M}} \right) \text{.}$$
Flattening the expressions, we obtain the second loop invariant, $P_\text{inv}'$:
$$P_\text{inv}' = \myFlaOneByThree{X_L = \Psi(A, B_{TL}, \hat{C}_L)}	{X_M = \Psi(A, B_{MM}, \hat{C}_M - X_L B_{TM})}	{\neq}$$
\subsubsection{Equivalence of Loop Invariants of Different Granularities}

The presented method of term rewriting can also be used to decide whether an algorithm or loop invariant is equivalent to a different one obtained with a finer or coarser partitioning. So far, the quantities on both sides of the replacements originated from the same partitioned object. Now, each operand is partitioned twice, using partitionings of two different granularities. Parts obtained from one are replaced with parts of the other. 
Apart from that, exactly the same rules for the replacement hold.

We demonstrate this with a more involved example: The inverse of a lower triangular matrix. A $2 \times 2$ partitioning yields the following PME, using $\tilde{X}$ and $\tilde{L}$ to avoid confusion:
\begin{align}
\myFlaTwoByTwo	{\tilde{X}_{TL} := \tilde{L}_{TL}^{-1}}				{0}
				{\tilde{X}_{BL} := -\tilde{L}_{BR}^{-1} \tilde{L}_{BL} \tilde{X}_{TL}}	{\tilde{X}_{BR} := \tilde{L}_{BR}^{-1}} \label{eq:PME:loTriInv}
\end{align}
This is the $3 \times 3$ PME:
$$
\myFlaThreeByThree{X_{TL} := L_{TL}^{-1}}				{0}					{0}
				{X_{ML} := -L_{MM}^{-1} L_{ML} X_{TL}}	{X_{MM} := L_{MM}^{-1}}	{0}
				{X_{BL} := -L_{BR}^{-1} L_{BM} X_{ML} -L_{BR}^{-1} L_{BL} X_{TL}}		{X_{BM} := -L_{BR}^{-1} L_{BM} X_{MM}}	{X_{BR} := L_{BR}^{-1}}
$$
We will show that the two loop invariants below are equivalent:

\vspace{8pt}
\begin{tabular}{cc}
$\myFlaTwoByTwo	{\tilde{X}_{TL} := \tilde{L}_{TL}^{-1}}				{0}
				{\tilde{X}_{BL} := -\tilde{L}_{BR}^{-1} \tilde{L}_{BL} \tilde{X}_{TL}}	{\neq}$ &
$\myFlaThreeByThree{X_{TL} := L_{TL}^{-1}}				{0}					{0}
				{X_{ML} := -L_{MM}^{-1} L_{ML} X_{TL}}	{\neq}	{0}
				{X_{BL} := -L_{BR}^{-1} L_{BM} X_{ML} -L_{BR}^{-1} L_{BL} X_{TL}}		{\neq}	{\neq}$
\end{tabular}
\vspace{8pt}

%
To do so, we first need to choose an appropriate replacement. The replacement used in this example will transform the $2 \times 2$ loop invariant into the $3 \times 3$ loop invariant, so some parts of the coarser loop invariant will be replaced with multiple parts of the finer one.
\begin{align}
\tilde{X}_{TL} &\rightarrow X_{TL} \label{eq:replacement:XTL}\\
\tilde{L}_{TL} &\rightarrow L_{TL} \label{eq:replacement:LTL}\\
\tilde{X}_{BL} &\rightarrow \myFlaTwoByOne{X_{ML}}{X_{BL}} \label{eq:replacement:XBL}\\
\tilde{L}_{BL} &\rightarrow \myFlaTwoByOne{L_{ML}}{L_{BL}} \label{eq:replacement:LBL}\\
\tilde{L}_{BR} &\rightarrow \myFlaTwoByTwo{L_{MM}}{0}{L_{BM}}{L_{BR}} \label{eq:replacement:LBR}
\end{align}
To make sure that it is a valid replacement, it is necessary to check if the dimensions of the objects on both sides of the replacements match. The sizes of $L$, and likewise $X$, are as follows:
$$
\kbordermatrix{
& b i & \omit & b j + c \\
b i & \tilde{L}_{TL} & \omit\vrule &  0 \\ \cline{2-4}
b j + c & \tilde{L}_{BL} & \omit\vrule & \tilde{L}_{BR}
}
\qquad
\kbordermatrix{
& b i & \omit & b & \omit & b j + c \\
b i & L_{TL} & \omit\vrule & 0 & \omit\vrule & 0 \\ \cline{2-6}
b & L_{ML} & \omit\vrule & L_{MM} & \omit\vrule & 0 \\ \cline{2-6}
b j + c & L_{BL} & \omit\vrule & L_{BM} & \omit\vrule & L_{BR}
}
$$
The replacements (\ref{eq:replacement:XTL}) and (\ref{eq:replacement:LTL}) are easily seen to be correct, as all parts have size $bi \times bi$. In case of (\ref{eq:replacement:XBL}) and (\ref{eq:replacement:LBL}), the left-hand sides have the dimensions $(bj+c) \times b$, and the right-hand sides $(bj+c + b) \times b$. The number of rows matches because both are a multiple of the block size $b$, plus a constant $c$. Similarly, the last replacement (\ref{eq:replacement:LBR}) is valid, as the sizes $(bj+c) \times (bj+c)$ and $(bj+c + b) \times (bj+c + b)$ are equivalent.

Applying this replacement to the coarser loop invariant, we obtain two assignments:
\begin{align}
X_{TL} &:= L_{TL}^{-1} \\
\myFlaTwoByOne{X_{ML}}{X_{BL}} &:= -\myFlaTwoByTwo{L_{MM}}{0}{L_{BM}}{L_{BR}}^{-1} \myFlaTwoByOne{L_{ML}}{L_{BL}} X_{TL} \label{eq:loTriInv1}
\end{align}
Clearly, the first one already has the exact same shape as in the finer loop invariant. The second one, however, requires some rewriting. To be able to compute the product on the right-hand side, we first have to find the symbolic inverse of
$$\myFlaTwoByTwo{L_{MM}}{0}{L_{BM}}{L_{BR}}\text{.}$$
We obtain it from the PME of the inverse of a lower triangular matrix, (\ref{eq:PME:loTriInv}), by eliminating all occurrences of parts of $X$ and replacing the corresponding quantities:
$$\myFlaTwoByTwo{L_{MM}}{0}{L_{BM}}{L_{BR}}^{-1} = \myFlaTwoByTwo{L_{MM}^{-1}}{0}{-L_{BR}^{-1}L_{BM} L_{MM}^{-1}}{L_{BR}^{-1}}$$
Assignment (\ref{eq:loTriInv1}) then becomes
\begin{align}
X_{ML} &:= - L_{MM}^{-1} L_{ML} X_{TL} \label{eq:loTriInv2}\\
X_{BL} &:= L_{BR}^{-1}L_{BM} L_{MM}^{-1} L_{ML} X_{TL} - L_{BR}^{-1} L_{BL} X_{TL} \label{eq:loTriInv3}
\text{.}
\end{align}
The assignment for $X_{ML}$ is now identical to the one in the $3 \times 3$ loop invariant. While the assignments for $X_{BL}$ still differ, we observe that it is possible to replace a subexpression of (\ref{eq:loTriInv3}) with $- X_{ML}$:
$$X_{BL} := L_{BR}^{-1}L_{BM} \underbrace{L_{MM}^{-1} L_{ML} X_{TL}}_{- X_{ML}} - L_{BR}^{-1} L_{BL} X_{TL}$$
Thus, we obtain the following expression, which is the same as the assignment for $X_{BL}$ in the finer loop invariant:
$$X_{BL} := - L_{BR}^{-1}L_{BM} X_{ML} - L_{BR}^{-1} L_{BL} X_{TL}$$

\chapter{Derivation of Algorithms for Iterative Methods: Foundations}
\label{chap:iterativeMethodsIntro}



To be able to use a FLAME-like methodology to derive algorithms for iterative methods, a matrix representation of those methods is required. Such a representation was used in \cite{eijkhout:CGderivation} for CG and the Krylov sequence, and in \cite{eijkhout:CGvariants2} for some CG variants.

This chapter begins with the introduction of an additional notation. The details of the matrix representation are explained in the first section. In the second part, a systematic method for deriving properties of matrices from this representation is presented. Finally, different possible partitionings are discussed, in addition to their implications for the derivation of algorithms.

\subsubsection{Notation}

In this thesis, we use a notation that deviates slightly from the one used in \cite{eijkhout:CGvariants2} and \cite{eijkhout:CGderivation}. We use $e_0$ to denote the unit vector that is one in the first position and $e_r$ for the unit vector that is one in the last position. Both are column vectors. The matrix $J$ is a square matrix with ones on the lower diagonal:
$$
J = \left(
\begin{array}{c@{\;\;}@{\;\;}c@{\;\;}@{\;\;}c@{\;\;}@{\;\;}c} 
0 & 0 & \hdots \\
1 & 0 & \\
0 & 1 & \ddots \\
\vdots & 0 & \ddots
\end{array} 
\right)
$$
 As usual, $I$ is the identity matrix. The dimension of those matrices and vectors will not be given explicitly if they are clear from the context. If $X$ is a matrix, we use $\underline{X}$ to indicate that the right-most column of this matrix is omitted. Thus, $\underline{I}$ and $\underline{J}$ are both lower trapezoidal matrices with one more row than columns.

\section{Matrix Representation}
\label{sec:matrixRepresentationIntroduction}

In this section, we discuss the details of the matrix representation for iterative methods used in this thesis. Note that the representation for CG, which will serve as an example, is a slightly modified version of the one introduced in \cite{eijkhout:CGderivation}. The difference lies in the use of the underline.
%
%
%
%
The three governing equations are shown below. Deriving those equations is not trivial and beyond the scope of this thesis. 
\begin{align}
A P D &= R \left( \underline{I} - \underline{J}  \right) \label{eq:CGrr1} \\
P \left( I - U \right) &= \underline{R} \label{eq:CGrr2} \\
P D &= X \left( \underline{I} - \underline{J} \right) \label{eq:CGrr3}
\end{align}
%
The operands have the following properties:
\begin{itemize}
\item[-] $A \in \mathbb{R}^{n \times n}$ is the coefficient matrix of the linear system that is supposed to be solved. It is nonsingular. Depending on whether $A$ is symmetric or not, different algorithms can be derived.
\item[-] $P \in \mathbb{R}^{n \times m}$ is the matrix of search directions, that is, each column represents the search direction vector during one iteration. It is initially unknown.
\item[-] $D \in \mathbb{R}^{m \times m}$ is an unknown diagonal matrix.
\item[-] $R \in \mathbb{R}^{n \times (m + 1)}$ is the residual matrix. Initially, only the first column $r_0$ is known. It is computed as $r_0 = A x_0 - b$, where $x_0$ is an initial guess for the solution. Additionally, it is orthogonal.
\item[-] $U \in \mathbb{R}^{m \times m}$ is unknown and upper diagonal\footnote{This property is defined in Appendix \ref{chap:appendixProperties}.} if $A$ is symmetric. Otherwise, it is strictly upper triangular.
\item[-] $X \in \mathbb{R}^{n \times (m + 1)}$ is the matrix of approximated solution vectors. Similar to $R$, only the first column is initially known.
\end{itemize}

While it might seem unusual that the same matrix ($R$) appears twice in the governing equations with varying sizes, this is necessary to ensure the correctness of the last column of $R$ and $X$.  The formula for computing the residual, in indexed notation, is $r_{i+1} = r_i - A p_i \delta_i$. Without the additional column of $R$ in equation (\ref{eq:CGrr1}), the incorrect equation $A p_i \delta_i = r_i$ would be obtained for the last iteration. Similarly, from equation (\ref{eq:CGrr3}), we would obtain $p_i \delta_i = x_i$, which is not correct either.

One of the fundamental differences compared to direct methods is that the dimensions of some matrices are not fixed. 
Usually, all dimensions are determined by the sizes of input operands. In case of the LU-factorization, for example, we know that $L$ and $U$ have the same size as $A$. If not, the equation $LU = A$ is not valid. If $A$ is of size $n \times n$, and the $LU$ factorization of a $k \times k$ block of $A$, with $k < n$, is computed, then the postcondition is not rendered true.

Due to the orthogonality of $R$, $m + 1$ can not be larger than $n$, as the number of $n$-dimensional, orthogonal vectors is at most $n$. However, for every $m < n$, the equations above can be satisfied by performing the corresponding number of iterations.

For the systematic derivation of algorithms, this introduces the problem that it is not possible to derive a loop guard exclusively from the equations. It makes no sense to compare the size of a block of $R$ to $R$ itself, because $R$ grows too. For those iterative methods that are used to find solutions for linear systems, the goal is usually to minimize the residual in some norm \cite{barrett:templates}. 
Thus, the loop guard typically is a predicate like ``$\| r_i \| \geq \varepsilon$'', where $\varepsilon$ is a threshold chosen by the user. We make sure that the postcondition of CG correctly represents the situation at the end of the operation by adding $\| R e_r^T \| < \varepsilon$ to it. This also allows us to derive a suitable loop guard from it.

For stationary iterative methods, the stopping criterion quite often is $\| x_i - x_{i-1} \| < \varepsilon$. Translating that in our notation, we obtain $\| X e_r^T - X e_{r-1}^T \| < \varepsilon$. 
%
%
For iterative methods where those criteria are not applicable, we will add an expression to the postcondition that fixes the number of columns of a matrix with variable size. Such a predicate could be ``$\text{size}(Y) = n \times k$''.\footnote{It is not uncommon to combine multiple criteria \cite{barrett:templates}. For the sake of simplicity, we do not use more than one at a time.}

\section{Systematic Derivation of Matrix Properties}
\label{sec:propertyDerivation}

In order to derive algorithms for CG from its description in matrix form, it is necessary to use properties of matrices and expressions that are not explicitly part of the initial description. Those properties have to be derived from the description by means of algebraic manipulation and deductive reasoning. To automate the process of finding algorithms, an automatic method for the derivation of properties is necessary. Thus, this systematic approach should replicate the steps performed by a human expert, without actually requiring any human guidance. In this section, we describe such a method.



\subsection{Preliminaries}

To begin, it is useful to formalize the notion of properties and equations. We start with the most basic building blocks, terms and expressions.

\begin{definition}[Terms and Expressions]
\emph{Terms} are inductively defined as follows:
\begin{enumerate}
\item Every matrix, vector and scalar is a term.
\item If $t$ is a term, then $t^T$, $(-t)$ and $t^{-1}$ are terms.
\item If $t_1$ and $t_2$ are terms, then $(t_1 + t_2)$ and $(t_1 \cdot t_2)$ are terms.
\end{enumerate}
Every term is also an \emph{expression}. Furthermore, if $t_1$ and $t_2$ are terms, then $t_1 = t_2$ is an \emph{expression}. \qed
\end{definition}

Note that this is a simplified definition, as it does not include functions,
nor does it make any statement about the validity of terms. A product of two matrices with dimensions that do not match is still a valid term. Nonetheless, this definition is sufficient for our purposes. Furthermore, we usually use a simplified notation, omitting parentheses, if they are unnecessary, as well as the multiplication dot. We can now define equations and properties:

\begin{definition}[Equations]
Let $t_1$, $t_2$ be terms. An \emph{equation} is an expression of the form $t_1 = t_2$. \qed
\end{definition}

\begin{definition}[Properties]
Let $t$ be a term and $\text{\ttfamily P}$ be a boolean predicate. Then $\text{\ttfamily P}[t]$ is a \emph{property}. \qed
\end{definition}

Note that properties are always predicates, even if they can also be expressed as equations. Take the property $\text{\ttfamily Orthogonal}[R]$ as an example. It implies $\text{\ttfamily Diagonal}[R^T R]$. Applying the partitioning that is also applied to the postcondition,
$$R \rightarrow \myFlaOneByThree{R_L}{r_M}{R_R}\text{,}$$
to $R^T R$, we find out, among other things, that $R_L^T R_L$ is diagonal and $r_M^T R_L$ equals zero. 
Instead of considering the equation $r_M^T R_L = 0$ as a property, we use $\text{\ttfamily Zero}[r_M^T R_L]$. This allows us to use a more consistent notation and simplifies the systematic derivation.

%

\subsection{Representing Knowledge about Linear Algebra}

A human expert who derives properties of matrices inevitably applies some basic knowledge about linear algebra. To allow a system to replicate the expert reasoning, it needs a knowledgebase that encodes this knowledge. We define five different types of implications that will be included in the knowledgebase.
\begin{enumerate}
\item $\text{\ttfamily P}_1[t] \land \ldots \land \text{\ttfamily P}_i[t] \rightarrow \text{\ttfamily P}[t]$

This type of implication allows to reason about the combination of properties of one single term. One example is:
\begin{align}
\text{\ttfamily LowerTriangular}[t] \land \text{\ttfamily Symmetric}[t] &\rightarrow \text{\ttfamily Diagonal}[t] 
\end{align}
\item $\text{\ttfamily P}_1[t_1] \land \ldots \land \text{\ttfamily P}_i[t_i] \land \exists t \rightarrow \text{\ttfamily P}[t] $

Here, $t_1,\ldots, t_i$ are subterms of $t$. Thus, it allows the system to infer the properties of a product or sum of multiple quantities with different properties. The $\exists t$ is used to avoid deriving properties for terms that do not occur anywhere. Consider two examples:
\begin{align}
\text{\ttfamily Diagonal}[t_1] \land \text{\ttfamily LowerTriangular}[t_2] \land \exists t_1 t_2 &\rightarrow \text{\ttfamily LowerTriangular}[t_1 t_2] \\
\text{\ttfamily StrictlyUpperTriangular}[t] \land \exists I + t &\rightarrow \text{\ttfamily UpperTriangular}[I + t]
\end{align}
\item $\text{\ttfamily P}[t_1] \land t_1 = t_2 \rightarrow \text{\ttfamily P}[t_2]$

This implication enables the system to propagate properties across equalities. $\text{\ttfamily P}$ now is a pattern that matches any property.
\item $\text{\ttfamily P}_1[t] \rightarrow \text{\ttfamily P}_2[f(t)]$

$f$ is a function, for example transposition. Thus, this kind of implication allows to derive properties of transposed or inverted quantities. In addition to that, it is used reason about orthogonal or orthonormal matrices. Let us look at three examples:
\begin{align}
\text{\ttfamily LowerTriangular} \left[t\right] &\rightarrow \text{\ttfamily UpperTriangular} \left[ t^T \right] \label{eqn:loTriTransposed}\\
\text{\ttfamily LowerTriangular} \left[t\right] &\rightarrow \text{\ttfamily LowerTriangular} \left[ t^{-1} \right] \\
\text{\ttfamily Orthogonal}[t] &\rightarrow \text{\ttfamily Diagonal}\left[t^T t \right]
\end{align}
\item $t \rightarrow f(t)$

For one of the steps of the method presented in this chapter, it is necessary that properties have a canonical form: The unary operators $^{-1}$ and $^{T}$ will not be applied to products. However, some of the other types of implications may produce such terms. To transform those terms into the canonical form, implications are necessary that distribute those unary operators across products. Those implications are:
%
\begin{align}
(t_1 t_2)^T &\rightarrow t_2^T t_1^T \label{eqn:productTransposed} \\
(t_1 t_2)^{-1} &\rightarrow t_2^{-1} t_1^{-1} 
\end{align}
\end{enumerate}
%
%
Note that all the terms in the implications above may match not only a single op\-er\-and, but every term that has the required property. So if the product $AB$ is lower triangular, implication (\ref{eqn:loTriTransposed}) tells us that $(AB)^T$ is upper triangular, and using (\ref{eqn:productTransposed}), we find out that $B^T A^T$ is upper triangular.


\subsection{Derivation of Properties}
\label{sec:propertyDerivation:Derivation}

In the following, we will assume that all equations have the form
$$\prod_{1 \leq i \leq n} t_i =\prod_{1 \leq j \leq m} t^\prime_j \text{.}$$
%
In case of a sum $X + Y$, the product consists of just one term $t = X + Y$. This assumption is no restriction as the only sums that appear in the matrix representations of iterative methods are of the form $I - Y$ and are used exclusively to emphasize the structure of those matrices.

\begin{description}
\item[Initialization] We start with two sets, $\mathcal{P}$ and $\mathcal{E}$. $\mathcal{P}$ is a subset of the precondition $P_\text{pre}$ of the description of an operation. It only contains those properties that describe operand types, like $\text{\ttfamily Square}[X]$, $\text{\ttfamily Matrix}[X]$ or $\text{\ttfamily Diagonal}(X)$. Not included are properties that specify what is input and output, as they are superfluous for the derivation. $\mathcal{E}$ is a set of expressions which contains the equations of the corresponding postcondition $P_\text{post}$.

\item[Derivation of Properties] The implications of the knowledgebase $\mathcal{K}$ are used to derive all possible properties at this stage. This is done as follows:

\begin{itemize}
\item[-] Given $\text{\ttfamily P}_1[t] \land \ldots \land \text{\ttfamily P}_i[t] \rightarrow \text{\ttfamily P}[t] \in \mathcal{K}$, the property $\text{\ttfamily P}[t]$ is added to $\mathcal{P}$ if $\text{\ttfamily P}_k[t] \in \mathcal{P}$ for all $k \in \{1,\ldots, i\}$.

\item[-] Given $\text{\ttfamily P}_1[t] \land \ldots \land \text{\ttfamily P}_i[t] \land \exists t \rightarrow \text{\ttfamily P}[t] \in \mathcal{K}$, the property $\text{\ttfamily P}[t]$ is added to $\mathcal{P}$ if 
\begin{enumerate}
\item there is an equation $e \in \mathcal{E}$ that contains the term $t$, and
\item $\text{\ttfamily P}_k[t] \in \mathcal{P}$ for all $k \in \{1,\ldots, i\}$.
\end{enumerate}

\item[-] Given $\text{\ttfamily P}[t_1] \land t_1 = t_2 \rightarrow \text{\ttfamily P}[t_2] \in \mathcal{K}$, the property $\text{\ttfamily P}[t_2]$ is added to $\mathcal{P}$ if $\text{\ttfamily P}[t_1] \in \mathcal{P}$ and $(t_1 = t_2) \in \mathcal{E}$.

\item[-] Given $\text{\ttfamily P}_1[t] \rightarrow \text{\ttfamily P}_2[f(t)] \in K$, the property $\text{\ttfamily P}_2[f(t)]$ is added to $\mathcal{P}$ if $\text{\ttfamily P}_1[t] \in \mathcal{P}$.

\item[-] Given $t \rightarrow f(t) \in \mathcal{K}$ and $\text{\ttfamily P}[t] \in \mathcal{P}$, the property $\text{\ttfamily P}[f(t)]$ is added to $\mathcal{P}$.

\end{itemize}

\item[Matrix Inversion] Then, new expressions are added to $\mathcal{E}$.
\begin{enumerate}
\item For every equation $e \in \mathcal{E}$ with $e = ( t_1 \cdots t_i = u_1 \cdots u_j )$, it is checked whether $t_1$, $t_i$, $u_1$ or $u_j$ is invertible. 
\item If $t_1$ or $u_1$ is invertible, $t_1^{-1}$ or $u_1^{-1}$, respectively, is multiplied from the left to $e$, eliminating the invertible term on the one side and adding its inverse on the other. We proceed analogously with $t_i$ and $u_j$, then multiplying the inverse of the term from the right.
\item The resulting, new equation $e'$ is added to $\mathcal{E}$.
\end{enumerate}

Those three steps are repeated until no new expressions are found.


\item[Application of Properties] In the final step, we apply known properties to expressions to derive new properties. Intuitively, we multiply quantities to both sides of an equation in order to recreate subexpressions that are known to present some property, which are then used to infer new properties.
\begin{enumerate}
\item For every property $\text{\ttfamily P}[t] \in \mathcal{P}$ with $t = t_1 \cdots t_i$, $i > 1$ a set of tuples $$S(t) = \{ (t_1 \cdots t_k, t_{k + 1} \cdots t_i) \mid 1 \leq k < i \} $$ is generated. Intuitively, this set contains the term $t$, split into two parts in all possible ways.
\item Let $(t_L, t_R) \in S(t)$. If there is an equation $e \in \mathcal{E}$ where $t_L$ is the rightmost (sub)term in any term, then $t_R$ is multiplied from the right. $t_L$ is multiplied from the left in the analogous case. Let $e'$ be the resulting equation. As an example, let $e$ be $ABC = D$ and $(BC, F) \in S(t)$. Since $BC$ appears as the rightmost subterm in $e$, $t_R = F$ is multiplied from to right, resulting in $ABCF = DF$.
\item The knowledgebase is used to derive new properties using $e'$ that are added to $\mathcal{P}$ (see step ``Derivation of Properties''). $e'$ is \emph{not} added to $\mathcal{E}$.
\end{enumerate}

\end{description}

\subsection{Design Considerations}

One observes that the presented method is not goal-oriented. The derived properties are mainly used to solve equations, so it might seem more natural to derive properties starting with an equation that has to be solved. Based on this equation, expressions would be selected, and in a second step, the properties of those expression would be derived. Those properties would then be used to solve the equation. Instead, the presented method derives a large number or properties, irrespective of the question whether they might be useful or not.

The problem with a goal-oriented approach is that quite often, it is not obvious from the equation which property could be used to solve it. This leaves us with the much more challenging task of identifying which expression might have relevant properties. Take the following equation as an example, which appears when deriving algorithms for nonsymmetric CG. $P_L$, $p_M$, $r_M$ and $A$ are known.
$$-P_L u_{TM} + p_M = r_M$$
It is solved for $u_{TM}$ by using the fact that $P^T A P$ is lower triangular. Thus, $P_L^T A P_L$ is lower triangular as well and $P_L^T A p_M$ is zero. Multiplying $P_L^T A$ from the left to both sides of the equation gives us
$$-P_L^T A P_L u_{TM} = P_L^T A r_M \text{.}$$
This now is a triangular system that can easily be solved. While it is possible to individually derive that $P_L^T A P_L$ is lower triangular and $P_L^T A p_M$ is zero, the initial equation gives us little to no indication to inspect the properties of those expressions in the first place.

One the other hand, the advantage of a method that is not goal-oriented is that it may find properties that we do not expect to find.

\subsubsection{Orthogonality of the Residual Matrix}

For some iterative methods, for example CG, the residual matrix $R$ is orthogonal. In the postcondition of those methods, $R$ usually appears multiple times, either as a whole, or without the last column ($\underline{R}$). Unfortunately, for the derivation of properties, this poses a problem.

If $R$ is orthogonal, then $R^T R$ and $\underline{R}^T \underline{R}$ are diagonal. $R^T \underline{R}$ and its transpose are rectangular and all entries except for the ones on the main diagonal are zero. Thus, in some sense, they are diagonal as well.

It turns out that for deriving certain properties, $R^T \underline{R}$ is needed. Unfortunately, using the described method, neither $\text{\ttfamily Orthogonal}[R]$ nor $\text{\ttfamily Orthogonal}[\underline{R}]$ implies any property of $R^T \underline{R}$. One way to solve this would be to treat any matrix that appears with and without the last column in a special way, such that properties of $R^T \underline{R}$ and $\underline{R}^T R$ are found as well. This, however, would require significant modifications of the derivation process.

The simpler solution is to consider $R$ and $\underline{R}$ to be two distinct objects, and properties of $R^T \underline{R}$ and its transpose are added to the precondition. 

\subsubsection{Substituting Equations}

Note that we deliberately avoid substituting quantities in one equation by expressions obtained from others. While this might be a very natural approach if deriving properties by hand, doing this systematically is difficult. Plugging in equations quickly leads to arbitrarily large expressions unless some heuristics are applied to terminate this process. Apart from that, it is possible to achieve the same results using the approach presented in this section. Consider a short example to get an intuition why this is the case. Let us assume we have two equations
\begin{align*}
t_1 &= t_2 {\color{green!50!black} t_3} \\
{\color{green!50!black}t_4} t_2 &= {\color{green!50!black}t_5}
\end{align*}
and we want to derive a property for $t_1$. Properties of $t_3$, $t_4$ and $t_5$ are known (colored green). As the properties of $t_2$ are not known as well, we have to use the second equation to proceed. If $t_4$ is invertible, we can solve to
$$t_2 = {\color{green!50!black}t_4^{-1}} {\color{green!50!black}t_5} \text{.}$$
Instead of substituting $t_2$ in $t_1 = t_2 t_3$, yielding
$$t_1 = {\color{green!50!black}t_4^{-1}} {\color{green!50!black}t_5} {\color{green!50!black} t_3} \text{,}$$
and then reason about properties of $t_4^{-1} t_5$, we first derive all properties of $t_4^{-1} t_5$, which are also properties of $t_2$. In the final step, we derive all properties of $t_2 t_3$, obtaining the same properties for $t_1$ we would find by plugging one equation into the other.

\subsubsection{Termination}

The disadvantage of this approach is that it may not terminate either, and increasingly long properties are derived. This is not unexpected, since this approach aims at replicating the process of substituting equations. The difference of the presented method is that the set of equations is finite, its size does not even change anymore after the matrix inversion step.\footnote{This is why we refrain from adding $e'$ to $\mathcal{E}$ in ``Application of Properties'', step 3.} Furthermore, for most iterative methods, no properties of products of  more than three quantities are used. While it might be possible to use significantly longer properties to solve equations, most likely, they result in algorithms that use unnecessarily large expressions to compute certain quantities. This naturally leads to the solution of introducing a maximum length for properties, similar to a recursion limit, with a reasonable default value that can be changed by the user. This way, there is only a finite number of properties that can be derived. 

Unfortunately, if we just refrained from adding properties of products of more than three quantities to $\mathcal{P}$, the derivation process would cease to work in certain cases (this will be explained in the following section). To avoid that, there are multiple options. In both cases, we initially add longer properties to $\mathcal{P}$ as well. Then, one solution is to never use them to derive further properties, that is, we never construct $S(t)$ if $t$ is a product of more than three quantities. Alternatively, those longer properties are removed from $\mathcal{P}$ when $e'$ is discarded. A third option would be to construct an additional set of temporary properties.

\subsection{Example: Nonsymmetric CG}
\label{sec:derivationOfPropertiesExample}

We demonstrate the method presented in this chapter by deriving some properties for nonsymmetric CG. In the interest of brevity, we only derive a small number of selected properties, in addition to limiting properties to a maximum length of three quantities. The knowledgebase $\mathcal{K}$ is not shown here due to its size. Properties used in the following which are not self-explanatory are defined in Appendix \ref{chap:appendixProperties}. The pre- and postcondition of nonsymmetric CG are shown below. For the sake of simplicity, we treat $\underline{I} - \underline{J}$ as one distinct matrix, as $\underline{I}$ and $\underline{J}$ do not appear separately.
\begin{align*}
P_\text{pre}: \{ &\text{\ttfamily Input}[A] \land \text{\ttfamily Matrix}[A] \land \text{\ttfamily NonSingular}[A] \land \\
		&\text{\ttfamily Output}[P] \land \text{\ttfamily Matrix}[P] \land \\
		&\text{\ttfamily Output}[D] \land \text{\ttfamily Matrix}[D] \land \text{\ttfamily Diagonal}[D] \land \\
		&\text{\ttfamily FirstColumnInput}[R] \land \text{\ttfamily Matrix}[R] \land \text{\ttfamily Orthogonal}[R] \land \\
		&\text{\ttfamily FirstColumnInput}[\underline{R}] \land \text{\ttfamily Matrix}[\underline{R}] \land \text{\ttfamily Orthogonal}[\underline{R}] \land \\
		&\text{\ttfamily DiagonalR}[R^T \underline{R}] \land \text{\ttfamily DiagonalR}[\underline{R}^T R] \land \\
		&\text{\ttfamily FirstColumnInput}[X] \land \text{\ttfamily Matrix}[X] \land \\
		&\text{\ttfamily Output}[U] \land \text{\ttfamily Matrix}[U] \land \text{\ttfamily StrictlyUpperTriangular}[U] \land \\
		& \text{\ttfamily Matrix}[\underline{I} - \underline{J}] \land \text{\ttfamily LowerTrapezoidal}[\underline{I} - \underline{J}] \}
\end{align*}
\begin{align*}
P_\text{post}: \{ &APD = R \left( \underline{I} - \underline{J}  \right) \\
			&P \left( I - U \right) = \underline{R} \\
			&PD = X \left( \underline{I} - \underline{J} \right) \\
			&\| R e_r^T \| < \varepsilon\}
\end{align*}
\subsubsection{Initialization}

The first step consists of initializing $\mathcal{P}$ and $\mathcal{E}$. The former contains all the properties of the precondition $P_\text{pre}$ that describe operand types.
\begin{align*}
\mathcal{P} = \{ & \text{\ttfamily Matrix}[A], \text{\ttfamily NonSingular}[A], \\
		& \text{\ttfamily Matrix}[P], \\
		& \text{\ttfamily Matrix}[D], \text{\ttfamily Diagonal}[D], \\
		& \text{\ttfamily Matrix}[R],  \text{\ttfamily Orthogonal}[R], \\
		& \text{\ttfamily Matrix}[\underline{R}], \text{\ttfamily Orthogonal}[\underline{R}], \\
		&\text{\ttfamily DiagonalR}[R^T \underline{R}], \text{\ttfamily DiagonalR}[\underline{R}^T R] \land \\
		& \text{\ttfamily Matrix}[X], \\
		& \text{\ttfamily Matrix}[U], \text{\ttfamily StrictlyUpperTriangular}[U], \\
		& \text{\ttfamily Matrix}[\underline{I} - \underline{J}], \text{\ttfamily LowerTrapezoidal}[\underline{I} - \underline{J}] \} \\
\end{align*}
$\mathcal{E}$ contains the equations of the postcondition $P_\text{post}$.
\begin{align*}
\mathcal{E} = \{ &APD = R \left( \underline{I} - \underline{J}  \right), \\
	&P \left( I - U \right) = \underline{R}, \\
	&PD = X \left( \underline{I} - \underline{J} \right) \}
\end{align*}
\subsubsection{Derivation of Properties}
During this step, only a small number of new properties can be derived. The implication
$$\text{\ttfamily StrictlyUpperTriangular}[t_1] \land I + t_1 \rightarrow \text{\ttfamily UpperTriangular}[I + t_1]$$
is used to infer that $(I - U)$ is upper triangular. $\text{\ttfamily Orthogonal}[R]$ and $\text{\ttfamily Orthogonal}[\underline{R}]$ imply that $R^T R$ and $\underline{R}^T \underline{R}$ are diagonal. Thus, this step yields
\begin{align*}
\mathcal{P} = P \cup \{ &\text{\ttfamily Diagonal}[R^T R], \\
				&\text{\ttfamily Diagonal}[\underline{R}^T \underline{R}], \\
				& \text{\ttfamily UpperTriangular}[I - U] \}\text{.}
\end{align*}
\subsubsection{Matrix Inversion}

From the set of properties $\mathcal{P}$ it follows that $A$, $D$ and $(I - U)$ are nonsingular. The equation $APD = R \left( \underline{I} - \underline{J} \right)$ is inspected first. Two invertible objects, namely $A$ and $D$, occur in it. Since $A$ is the leftmost quantity of the product $APD$, it can be eliminated by multiplying its inverse from the left-hand side to both sides of the equation. This yields $PD = A^{-1} R \left( \underline{I} - \underline{J} \right)$, which is added to $\mathcal{E}$. By multiplying $D^{-1}$ from the right to this new equation and the original $APD = R \left( \underline{I} - \underline{J} \right)$, two additional equations are obtained. By applying the same procedure to the two remaining equations, $\mathcal{E}$ becomes the following set:
\begin{align*}
\mathcal{E} = \{ &APD = R \left( \underline{I} - \underline{J} \right), &&AP = R \left( \underline{I} - \underline{J} \right)D^{-1}, \\
	&PD = A^{-1}R \left( \underline{I} - \underline{J} \right), &&P = A^{-1} R \left( \underline{I} - \underline{J} \right) D^{-1}, \\
	&P \left( I - U \right) = \underline{R}, &&P = \underline{R} \left( I - U \right)^{-1}, \\
	&PD = X \left( \underline{I} - \underline{J} \right), &&P = X \left( \underline{I} - \underline{J} \right) D^{-1} \}
\end{align*}
At this point, it is not possible to find any new equations by multiplying inverted quantities, so we proceed to the next step.

\subsubsection{Application of Properties}

Initially, the only properties $\text{\ttfamily P}[t] \in \mathcal{P}$ with $t = t_1 \cdots t_i$, $i > 1$ are
\begin{align*}
&\text{\ttfamily Diagonal}[R^T R] &&\text{\ttfamily Diagonal}[\underline{R}^T \underline{R}]\\
&\text{\ttfamily DiagonalR}[R^T \underline{R}] && \text{\ttfamily DiagonalR}[\underline{R}^T R] \text{.}
\end{align*}
From $\text{\ttfamily Diagonal}[R^T R]$, the set $S(R^T R) = \{(R^T, R)\}$ is obtained. For every $(t_L, t_R) \in S(R^T R)$, it now has to be checked if $t_L$ or $t_R$ appears in any equation contained in $E$. Since $(R^T, R)$ is the only element in $S(R^T R)$, the system just searches for $t_L = R^T$ and $t_R = R$. The following four equations, all containing $R$, are found:
\begin{align*}
&APD = R \left( \underline{I} - \underline{J} \right) &&AP = R \left( \underline{I} - \underline{J} \right)D^{-1} \\
	&PD = A^{-1}R \left( \underline{I} - \underline{J} \right) &&P = A^{-1} R \left( \underline{I} - \underline{J} \right) D^{-1}
\end{align*}
Only in the first two equations, $t_L = R$ appears at the rightmost position in a product. Hence, $t_R = R^T$ is multiplied just to those two:
\begin{align*}
&R^TAPD = R^TR \left( \underline{I} - \underline{J}  \right) &&R^TAP = R^TR \left( \underline{I} - \underline{J}  \right)D^{-1}
\end{align*}
The knowledgebase is now used to derive new properties. In this example, we just look at the equation on the right-hand side, and show the steps that the system would perform.
\begin{enumerate}
\item $R^TR \left( \underline{I} - \underline{J}  \right)$ is a product of a diagonal and a lower trapezoidal matrix, so it is lower trapezoidal.
\item $D$ is diagonal, so $D^{-1}$ is diagonal as well.
\item $R^TR \left( \underline{I} - \underline{J}  \right)D^{-1}$ is a product of a lower trapezoidal matrix ($R^TR \left( \underline{I} - \underline{J}  \right)$) and a diagonal matrix ($D^{-1}$), so it is lower trapezoidal too.
\item The right-hand side of the equation is lower trapezoidal, so the left-hand side of the equation, $R^TAP$, is lower trapezoidal as well.
\end{enumerate}
With every step, the new properties are added to $\mathcal{P}$. Thus, the set becomes:
\begin{align*}
\mathcal{P} := \mathcal{P} \cup \{ & \text{\ttfamily LowerTrapezoidal}[R^TR \left( \underline{I} - \underline{J}  \right) ], \\
	& \text{\ttfamily Diagonal}[D^{-1} ], \\
	& \text{\ttfamily LowerTrapezoidal}[R^TR \left( \underline{I} - \underline{J}  \right)D^{-1} ], \\
	& \text{\ttfamily LowerTrapezoidal}[R^T A P ] \}
\end{align*}
Now, it becomes apparent why it is not possible to simply set a limit on the length of properties and never add any longer properties to $\mathcal{P}$. $R^TR \left( \underline{I} - \underline{J}  \right)D^{-1}$ is a product of four quantities. To infer that $R^T A P$ is lower trapezoidal, the property $$\text{\ttfamily LowerTrapezoidal}[R^TR \left( \underline{I} - \underline{J}  \right)D^{-1}]$$ has to be derived. If this property is never added to $\mathcal{P}$, it is not possible to derive that $R^T A P$ is lower trapezoidal. Hence, it must be possible to derive properties of any length, even though they are only needed temporarily.

Using the newly added properties, immediately some more are found because of the implications 
\begin{align*}
\text{\ttfamily LowerTrapezoidal} \left[t\right] &\rightarrow \text{\ttfamily UpperTrapezoidal} \left[ t^T \right] \\
(t_1 t_2)^T &\rightarrow t_2^T t_1^T \text{.}
\end{align*}
The following properties are added to the set:
\begin{align*}
\mathcal{P} := \mathcal{P} \cup \{ & \text{\ttfamily UpperTrapezoidal}[ \left( \underline{I} - \underline{J}  \right)^T R^TR ], \\
	& \text{\ttfamily Diagonal} [D^{-T} ], \\
	& \text{\ttfamily UpperTrapezoidal} [ D^{-T} \left( \underline{I} - \underline{J} \right) R^TR ], \\
	& \text{\ttfamily UpperTrapezoidal} [P^T A R ] \}
\end{align*}
At this point, we cut this derivation short. In practice, many more properties would be derived. Most of them may not be of any use for the derivation of algorithms, but a few are crucial. For nonsymmetric CG, such an important property is that $P^T A P$ is lower triangular. In \cite{eijkhout:CGderivation}, this property is derived manually. Let us shortly illustrate how it is derived.

\begin{enumerate}
\item Because of $\text{\ttfamily DiagonalR}[\underline{R}^T R]$, $R^T$ is multiplied from the left to $P = \underline{R} \left( I - U \right)^{-1}$, yielding
$$R^T P = R^T \underline{R} \left( I - U \right)^{-1}\text{.}$$
$\left( I - U \right)^{-1}$ is upper triangular and $R^T \underline{R}$ is rectangular diagonal, so $R^T P$ is upper triangular and rectangular. It follows that its transpose $P^T R$ is lower triangular and rectangular.
\item Because of the property $\text{\ttfamily LowerTriangularR}[P^T R]$, $P^T$ is multiplied from the left to $AP D = R \left( \underline{I} - \underline{J}  \right)$. Based on the resulting equation
$$P^T A P D = P^T R \left( \underline{I} - \underline{J}  \right)\text{,}$$
we derive that $P^T R \left( \underline{I} - \underline{J}  \right)$ is square.
\item For the same reason, $P^T$ is multiplied from the left to $AP = R \left( \underline{I} - \underline{J}  \right)D^{-1}$, resulting in
$$P^T A P = P^T R \left( \underline{I} - \underline{J}  \right)D^{-1}\text{.}$$
$P^T R$ is lower triangular and rectangular and $\underline{I} - \underline{J}$ is lower trapezoidal. $P^T R \left( \underline{I} - \underline{J}  \right)$ is square, so it is lower triangular. Since $D^{-1}$ is diagonal, $P^T A P$ is lower triangular as well. 
\end{enumerate}

\section{Initial Partitionings}
\label{sec:MatrixRepresentationPartitionings}

The first step towards the derivation of algorithms in the FLAME methodology is to partition the operands. The matrix representation of iterative methods gives rise to a new type of operands, namely those which are initially partially known and partially unknown. Thus, we have to address the question of how to partition them. To some extent, it is easy to answer. Standard CG algorithms always compute full vectors $r_i$ and $x_i$ \cite{barrett:templates, saad2000iterative, vanderVorst:book}, so refraining from partitioning $R$ and $X$ horizontally, that is, into a top and a bottom part, is a very natural choice. For partitioning vertically, there are multiple possibilities. Their advantages and disadvantages will be discussed in this section.

While not strictly necessary to derive algorithms, describing the partitioned op\-er\-ands in terms of functions proved to be useful for the systematic derivation \cite{Fabregat-Traver2011:54}. Hence, if and how different partitionings permit to match functions is an important criterion for their evaluation. 


We begin this section with a discussion of what such a function should look like. Then, we investigate different possible partitionings, starting with partitionings that are also used for direct methods and continuing with some that are tailored to properties of algorithms for iterative methods.

\subsubsection{Functions}

On the highest level, the input of a CG algorithm are vectors $b$ and $x_0$, as well as the matrix $A$, so a function would have the form $x_i :=  \text{CG}(A, b, x_0)$. As we are using the recurrence relation to derive algorithms, which does not involve $b$, but a number of other quantities, for example the residual $r$, a function like this is not helpful.

From direct methods, we remember that the function is  uniquely\footnote{Except for the ordering of input and output.} defined by the precondition. Every quantity that is initially known has the property $\text{\ttfamily Input}[X]$, initially unknown quantities have the property $\text{\ttfamily Output}[X]$. Take the LU factorization as an example. The governing equation is $LU = A$. $A$ is input, $L$ and $U$ are output. This naturally leads to a function $\{L, U\} = \Psi (A)$. It is important to note that this function is defined prior to any derivation steps, solely based on the abstract description of the operation. It then happens, due to the recursive nature of the operations, \emph{and a partitioning that reveals it}, that this function matches expressions obtained by partitioning the postcondition. In some cases, the expressions have to be rewritten first. The remaining expressions can be decomposed into basic buildings blocks.

Applying this scheme to CG, the function $\{R, U, P, D, X\} := \text{CG} (A, R e_0, X e_0)$ is obtained. For simplicity, we omit $I$ and $J$ as input, since they are constant. 

\subsubsection{Standard 2 x 2 Partitioning}

Naturally, the first choice for a partitioning is the one that is usually used for direct methods. Partitioning $R$ into $\myFlaOneByTwo{R_L}{R_R}$ implies the following partitioning for equations (\ref{eq:CGrr1}) and (\ref{eq:CGrr2}). Because of its similarity to the first equation, we will usually omit equation (\ref{eq:CGrr3}) in this section.
\begin{align*}
A \myFlaOneByTwo{P_L}{P_R} \myFlaTwoByTwo{D_{TL}}{0}{0}{D_{BR}} &= \myFlaOneByTwo{R_L}{R_R} \myFlaTwoByTwo{I - J}{0}{-H}{\underline{I} - \underline{J}} \\
\myFlaOneByTwo{P_L}{P_R} \myFlaTwoByTwo{I - U_{TL}}{- U_{TR}}{0}{I - U_{BR}} &= \myFlaOneByTwo{R_L}{\underline{R}_R}
\end{align*}
Here, $H$ is a matrix with one more row than columns that is one in the top right corner and zero everywhere else. Flattening the expressions, we obtain
\begin{align}
\myFlaOneByTwo{A P_L D_{TL} = R_L \left( I - J \right) - R_R H}{ A P_R D_{BR} = R_R \left( \underline{I} - \underline{J} \right) } \label{eq:CG:2x2s1}\\
\myFlaOneByTwo{P_L \left( I - U_{TL} \right) = R_L}{- P_L U_{TR} + P_R \left( I - U_{BR} \right) = \underline{R}_R} \text{.} \label{eq:CG:2x2s2}
\end{align}
%
%
In this form, the expressions are not matched by the pattern of the CG function. While both the right equation of (\ref{eq:CG:2x2s1}) and the left one of (\ref{eq:CG:2x2s2}) have the correct shape, one contains $R_L$ and $P_L$, while $R_R$ and $P_R$ appear in the other. In the original equations, those quantities are the same. If it is possible to match the function at all, then either all parts on the left or all parts on the right match (or both). Rewriting the equation on the left in (\ref{eq:CG:2x2s1}) as
$$A P_L D_{TL} = \myFlaOneByTwo{R_L}{R_R} \myFlaTwoByOne{ I - J }{ -H}\text{,}$$
we recognize similarities to the corresponding equation of the recurrence relation (\ref{eq:CGrr1}). However,
$$\myFlaTwoByOne{ I - J }{ -H}$$
does not have one more row than columns. While we know that only the first row of $H$ has a nonzero entry, formally, this equation does not match the pattern.

\subsubsection{3 x 3 Partitioning}

The problem above can be solved by applying a $1 \times 3$ partitioning to $R$ where the middle part is a single column. This has the effect that the first row of $H$ becomes a separate block. As CG proceeds by one column per iteration, exposing a single column appears to be a suitable choice. The partitioned operands and the resulting expressions are shown below.
\begin{gather*}
A \myFlaOneByThree{P_L}{p_M}{P_R}
\myFlaThreeByThree{D_{TL}}{0}{0}
				{0}{\delta_{MM}}{0}
				{0}{0}{D_{BR}}
=\myFlaOneByThree{R_L}{r_M}{R_R}
\myFlaThreeByThree{I - J}{0}{0}
				{-e_r^T}{1}{0}
				{0}{-e_0}{\underline{I} - \underline{J}} \\
\myFlaOneByThree{P_L}{p_M}{P_R}
\myFlaThreeByThree{I - U_{TL}}{-u_{TM}}{-U_{TR}}
				{0}{1}{-u_{MR}}
				{0}{0}{I - U_{BR}}
= \myFlaOneByThree{R_L}{r_M}{\underline{R}_R}
\end{gather*}
\begin{align*}
\myFlaOneByThree{A P_L D_{TL} = R_L \left(I - J\right) - r_M e_r^T}{A p_M \delta_{MM} = r_M - R_R e_0}{A P_R D_{BR} = R_R \left( \underline{I} - \underline{J}  \right)} \quad \\
\myFlaOneByThree{P_L \left( I - U_{TL} \right) = R_L}{-P_L u_{TM} + p_M = r_M}{-P_L U_{TR} - p_M u_{MR} + P_R \left( I - U_{BR} \right) = \underline{R}_R}
\end{align*}
Now, a similar rewriting as in the previous section yields
$$A P_L D_{TL} = \myFlaOneByTwo{R_L}{r_M} \myFlaTwoByOne{ I - J }{ -e_r^T}\text{.}$$
Here, the constant matrix has one more row than columns, and $\myFlaOneByTwo{R_L}{r_M}$ has one more column than $R_L$ in $P_L \left( I - U_{TL} \right) = R_L$. Hence, those equations, together with the additional one for $X$, are matched by the pattern of the CG function. The resulting assignment is
\begin{gather*}
\left\{ \myFlaOneByTwo{R_L}{r_M}, U_{TL}, P_L, D_{TL}, \myFlaOneByTwo{X_L}{x_M} \right\} := \qquad \qquad \qquad \qquad \qquad \\ \qquad \qquad \qquad \qquad \qquad \text{CG} \left(A, \myFlaOneByTwo{R_L}{r_M} \myFlaTwoByOne{e_0}{0}, \myFlaOneByTwo{X_L}{x_M} \myFlaTwoByOne{e_0}{0} \right) \text{.}
\end{gather*}
While the equations in the middle are not matched by the function, they can be solved to $u_{TM}$, $p_M$, $\delta_{MM}$ and $R_R e_0$, leading to computable assignments. Unfortunately, the function does not match the equations on the right-hand side. For direct methods, in a comparable situation, an equation like
$${-P_L U_{TR} - p_M u_{MR} + P_R \left( I - U_{BR} \right) = \underline{R}_R}$$
would be rewritten as 
$${P_R \left( I - U_{BR} \right) = \underline{R}_R + P_L U_{TR} + p_M u_{MR}}\text{,}$$
such that $\underline{R}_R$ is updated. Here, this is not possible. While the first column of $\underline{R}_R$ can be considered known at this point, neither the first column of $P_L U_{TR}$ nor $p_M u_{MR}$ is known, because neither $U_{TR}$ nor $u_{MR}$ are known. Apart from that, $\underline{R}_R$ can not be updated, as this would have influences on the other two equations. Finally, if quantities are updated, they are usually updated in their entirety, before they are used as input for a function. Here, $U_{TR}$ and $u_{MR}$ are not known, and we would expect them to be the output of said function, leading to circular data dependencies.

\subsubsection{Splitting off the First Column}

Since the first column of $R$ and $X$, respectively, plays a special role, an entirely different approach could be to apply a partitioning that splits off this first column. $R$ is partitioned into $\myFlaOneByTwo{r_0}{R'}$, and for the actual derivation of algorithms, $R'$ is used. The advantage is that there is a clear distinction between input and output, and it would be possible to write CG as $\{R', \ldots, X'\} := \text{CG} (A, r_0, x_0)$. As usual, partitioning $R$ and $X$ like that also implies a matching partitioning for the remaining operands:
\begin{align*}
A \myFlaOneByTwo{p_0}{P'}
\myFlaTwoByTwo	{\delta_{0}}{0}
				{0}{D'}
&=\myFlaOneByTwo{r_0 }{ R'}
\myFlaTwoByTwo	{1}{0}
				{-e_0}{\underline{I} - \underline{J}} \\
\myFlaOneByTwo{p_0}{P'}
\myFlaTwoByTwo	{1}{-u'}
				{0}{I - U'}
&= \myFlaOneByTwo{r_0}{\underline{R}'}
\end{align*}
Flattening those expression, we immediately obtain a value for $p_0$:
\begin{gather*}
\myFlaOneByTwo{A p_0 \delta_0 = r_0 - R' e_0}{ A P' D' = R' \left( \underline{I} - \underline{J} \right)} \\
\myFlaOneByTwo{p_0 = r_0}{ - p_0 u' + P' \left(I - U' \right) = \underline{R}'}
\end{gather*}
Using the orthogonality of $R$, it is also possible to find an assignment for $\delta_0$. Thus, this partitioning allows to compute all quantities of the first iteration and declare them as known. Unfortunately, splitting off $r_0$ forces us to also split off the first row of $U$. At this point, it is not possible to compute it in its entirety. Furthermore, because of $-p_0 u'$, the equations on the right-hand side do not have the same shape as the original description of CG. While it is of course still possible to compute $u'$, it is only possible entry by entry, adding an additional assignment to any update we can derive. Consequently, it will not be possible to derive the updates for CG usually found in textbooks.


\subsubsection{Divide and Conquer}

So far, no partitioning enabled us to describe one CG operation as multiple, smaller CG operations, if necessary, with updated quantities as input, similar to how the PME of the triangular system in Section \ref{sec:triLS} contains the function $\Phi$ two times. This is, however, possible, if we define a generalized version of CG. The disadvantage is that it requires a deeper understanding of the algorithm, in addition to some knowledge that is initially not available when deriving algorithms solely based on their matrix representation.

To derive such a representation, we partition $R$ into \smash{$\myFlaOneByFour{R_0}{r_1}{r_2}{R_3}$}. Since the resulting partitioned postcondition is very large, it is not shown here. The generalized version of CG requires some parts of $P$ as an additional argument. Initially, it is only the first column $P$, here denoted with $P e_0$: $\{R, U, P, D, X\} := \text{CG} (A, R e_0, X e_0, P e_0)$. Since it is equal to the first column of $R$, it can be considered known. For now, we will assume that $A$ is nonsymmetric.

Our goal is now to write this CG operation as two separate ones, one covering \smash{$\myFlaOneByTwo{R_0}{r_1}$}, and one for \smash{$\myFlaOneByTwo{r_2}{R_3}$}. Clearly, the first one is
$$\left\{\myFlaOneByThree{R_0}{r_1}{r_2}, \ldots, \myFlaOneByTwo{P_0}{p_1}, \ldots, \myFlaOneByThree{X_0}{x_1}{x_2}\right\} := \text{CG} \left(A, R_0 e_0, X_0 e_0, P_0 e_0 \right) \text{,}$$
omitting some of the output in the interest of legibility. Now, to find the correct arguments for the second one, we need to know how $p_2$ is computed:
$$p_2 = r_2 + P_0 u_{02} + p_1 \nu_{12}$$
$u_{02}$ and $\nu_{12}$ are not part of the output of the function above, but they can in turn be computed using known quantities only:
\begin{align*}
u_{02} &= \left(- P_0^T A P_0 \right)^{-1} \cdot P_0^T A r_2 \\
\nu_{12} &= - \frac{p_1^T A r_2 + p_1^T A P_0 u_{02}}{p_1^T A p_1}
\end{align*}
To eliminate them entirely, we can also write 
$$p_2 = r_2 - \myFlaOneByTwo{P_0}{p_1} \myFlaTwoByTwo{P_0^T A P_0}{0}{p_1^T A P_0}{p_1^T A p_1}^{-1} \myFlaTwoByOne{P_0^T A r_2}{p_1^T A r_2} \text{.}$$
Now, we can write the second part as
$$\{R_3, \ldots, P_{3}, \ldots, X_3\} := \text{CG} \left(A, r_2, x_2, \myFlaOneByThree{P_0}{p_1}{p_2} \right) \text{.}$$
To ensure that this function computes the correct sequence of search directions, it also has to use the ones computed by the first function, which are $\myFlaOneByTwo{P_0}{p_1}$. This is the reason why not just $p_2$, but $\myFlaOneByThree{P_0}{p_1}{p_2}$ is input.

This situation is slightly different if $A$ is symmetric. In that case, $U$ is upper diagonal and $P^T A P$ is diagonal, so $p_2$ is computed as
$$p_2 = r_2 - p_2 \frac{p_1^T A r_2}{p_1^T A p_1}\text{.}$$
Now, each search direction is computed using only the last one, so $\myFlaOneByTwo{P_0}{p_1}$ is not needed as input. Thus, the second function simplifies to
$$\left\{R_3, \ldots, P_{3}, \ldots, X_3 \right\} := \text{CG} \left(A, r_2, x_2, p_2\right) \text{.}$$
As mentioned before, the disadvantage is that we already need to know how some quantities are computed to derive this representation, while it is actually our goal to find those updates.

\subsubsection{Splitting off the Last Column}

Among those presented in this section, the $3 \times 3$ partitioning that exposes a single column is the only one that resulted in expressions that were naturally matched by the CG function. The problem of this partitioning is that there is no obvious way how to deal with the right-hand side parts.

To find one that better suits iterative methods, it is helpful to again inspect the differences to direct methods. After all, the $3 \times 3$ partitioning came to existence as a modification of the standard partitioning used for direct methods. As mentioned before, with direct methods, the sizes of all operands are initially known. Thus, at any point during the computation, there are (potentially empty) parts of operands that are already computed, and (potentially empty) parts that are not computed yet. The loop invariant describes those parts that are already computed at the beginning of the loop body. Conversely, one can think of those parts of the PME that are not part of the loop invariant as those parts that are not computed yet.

To see that this is consistent, recall that the reason that the full set of nodes can never be a feasible loop invariant is that it would imply that the solution is already computed before the loop is entered. Clearly, this is equivalent to saying that there are no parts left that are not computed yet, even before the loop is entered.

In case of iterative methods, there is little use in talking about parts of operands that are not yet computed beyond the current iteration, as each iteration might be the last.
We can conclude that it makes little sense to use a partitioning where the right-hand side is more than a single column.

The solution is to use a partitioning that is a hybrid of the standard $2 \times 2$ partitioning and the $3 \times 3$ partitioning that exposes a single column: A partitioning that splits off the last column: $P$ is partitioned into $\myFlaOneByTwo{P_L}{p_R}$. Because of the additional column of $R$, it is partitioned into $\myFlaOneByThree{R_L}{r_R}{r_+}$. Thus, for CG, we obtain
\label{eq:partitionedPostconditionNonsymCG}
\begin{align*}
A \myFlaOneByTwo{P_L}{p_R} \myFlaTwoByTwo{D_{TL}}{0}{0}{\delta_{BR}} &= \myFlaOneByThree{R_L}{r_R}{r_+} \myFlaThreeByTwo{I - J}{0}{-e_r^T}{1}{0}{-1} \\
\myFlaOneByTwo{P_L}{p_R}\myFlaTwoByTwo{I- U_{TL}}{- u_{TR}}{0}{1} &= \myFlaOneByTwo{R_L}{r_R} \\
\myFlaOneByTwo{P_L}{p_R} \myFlaTwoByTwo{D_{TL}}{0}{0}{\delta_{BR}} &= \myFlaOneByThree{X_L}{x_R}{x_+} \myFlaThreeByTwo{I - J}{0}{-e_r^T}{1}{0}{-1} \text{.}
\end{align*}
Flattening those expressions yields
\begin{gather*}
\myFlaOneByTwo{A P_L D_{TL} = R_L \left( I - J \right) - r_R e_r^T}{ A p_R \delta_{BR} = r_R - r_+ } \\
\myFlaOneByTwo{P_L \left( I - U_{TL} \right) = R_L}{ - P_L u_{TR} + p_R = r_R} \\
\myFlaOneByTwo{ P_L D_{TL} = X_L \left( I - J \right) - x_R e_r^T}{ p_R \delta_{BR} = x_R - x_+ } \text{.}
\end{gather*}
The left-hand side parts are now matched by
\begin{gather*}
\left\{ \myFlaOneByTwo{R_L}{r_R}, U_{TL}, P_L, D_{TL}, \myFlaOneByTwo{X_L}{x_R} \right\} := \qquad \qquad \qquad \qquad \qquad \\ \qquad \qquad \qquad \qquad \qquad \text{CG} \left(A, \myFlaOneByTwo{R_L}{r_R} \myFlaTwoByOne{e_0}{0}, \myFlaOneByTwo{X_L}{x_R} \myFlaTwoByOne{e_0}{0} \right) \text{.}
\end{gather*}
Using some of the properties derived with the approach presented in Section \ref{sec:propertyDerivation}, the equations on the right can be solved to all remaining unknowns. Thus, we obtain a PME with assignments for every unknown quantity. The systematic derivation of loop-based algorithms, using this partitioning, is presented in the next chapter.

Note that this PME can be interpreted as an ``inductive PME'': Assuming it is possible to compute an arbitrary number of previous iterations, it is possible to compute one additional iteration. The previous iterations are represented by the CG function, and the additional iteration is computed using the remaining, explicit assignments of the PME. The base case is obtained by assuming all left and top left parts to be empty.


\chapter{Derivation of Algorithms for Iterative Methods}
\label{chap:derivationIM}

After having laid the foundations in the previous chapter, the actual approach for deriving algorithms for iterative methods is presented in this chapter. The approach itself can be found in Section \ref{sec:derivationIterativeMethods}, followed by two examples in Section \ref{sec:derivationOfAlgorithmsExample} and \ref{sec:BiCGexample}. In the final Section \ref{sec:scopeLimitations}, scope and limitations of the approach are discussed.

\section{Derivation of Algorithms}
\label{sec:derivationIterativeMethods}

The derivation of algorithms for iterative methods mainly follows the same three basic steps as for direct methods. First, one or more PMEs are generated. In the second step, dependency graphs are constructed, which are then used to select loop invariants. In the third and final step, from each loop invariant, an algorithm is constructed. In this section, we will present this process for iterative methods. There is, however, a new fourth step. In this step, some postprocessing is applied to the derived algorithms to generate a number of variants that may behave differently in floating point arithmetic or vary in their performance. If applicable, we follow the same structure as \cite{Fabregat-Traver:thesis}.

\subsection{PME Generation}

The first stage towards the generation of algorithms is to find PMEs. The necessary steps are explained in the following. There are a number of differences compared to direct methods: It is necessary to derive properties of matrices, which then enable the system to solve equations. Additionally, different partitionings are needed to deal with new types of operands.

\subsubsection{Derivation of Properties}

The approach presented in Section \ref{sec:propertyDerivation} is used to derive properties.

\subsubsection{Initial Partitioning}

Operands of the postcondition are partitioned depending on their shape and properties. For objects that are completely known or unknown, the applicable partitionings are similar to the ones used for direct methods. The difference lies in the sizes of the resulting objects. Only the top, left, and top left parts are matrices, the remaining ones are either vectors or scalars:
\begin{align*}
B &\rightarrow B	&	B &\rightarrow \myFlaOneByTwo{B_L}{b_R} \\
B &\rightarrow \myFlaTwoByOne{B_T}{b_B}	&	B &\rightarrow \myFlaTwoByTwo{B_{TL}}{b_{TR}}{b_{BL}}{\beta_{BR}}
\end{align*}
%
Just as with direct methods, triangular or symmetric matrices are either not partitioned at all, or the $2 \times 2$ partitioning is used. In case of the latter, the top left part is required to be square, such that it inherits the property of the matrix.

Separate partitionings are necessary to deal with matrices where initially, only the first column is known. If the last column of those matrices is omitted, the usual $1 \times 2$ partitioning is applied:
$$\underline{B} \rightarrow \myFlaOneByTwo{B_L}{b_R}$$
In case of the complete matrix, an additional column is obtained:
$$B \rightarrow \myFlaOneByThree{B_L}{b_R}{b_{+}}$$
The constant matrices $\underline{J}$ and $\underline{I}$ pose a special case. They have one more row than columns, so $\underline{J}$ is not lower diagonal and $\underline{I}$ is not an identity matrix. To derive algorithms, it is important to utilize their specific structure, so we will provide the partitionings explicitly:
\begin{align*}
\underline{J} \rightarrow \myFlaThreeByTwo{J}{0}{e_r^T}{0}{0}{1} \qquad \underline{I} \rightarrow \myFlaThreeByTwo{I}{0}{0}{1}{0}{0}
\end{align*}
The same partitioning that is applied to the operands in the postcondition is also applied to the operands in the set of properties. Then, properties of expressions of partitioned operands are derived, similar to how partitioned operands inherit properties. If for example $B$ is partitioned into $\myFlaOneByTwo{B_L}{b_R}$, and $\text{\ttfamily Diagonal} \left[B^T B\right]$ is contained in the set of properties, the system obtains
$$\text{\ttfamily Diagonal} \left[\myFlaTwoByOne{B_L^T}{b_R^T} \myFlaOneByTwo{B_L}{b_R} \right] = \text{\ttfamily Diagonal} \left[\myFlaTwoByTwo{B_L^T B_L}{B_L^T b_R}{b_R^T B_L}{b_R^T b_R} \right] \text{.}$$
Thus, it is possible to derive that $B_L^T B_L$ is diagonal as well, and $B_L^T b_R$ and $b_R^T B_L$ are zero.

\subsubsection{Finding the PME}

The first part of this step consists of performing symbolic arithmetic and distributing equalities across the partitionings. Consider the Krylov sequence as an example:
\begin{align*}
K:= \text{KS} ( A, K e_0) \equiv
\left\{
\begin{aligned}
P_\text{pre}: \{ &\text{\ttfamily Input}(A) \land \text{\ttfamily Matrix}(A) \\
		&\text{\ttfamily Matrix}[\underline{J}] \land \text{\ttfamily LowerDiagonalR}[\underline{J}] \land \\
		&\text{\ttfamily FirstColumnInput}(K) \} \\
P_\text{post}: \{ &A \underline{K} = K \underline{J} \\
			&\text{size}(K) = n \times m \}
\end{aligned}
\right.
\label{eq:KrylovSeq:Description}
\end{align*}
The partitioned postcondition looks as follows:
$$A \myFlaOneByTwo{K_L}{k_R} = \myFlaOneByThree{K_L}{k_R}{k_{+}} \myFlaThreeByTwo{J}{0}{e_r^T}{0}{0}{1}$$
It can be rewritten as the following expression:
$$\myFlaOneByTwo{A K_L = K_L J + k_R e_r^T}{ A k_R = k_{+}}$$
%
In the second part, the goal is to find a representation of this expression where the value of each unknown quantity is determined by an assignment, using known operations. The quantities on the right-hand side of that assignment either have to be known, or their value is determined by another assignment. Intuitively, one could say the goal is to make this expression computable. This is done by matching patterns of known functions and operations.

Rewriting the equation on the left-hand side as
$$A K_L = \myFlaOneByTwo{K_L}{k_R} \myFlaTwoByOne{J}{e_r^T} \text{,}$$
it is easy to see that it describes the computation of a Krylov sequence, so it is possible to use the function from the description of the operation:
%
$$\myFlaOneByTwo{K_L}{k_R} := \text{KS} \left(A, \myFlaOneByTwo{K_L}{k_R} \myFlaTwoByOne{e_0}{0} \right)$$
Now, $K_L$ and $k_R$ can be considered known, so all that remains is to find an assignment for $k_+$. In this example, the equation on the right hand side is already solved to $k_+$, so
$$k_+ := A k_R$$
is immediately obtained. The PME then is 
$$\myFlaOneByTwo{ \myFlaOneByTwo{K_L}{k_R} := \text{KS} \left(A, \myFlaOneByTwo{K_L}{k_R} \myFlaTwoByOne{e_0}{0} \right) }{ k_+ := A k_R} \text{.}$$
In general, finding assignments is not that straightforward. Usually, it is necessary to apply known properties. Similar to the ``Application of Properties'' step in Section \ref{sec:propertyDerivation:Derivation}, parts of expressions with known properties are multiplied to both sides of equations, with the intention to recreate said expressions. Unfortunately, it is not easy to determine beforehand if the application of a property allows to solve an equation. Thus, in the manner of an exhaustive search, all matching properties are applied. If by this means, multiple assignments for the same expression are found, separate PMEs are derived for each variant.

\subsubsection{Remark on Recursive Algorithms}
The PMEs for direct methods immediately lead to recursive algorithms \cite{Bientinesi:thesis, Fabregat-Traver:thesis}. This is also true for iterative methods. The difference is that in case of direct methods, those algorithms are usually divide and conquer algorithms. For iterative methods, the PME naturally leads to a ``head recursive'' implementation, that is, the recursive function call is the first operation in the function body.

\subsection{Loop Invariant Identification}

For the second step, the identification of loop invariants, constraints for the feasibility of loop invariants are introduced that differ from the ones used for direct methods. To understand why they are introduced, it is helpful to get an intuition for what those loop invariants express, and how they differ from the ones for direct methods.

\paragraph{Loop Invariants for Iterative Methods}

For iterative methods, there is usually a (partial) ordering in which parts of different operands can be computed. Let us assume there are three matrices $B$, $C$ and $D$. During the execution of the algorithm, the following sequence is computed, where $b_k$, $c_k$ and $d_K$ are columns of $B$, $C$ and $D$, respectively:
$$\ldots b_i \; c_i \; d_i \; b_{i+1} \; c_{i+1} \; d_{i+1} \ldots$$
A variant derived from one loop invariant may compute $c_i \; d_i \; b_{i+1}$ in one iteration. In contrast, one derived from another loop invariant may compute $d_i \; b_{i+1} \; c_{i+1}$.
%
%
Since no quantities are overwritten, a loop invariant only expresses at which point in the sequence above an iteration starts.

In contrast, different loop invariants for direct methods may result in algorithms that compute the solution by row or by column.

\subsubsection{Graph of Dependencies}

The construction of the dependency graph is even simpler compared to direct methods. The assignments of the PME are not decomposed into basic building blocks. Each assignment is represented by one node in the dependency graph. The dependencies are established as usual.

\subsubsection{Subset Selection}

Just as with direct methods, subsets of nodes of the dependency graph are selected as candidates for loop invariants. Again, for every node that is contained in a subset, all preceding nodes have to be in that set, too. To assess the feasibility of loop invariants, however, slightly modified constraints have to be imposed:
%
\begin{enumerate}
\item There must exist a basic initialization of the operands, that is, an initial partitioning, followed by some preprocessing operations, that renders the predicate $P_\text{\textrm{inv}}$ true:
\begin{align*}
&\{P_\text{\textrm{pre}}\} \\
&\text{\textbf{Partition}} \\
&\text{Preprocessing} \\
&\{P_\text{\textrm{inv}}\}
\end{align*}
\item $P_\text{\textrm{inv}}$ and the negation of the loop guard, $G$, must imply the postcondition, $P_\text{\textrm{post}}$:
$$P_\text{\textrm{inv}} \land \neg G \Rightarrow P_\text{\textrm{post}}$$
\end{enumerate}

We begin with discussing the first condition. How the operands are partitioned was already established in the Section ``Initial Partitioning'', but the initial sizes of the partitioned operands were not specified. The iterative methods covered in this thesis proceed through those matrices where initially, only the first column is known, from the left to the right. As a consequence, the left parts of those matrices are initially empty. The initial block sizes for all partitionings are given in Table \ref{tab:initialBlockDimensions}.
\begin{table}[htp]
\begin{center}
\begin{tabular}{l l}
\toprule
Initial Partitioning	&	Dimensions \\ \midrule
$B \rightarrow \myFlaOneByTwo{B_L}{b_R}$	&	$B_L$ is $n \times 0$ \\
$B \rightarrow \myFlaTwoByOne{B_T}{b_R}$	&	$B_T$ is $0 \times n$ \\
$B \rightarrow \myFlaTwoByTwo{B_{TL}}{b_{TR}}{b_{BL}}{\beta_{BR}}$	&	$B_{TL}$ is $0 \times 0$ \\
$\underline{B} \rightarrow \myFlaOneByTwo{B_L}{b_R}$	&	$B_L$ is $n \times 0$ \\
$B \rightarrow \myFlaOneByThree{B_L}{b_R}{b_+}$		&	$B_L$ is $n \times 0$ \\ \bottomrule
\end{tabular}
\end{center}
\caption{Initial sizes of partitioned operands.}
\label{tab:initialBlockDimensions}
\end{table}%

Recall that for direct methods, the full set of nodes of the dependency graph can never be a feasible loop invariant. Formally, the reason is that there is no initial partitioning that renders this loop invariant true. Alternatively, one can say that it implies that the complete solution is already computed even before the loop is entered. For iterative methods, this is different. Due to the initial partitioning, which does not expose blocks that represent the remaining iterations, there is no subset implying that the entire solution is already computed. Thus, the full set is a feasible loop invariant.

However, further loop invariant candidates for iterative methods are not rendered true by the initial partitioning for a different reason. This is the case for some subsets that contain more than the first node. The reason is that the additional nodes may compute quantities that are not empty in the initial partitioning. Consider the Krylov sequence as an example: $k_R$ and $k_+$ remain $n \times 1$ vectors in the initial partitioning. Thus, the assignment $k_+ := A k_R$ computes a non-empty quantity, even if the initial partitioning is applied. However, at the beginning of the operation, only $k_R$ is known, which is the fist column of $K$ in the initial partitioning, while $k_+$ is the second column, which is unknown.

Fortunately, the equations from those additional nodes can be rendered true by computing those quantities with some preprocessing operations. Naturally, the necessary preprocessing operations are obtained by applying the initial partitioning to the operations of those nodes. In this example, the preprocessing operation is $k_+ := A k_R$. Those preprocessing operations can consist of the same type of operations as the actual update operations.

To demonstrate how the first constraint is checked, we return to the example of the Krylov sequence. We start with the following loop invariant candidate:
$$\myFlaOneByTwo{ \myFlaOneByTwo{K_L}{k_R} := \text{KS} \left(A, \myFlaOneByTwo{K_L}{k_R} \myFlaTwoByOne{e_0}{0} \right) }{ \neq}$$
If $K_L$ has size $n \times 0$, the expression on the left-hand side becomes
$k_R := \text{KS} \left(A, k_R \right)$.
Initially, the first column of $K$, which is now $k_R$, is known. Thus, since both sides of this assignment are known, this expression is considered to be true. This implies that the loop invariant satisfies the first constraint.

The second candidate for a loop invariant is 
$$\myFlaOneByTwo{ \myFlaOneByTwo{K_L}{k_R} := \text{KS} \left(A, \myFlaOneByTwo{K_L}{k_R} \myFlaTwoByOne{e_0}{0} \right) }{ k_+ := A k_R} \text{.}$$
Note that in the interest of simplicity, we usually just write $K_L e_0$ instead of
$$\myFlaOneByTwo{K_L}{k_R} \myFlaTwoByOne{e_0}{0}\text{,}$$
if it is sufficient. The following expression is obtained if the initial partitioning is applied, where $K_L$ is empty:
$$\myFlaOneByTwo{k_R := \text{KS} \left(A, k_R \right)}{k_+ := A k_R}$$
As $k_+$ is not known, the initial partitioning alone does not render this expression true. This, however, can be solved with a preprocessing operation. Here, the operation is $k_+ := A k_R$. It follows that this loop invariant satisfies the first constraint as well.

To check if the second condition for the feasibility of loop invariants is satisfied, we first need to determine the loop guard. As discussed in Section $\ref{sec:matrixRepresentationIntroduction}$, comparing the size of a growing block of a matrix to the size of the entire matrix does not work, as the matrix grows as well. For this reason, the additional predicates were added to the postcondition (by hand), which are now easily translated into loop guards (automatically). How this is done is shown in Table \ref{tab:loopGuards}.
\begin{table}[htp]
\begin{center}
\begin{tabular}{lll}
\toprule
 & \multicolumn{2}{c}{Loop guard} \\ \cmidrule{2-3}
Predicate & $\text{\ttfamily FirstColumnInput}[B]$ & $\text{\ttfamily Output}[B]$ \\ \midrule
$\| B e_r^T \| < \varepsilon$ & $\| b_R \| \geq \varepsilon$ & $\| B_L e_r^T \| \geq \varepsilon$ \\ 
$\text{size}(B) = n \times k$ & $\text{size}\left( \myFlaOneByTwo{B_L}{b_R} \right) < n \times k$ & $\text{size}\left( B_L \right) < n \times k$ \\
$\| B e_r^T - B e_{r-1}^T \| < \varepsilon$ & $\| b_R - B_L e_{r}^T \| \geq \varepsilon$ & $\| B_L e_r^T - B_L e_{r-1}^T \| \geq \varepsilon$ \\
\bottomrule
\end{tabular}
\end{center}
\caption{Look-up table for determining loop guards for iterative methods. The row is selected according to the additional predicate in the postcondition. The column is selected depending on the property of the operand that appears in the position of $B$ in that predicate. Example: The predicate is $\| R e_r^T \| < \varepsilon$. The precondition contains the property $\text{\ttfamily FirstColumnInput}[R]$. Thus, the loop guard is $\| r_R \| \geq \varepsilon$. To allow for other loop guards, this table has to be extended manually.}
\label{tab:loopGuards}
\end{table}%

Note that even though $B e_r^T$ in the predicate $\| B e_r^T \| < \varepsilon$ refers to the last column of $B$, $b_R$ in the loop guard $\| b_R \| \geq \varepsilon$ and $B_L e_r^T$ in $\| B_L e_r^T \| \geq \varepsilon$, respectively, are the second to last columns. Similarly, the other loop guards omit the last column as well. To understand why this is necessary, we have to look at the loop invariant candidates again. Just as with direct methods, the empty set can never be a valid loop invariant because it corresponds to the empty predicate. Thus, all remaining candidates contain at least the first node of the dependency graph. For iterative methods, this first node always represents the operation itself, that is, it contains the original function. Which parts of the operands are output of that function depends on the properties of the operands:
\begin{itemize}
\item[-] If initially, the first column of $B$ is known ($\text{\ttfamily FirstColumnInput}[B]$), $B_L$ and $b_R$ are output of the function.
\item[-] If $B$ is initially unknown ($\text{\ttfamily Output}[B]$), just $B_L$ is output.
\end{itemize}
Consider nonsymmetric CG as an example: The operation of the first node is
\begin{gather*}
\left\{ \myFlaOneByTwo{R_L}{r_R}, U_{TL}, P_L, D_{TL}, \myFlaOneByTwo{X_L}{x_R} \right\} :=  \text{CG} \left(A, R_L e_0, X_L e_0 \right) \text{.}
\end{gather*}
The first columns of $R$ and $X$, here denoted by $R_L e_0$ and $X_L e_0$, are initially known and the function computes $R_L$, $r_R$, $X_L$ and $x_R$. In contrast, the property of $P$ is $\text{\ttfamily Output}[P]$, and only $P_L$ is computed. Intuitively, the reason is that the function always computes the same number of columns of all operands. If initially, one column of an operand is already known, in the end, one additional column is obtained.

Since only the first node is guaranteed to be part of the loop invariant, only those parts computed in the first node are guaranteed to be known at the beginning and at the end of the loop. Thus, for $\text{\ttfamily FirstColumnInput}[B]$, checking the norm of $b_R$ in the loop guard is always possible, but using $b_+$ is not.

While there are variants which compute $b_+$ (or $b_R$ in case of $\text{\ttfamily Output}[B]$), and those would allow different loop guards, we refrain from using those to keep the derivation simple. This, however, means that if $b_+$ (or $b_R$) are computed, they will not be regarded as part of the solution. That is possible because $B$ does not have a fixed size. Unfortunately, this has the effect that some algorithms compute results that are subsequently discarded. If, however, the loop body of those algorithms has a reduced computational complexity, this is an acceptable tradeoff.

We demonstrate how the second condition for the feasibility of loop invariants is checked using the following loop invariant candidate $P_\text{inv}$ as an example:
\begin{align*}
\left\{ \myFlaOneByTwo{R_L}{r_R}, U_{TL}, P_L, D_{TL}, \myFlaOneByTwo{X_L}{x_R} \right\} &:= \text{CG} \left(A, R_L e_0, X_L e_0 \right) \\
u_{TR} &:= - \left( P_L^T A P_L \right)^{-1} P_L^T A r_R
\end{align*}
The loop guard $G$ is $\| r_R \| \geq \varepsilon$, so its negation $\neg G$ is $\| r_R \| < \varepsilon$. Now, it has to be checked whether $P_\text{inv}$ and $\neg G$ imply the postcondition $P_\text{post}$:
\begin{align*}
A P D &= R \left( \underline{I} - \underline{J}  \right)\\
P \left( I - U \right) &= \underline{R} \\
P D &= X \left( \underline{I} - \underline{J} \right)\\
\| R e_r^T \| &< \varepsilon
\end{align*}
To do that, we rewrite $P_\text{inv}$ as
\begin{align*}
A P_L D_{TL} &= \myFlaOneByTwo{R_L}{r_R} \myFlaTwoByTwo{I - J}{0}{-e_r^T}{1}\\
P_L \left( I - U_{TL} \right) &= R_L \\
P_L D_{TL} &= \myFlaOneByTwo{X_L}{x_R} \myFlaTwoByTwo{I - J}{0}{-e_r^T}{1} \\
u_{TR} &= - \left( P_L^T A P_L \right)^{-1} P_L^T A r_R \text{.}
\end{align*}
Clearly, by rewriting $\myFlaOneByTwo{R_L}{r_R}$ as $R$, $P_L$ as $P$ and so on, one can see that the first three equations above are the same as the equations in the postcondition. Furthermore, the negation of the loop guard, $\| r_R \| < \varepsilon$, refers to the last column of $\myFlaOneByTwo{R_L}{r_R}$, just like $\| R e_r^T \| < \varepsilon$ refers to the last column of $R$. Thus, $P_\text{inv} \land \neg G$ implies the postcondition. While $P_\text{inv} \land \neg G$ also implies the equation
$$u_{TR} = - \left( P_L^T A P_L \right)^{-1} P_L^T A r_R \text{,}$$
this equation is not needed to render the postcondition true. In the algorithm, $u_{TR}$ will be discarded.

Note that the initial partitioning exposed one additional column of each operand. $R$ for example was partitioned into $\myFlaOneByThree{R_L}{r_R}{r_+}$. In the postcondition, $R$ just consists of $R_L$ and $r_R$. As mentioned before, this is possible because the number of columns of $R$ is variable.







\subsection{Algorithm Construction}

The algorithm is constructed in the third step. In addition to the update, preprocessing operations have to be determined. Furthermore, in line with different rules for the initial partitioning, the repartitioning is modified. The process of identifying the update operations does not change.

\subsubsection{Preprocessing}

As mentioned before, the preprocessing operations are obtained by applying the initial partitioning (Table \ref{tab:initialBlockDimensions}) to all but the initial node contained in the loop invariant. The first node is excluded because when the initial partitioning is applied to it, it always reduces to an expression similar to $k_R := \text{KS} \left(A, k_R \right)$ for the Krylov sequence. Since the output of this function is already known, nothing has to be computed. A FLAME worksheet, extended by the preprocessing, is shown in Figure \ref{fig:ws:emptyP}.

\begin{figure}
\centering
\begin{minipage}[t]{2.35in}
	\resetsteps
	\setboolean{BlockedAlgQ}{false}
	
	\renewcommand{\WSoperation}{ $\ldots$}
	
	\renewcommand{\WSupdate}{\text{Update}}
	
	{
	\worksheetGrayNoNumbersEmptyP
	}
\end{minipage}
\caption{The skeleton of a FLAME worksheet, extended by some preprocessing.}
\label{fig:ws:emptyP}
\end{figure}

\subsubsection{Repartitioning the Operands}

To ensure that the resulting algorithms make progress, the operands have to be repartitioned. The sizes of some operands depend on the number of iterations that is computed, so with every iteration, their sizes have to grow. This is done by adding rows and/or columns in the ``Continue with'' repartitioning. The rules are shown in Table \ref{tab:repartitionIterativeMethods}.


%
\begin{table}[htp]
\begin{center}
\begin{tabular}{ccc}
\toprule
Initial Partitioning	&	Repartition	&	Continue with \\ \midrule
$\myFlaOneByTwoI{B_L}{b_R}$	&
$\myFlaOneByTwoI{B_0}{b_1}$	&
$\myFlaOneByThreeLI{B_0}{b_1}{b_2}$ \\
$\myFlaTwoByOneI{B_T}{b_R}$	&
$\myFlaTwoByOneI{B_0}{b_1}$ 	&
$\FlaThreeByOneTI{B_0}{b_1}{b_2}$ \\
$\myFlaTwoByTwoI{B_{TL}}{b_{TR}}{b_{BL}}{\beta_{BR}}$	&
$\myFlaTwoByTwoI{B_{00}}{b_{01}}{b_{10}}{\beta_{11}}$ &
$\myFlaThreeByThreeTLI	{B_{00}}	{b_{01}}		{b_{02}}
					{b_{10}}	{\beta_{11}}	{\beta_{12}}
					{b_{20}}	{\beta_{21}}	{\beta_{22}}$ \\
$\myFlaOneByThreeI{B_L}{b_R}{b_+}$		&
$\myFlaOneByThreeI{B_0}{b_1}{b_2}$ &
$\myFlaOneByFourTFF{B_0}{b_1}{b_2}{b_3}$\\ \bottomrule
\end{tabular}
\end{center}
\caption{``Repartition'' and ``Continue with'' rules for iterative methods.}
\label{tab:repartitionIterativeMethods}
\end{table}%

\subsubsection{Predicates $P_\text{\textrm{before}}$ and $P_\text{\textrm{after}}$}

To obtain the predicates $P_\text{\textrm{before}}$ and $P_\text{\textrm{after}}$, the repartitioned operands are plugged into the loop invariant. The resulting expressions are flattened, using the PME if necessary. What is obtained by applying the ``Repartition'' rules to the loop invariant becomes $P_\text{before}$. Applying the ``Continue with'' partitioning result in $P_\text{after}$.

We demonstrate this for the loop invariant
$$\myFlaOneByTwo{ \myFlaOneByTwo{K_L}{k_R} := \text{KS} \left(A, K_L e_0 \right) }{ \neq} \text{.}$$
The ``Repartition'' and ``Continue with'' rules for $K$ are shown below, introducing the newly added $k_3$:
\begin{align*}
\myFlaOneByThreeI{K_L}{k_R}{k_+} &\rightarrow \myFlaOneByThreeI{K_0}{k_1}{k_2} \\
\myFlaOneByThreeI{K_L}{k_R}{k_+} &\leftarrow \myFlaOneByFourTFF{K_0}{k_1}{k_2}{k_3} \text{.}
\end{align*}
Applying the ``Repartition'' rules to the loop invariant yields the following predicate $P_\text{before}$:
$$ \myFlaOneByTwo{K_0}{k_1} := \text{KS} \left(A, K_0 e_0 \right)$$
Using the ``Continue with'' repartitioning, we obtain the expression
$$\myFlaOneByThree{K_0}{k_1}{k_2} := \text{KS} \left(A, K_0 e_0 \right)\text{.}$$
To flatten this expression, the PME is used. This results in the following $P_\text{after}$:
$$\myFlaOneByTwo{ \myFlaOneByTwo{K_0}{k_1} := \text{KS} \left(A, K_0 e_0 \right) }{ k_2 := A k_1}$$
\subsubsection{Finding the Updates}

The difference between $P_\text{\textrm{after}}$ and $P_\text{before}$ now is the update. Identifying the differences is much easier compared to direct methods because they are always entire equations, not subexpression. The reason is that no quantities are updated.

For the Krylov sequence, $k_2 := A k_1$ can easily be identified as the difference between $P_\text{\textrm{after}}$ and $P_\text{before}$.

\subsection{Refinement}
\label{sec:postprocessing}


In practice, algorithms would not, and depending on the language, can not be implemented exactly like they are derived with the presented approach. Consider the following three assignments as an example. They are part of the update of one nonsymmetric CG algorithm.
\begin{align*}
r_2 &:= r_1 - A p_1 \delta_{11} \\
u_{02} &:= \left(- P_0^T A P_0 \right)^{-1} \cdot P_0^T A r_2 \\
\nu_{12} &:= - \frac{p_1^T A r_2 + p_1^T A P_0 u_{02}}{p_1^T A p_1}
\end{align*}
To translate the assignments above to C for example, they have to be decomposed into basic operations that are implemented in a library like BLAS. Since such libraries usually do not include functions for general products of more than two quantities, auxiliary variables have to be introduced.

While the assignments could immediately be translated into Matlab code, this would not result in efficient code. Clearly, some subexpressions appear multiple times, so it is preferable to introduce auxiliary variables for those and compute their values just once. This is referred to as common subexpression elimination. While this concept is well known in the domain of compiler construction \cite{steven1997advanced}, to the author's knowledge there is no research on the elimination of \emph{overlapping} common subexpressions. Overlapping common subexpression are common subexpression that can not be eliminated at the same time. Consider the three terms $P_0^T A P_0$, $p_1^T A r_2$ and $p_1^T A P_0$ for a simple example. $P_0^T A P_0$ and $p_1^T A P_0$ have $AP_0$ in common, $p_1^T A r_2$ and $p_1^T A P_0$ share $p_1^T A$. In $p_1^T A P_0$, those two common subexpression overlap since $A$ is part of both.

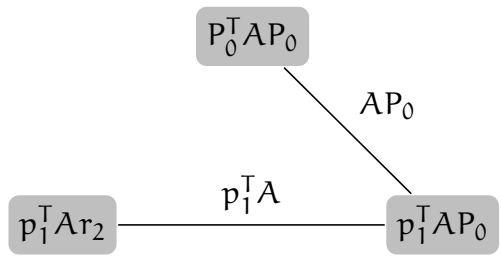
\begin{figure}[h]
\centering
\begin{tikzpicture}[node distance=2.5cm and 2.5cm]

\node[rect]	(A)							{$P_0^T A P_0$};

\node[rect]	(B)	[below left=of A]	{$p_1^T A r_2$};

\node[rect]	(C)	[below right=of A, xshift=0.0cm]		{$p_1^T A P_0$};

\path[-]	(A)		edge 		node {$AP_0$}		(C);
		
\path[-]	(B)		edge			node {$p_1^T A$}		(C);

\end{tikzpicture}
\caption{Elimination of common subexpressions represented as a graph.}
\label{fig:CSE}
\end{figure}

Finding a replacement of subexpressions that is optimal in the sense that it has the lowest computational cost is not trivial. Inspecting the problem, one observes that a simplified version of it can be mapped to a \emph{maximum weight matching} problem, which is solvable in polynomial time \cite{edmonds1965maximum}. Again, consider the terms $P_0^T A P_0$, $p_1^T A r_2$ and $p_1^T A P_0$ as an example. Every expression becomes a node in a graph. Every possible replacement of a common subexpression is represented by an edge. The resulting graph is shown in Figure \ref{fig:CSE}.
The weight of each edge would be the computational cost of the expression that is replaced. The problem then is to find a set of edges, such that each node is attached to at most one of the selected edges. At the same time, the sum of the weights should be minimized. The requirement that each node is attached to at most one of the selected edges represents the fact that multiple replacements are not possible because the expressions overlap.

Clearly, the actual problem is more complex. A common subexpression might be replaced in more than two expressions. A graph representing this is a hypergraph. In addition to that, in longer expressions, some subexpressions do not overlap, so they can both be replaced. 

In practice, trying to solve this problem is probably not necessary. The expressions encountered in most iterative methods are rarely as complex as in the example above. Furthermore, simple heuristics may already produce good results. In products of more than two quantities, matrix-vector products should always be computed first. Then, if one subexpression is replaced, all other occurrences of this expression should be replaced as well.

Alternatively, to generate as many variants as possible, all possible replacements could be constructed.

%
%
%

\section{Example: Nonsymmetric CG}
\label{sec:derivationOfAlgorithmsExample}

In this section, as a more elaborate example, we will show how the derivation of an algorithm for nonsymmetric CG proceeds.

\subsubsection{PME Generation}

The derivation of some properties for nonsymmetric CG was already shown in Section \ref{sec:derivationOfPropertiesExample} and will not be repeated here. The postcondition is shown below:
\begin{align*}
A P D &= R \left( \underline{I} - \underline{J}  \right)\\
P \left( I - U \right) &= \underline{R} \\
P D &= X \left( \underline{I} - \underline{J} \right)\\
\| R e_r^T \| &< \varepsilon
\end{align*}
The initial partitioning is determined based on the properties of the operands. It is applied to the postcondition as well as all derived properties. The expression that is obtained from the postcondition is flattened, yielding
\begin{gather*}
\myFlaOneByTwo{A P_L D_{TL} = R_L \left( I - J \right) - r_R e_r^T}{ A p_R \delta_{BR} = r_R - r_+ } \\
\myFlaOneByTwo{P_L \left( I - U_{TL} \right) = R_L}{ - P_L u_{TR} + p_R = r_R} \\
\myFlaOneByTwo{ P_L D_{TL} = X_L \left( I - J \right) - x_R e_r^T}{ p_R \delta_{BR} = x_R - x_+ } \text{.}
\end{gather*}
The left-hand side is now matched by the CG function, so
\begin{gather*}
\left\{ \myFlaOneByTwo{R_L}{r_R}, U_{TL}, P_L, D_{TL}, \myFlaOneByTwo{X_L}{x_R} \right\} :=  \text{CG} \left(A, R_L e_0, X_L e_0 \right)
\end{gather*}
is obtained. Thus, all quantities on the left-hand side of that assignment are considered to be known. To find assignments for the remaining unknown quantities, the three equations on the right have to be solved. One of the derived properties is that $P^T A P$ is lower triangular, so $P_L^T A P_L$ is lower triangular as well and $P_L^T A p_R$ is zero. The system would recognize that both $P_L$ and $p_R$ appear in $- P_L u_{TR} + p_R = r_R$ on the left-hand side of products. Hence, $P_L^T A$ is multiplied from the left to both sides of the equation to recreate those properties. The resulting equation $- P_L^T A P_L u_{TR} = P_L^T A r_R$ contains only one unknown quantity, so it is solvable. Since $P_L^T A P_L$ is lower triangular, a triangular system is identified, which can also be written as 
$$u_{TR} := - \left( P_L^T A P_L \right)^{-1} P_L^T A r_R \text{.}$$
Having found an assignment for $u_{TR}$, it is considered known as well. By instead using that $R^T A P$ is lower triangular, a different equation for $u_{TR}$ would have been found, resulting in a different algorithm. In practice, two separate derivation processes would be executed for both variants; here we continue just with the first one.

Now, there is only one unknown quantity left in $- P_L u_{TR} + p_R = r_R$, so the following formula is determined for $p_R$:
$$p_R := r_R + P_L u_{TR}$$
Because $r_R^T r_+$ is zero, $r_R^T$ is multiplied from the left to both sides of $A p_R \delta_{BR} = r_R - r_+$. There is only one unknown in the resulting equation $r_R^T A p_R \delta_{BR} = r_R^T r_R$, so another assignment is found:
$$\delta_{BR} := \frac{r_R^T r_R}{r_R^T A p_R}$$
Alternatively, the fact that $P^T R$ is lower triangular and rectangular could be used. Finally, $A p_R \delta_{BR} = r_R - r_+$ and $p_R \delta_{BR} = x_R - x_+$ can be solved to $r_+$ and $x_+$, respectively, completing the PME:
\begin{align*}
\left\{ \myFlaOneByTwo{R_L}{r_R}, U_{TL}, P_L, D_{TL}, \myFlaOneByTwo{X_L}{x_R} \right\} &:= \text{CG} \left(A, R_L e_0, X_L e_0 \right) \\
u_{TR} &:= - \left( P_L^T A P_L \right)^{-1} P_L^T A r_R \\
p_R &:= r_R + P_L u_{TR} \\
\delta_{BR} &:= \frac{r_R^T r_R}{r_R^T A p_R} \\
r_+ &:= r_R - A p_R \delta_{BR} \\
x_+ &:= x_R - p_R \delta_{BR}
\end{align*}
\subsubsection{Loop Invariant Identification}

\begin{figure}[]
\centering
\begin{center}
\begin{tikzpicture}

\node[default]	(R1)				{1};

\node[default]	(U1)	[below=of R1, yshift=-0.5cm]	{2};

\node[default]	(P1)	[below=of U1, yshift=-0.5cm]	{3};

\node[default]	(D1)	[below=of P1, yshift=-0.5cm]	{4};

\node[default]	(R2)	[below=of D1, yshift=-0.5cm, xshift=-1cm]	{5};

\node[default]	(X2)	[below=of D1, yshift=-0.5cm, xshift=1cm]	{6};

\path[->]	(R1)		edge							(U1)
		(R1)		edge		[bend right=30]			(P1)
		(R1)		edge		[bend right=30]			(D1)
		(R1)		edge		[bend right=30]			(R2);
		
\path[->]	(U1)		edge							(P1);

\path[->]	(P1)		edge							(D1)
		(P1)		edge		[bend right=20]			(R2)
		(P1)		edge		[bend left=20]			(X2);

\path[->]	(D1)		edge							(R2)
		(D1)		edge							(X2);

\end{tikzpicture}
\caption{Dependency graph for nonsymmetric CG.}
\label{fig:dg:nonSymCG}
\end{center}
\end{figure}
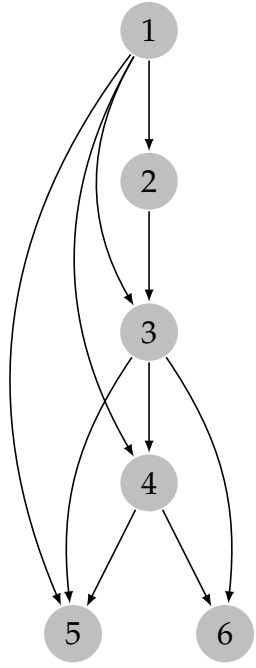

Based on this PME, the dependency graph is constructed. Since there are six assignments in the PME, there are six nodes:
\begin{enumerate}
\item $\left\{ \myFlaOneByTwo{R_L}{r_R}, U_{TL}, P_L, D_{TL}, \myFlaOneByTwo{X_L}{x_R} \right\} := \text{CG} \left(A, R_L e_0, X_L e_0 \right)$
\item $u_{TR} := - \left( P_L^T A P_L \right)^{-1} P_L^T A r_R$
\item $p_R := r_R + P_L u_{TR}$
\item $\delta_{BR} := \frac{r_R^T r_R}{r_R^T A p_R}$
\item $r_+ := r_R - A p_R \delta_{BR}$
\item $x_+ := x_R - p_R \delta_{BR}$
\end{enumerate}
The corresponding dependency graph is shown in Figure \ref{fig:dg:nonSymCG}.
All nonempty subsets of this graph that respect the dependencies are feasible loop invariants. Thus, the following seven loop invariants are obtained:
$$\{1\}, \{1, 2\}, \{1, 2, 3\}, \{1, 2, 3, 4\}, \{1, 2, 3, 4, 5\}, \{1, 2, 3, 4, 6\}, \{1, 2, 3, 4, 5, 6\}$$
%
%
To determine the loop guard $G$, the additional predicate in the postcondition is inspected. It is $\| R e_r^T \| < \varepsilon$. Since the precondition contains $\text{\ttfamily FirstColumnInput}[R]$, according to Table \ref{tab:loopGuards}, the loop guard $\| r_R \| \geq \varepsilon$ is selected.

In this example, we show the derivation for the set $\{1, 2, 3\}$, so the loop invariant is
\begin{align*}
P_\text{inv} = \big\{ &\left\{ \myFlaOneByTwo{R_L}{r_R}, U_{TL}, P_L, D_{TL}, \myFlaOneByTwo{X_L} {x_R} \right\} := \text{CG} \left(A, R_L e_0, X_L e_0 \right) \land \\
&\:u_{TR} := - \left( P_L^T A P_L \right)^{-1} P_L^T A r_R \land \\
&\:p_R := r_R + P_L u_{TR} \big\}\text{.}
\end{align*}
\subsubsection{Algorithm Construction}
In a first step, the preprocessing operations are determined. The relevant assignments are the ones contained in the loop invariant, except for the one obtained from the first node. They are shown below:
\begin{align*}
u_{TR} &:= - \left( P_L^T A P_L \right)^{-1} P_L^T A r_R \\
p_R &:= r_R + P_L u_{TR}
\end{align*}
According to the initial partitioning, $P_L$ is empty, that is, it has the size $n \times 0$. Consequently, the right-hand side of the first assignment is empty, too, so it disappears. The second assignment reduces to
$$p_R := r_R \text{,}$$
which is the only preprocessing operation. Next, the update is derived. The ``Repartition'' rules are determined using Table \ref{tab:repartitionIterativeMethods}:
\begin{align*}
\myFlaOneByThreeI{R_L}{r_R}{r_+} &\rightarrow \myFlaOneByThreeI{R_0}{r_1}{r_2} &
\myFlaOneByThreeI{X_L}{x_R}{x_+} &\rightarrow \myFlaOneByThreeI{X_0}{x_1}{x_2}\\
\myFlaTwoByTwoI{U_{TL}}{u_{TR}}{0}{0} &\rightarrow \myFlaTwoByTwoI{U_{00}}{u_{01}}{0}{0} &
\myFlaTwoByTwoI{D_{TL}}{0}{0}{\delta_{BR}} &\rightarrow \myFlaTwoByTwoI{D_{00}}{0}{0}{\delta_{11}} \\
\myFlaOneByTwoI{P_L}{p_R} &\rightarrow \myFlaOneByTwoI{P_0}{p_1} &&
\end{align*}
Applying those rules to the loop invariant yields the following predicate $P_\text{before}$:
\begin{align*}
P_\text{before} = \big\{ &\left\{ \myFlaOneByTwo{R_0}{r_1}, U_{00}, P_0, D_{00}, \myFlaOneByTwo{X_0} {x_1} \right\} := \text{CG} \left(A, R_0 e_0, X_0 e_0 \right) \land \\
&\:u_{01} := - \left( P_0^T A P_0 \right)^{-1} P_0^T A r_1 \land \\
&\:p_1 := r_1 + P_0 u_{01} \big\}
\end{align*}
The ``Continue with'' repartitioning is
\begin{align*}
\myFlaOneByThreeI{R_L}{r_R}{r_+} &\leftarrow \myFlaOneByFourTFF{R_0}{r_1}{r_2}{r_3} \\
\myFlaOneByThreeI{X_L}{x_R}{x_+} &\leftarrow \myFlaOneByFourTFF{X_0}{x_1}{x_2}{x_3} \\
\myFlaTwoByTwoI{U_{TL}}{u_{TR}}{0}{0} &\leftarrow
\myFlaThreeByThreeTLI	{U_{00}}	{u_{01}}	{u_{02}}
					{0}		{0}		{\nu_{12}}
					{0}		{0}		{0} \\
\myFlaTwoByTwoI{D_{TL}}{0}{0}{\delta_{BR}} &\leftarrow
\myFlaThreeByThreeTLI	{D_{00}}	{0}			{0}
					{0}		{\delta_{11}}	{0}
					{0}		{0}			{\delta_{22}}\\
\myFlaOneByTwoI{P_L}{p_R} &\leftarrow \myFlaOneByThreeLI{P_0}{p_1}{p_2} \text{.}
\end{align*}
Plugging that into the function, the following expression is obtained:
\begin{align*}
&\left\{ \myFlaOneByThree{R_0}{r_1}{r_2}, \myFlaTwoByTwoI{U_{00}}{u_{01}}{0}{0}, \myFlaOneByTwo{P_0}{p_1}, \right. \\&
\qquad \qquad \qquad \; \left. \myFlaTwoByTwoI{D_{00}}{0}{0}{\delta_{11}}, \myFlaOneByThree{X_0} {x_1}{x_2} \right\} := \text{CG} \left(A, R_0 e_0, X_0 e_0 \right)
\end{align*}
It it flattened by using the PME, resulting in six assignments:
\begin{align*}
\left\{ \myFlaOneByTwo{R_0}{r_1}, U_{00}, P_0, D_{00}, \myFlaOneByTwo{X_0}{x_1} \right\} &:= \text{CG} \left(A, R_0 e_0, X_0 e_0 \right) \\
u_{01} &:= - \left( P_0^T A P_0 \right)^{-1} P_0^T A r_1 \\
p_1 &:= r_1 + P_0 u_{01} \\
\delta_{11} &:= \frac{r_1^T r_1}{r_1^T A p_1} \\
r_2 &:= r_1 - A p_1 \delta_{11} \\
x_2 &:= x_1 - p_1 \delta_{11}
\end{align*}
For the second assignment of the loop invariant, the PME of a lower triangular system is needed (see Section \ref{sec:triLS}).
\begin{align*}
&\myFlaTwoByOne{u_{02}}{\nu_{12}}:= - \myFlaTwoByTwo{P_0^T A P_0}{0}{p_1^T A P_0}{p_1^T A p_1}^{-1} \myFlaTwoByOne{P_0^T A r_2}{p_1^T A r_2} \\
&\Rightarrow \left\{
\begin{aligned}
u_{02} &:= \left(- P_0^T A P_0 \right)^{-1} \cdot P_0^T A r_2 \\
\nu_{12} &:= - \frac{p_1^T A r_2 + p_1^T A P_0 u_{02}}{p_1^T A p_1}
\end{aligned} \right.
\end{align*}
The third equation, $p_1 := r_1 + P_0 u_{01}$, becomes $p_2 := r_2 + P_0 u_{02} + p_1 \nu_{12}$ after the application of the ``Continue with'' partitioning. Now that the complete predicate $P_\text{after}$ is determined, the update is found by comparing it to $P_\text{before}$. The assignments that are contained in $P_\text{after}$, but not in $P_\text{before}$, constitute the update. They are shown below.
\begin{align*}
\delta_{11} &:= \frac{r_1^T r_1}{r_1^T A p_1} \\
r_2 &:= r_1 - A p_1 \delta_{11} \\
x_2 &:= x_1 - p_1 \delta_{11} \\
u_{02} &:= \left(- P_0^T A P_0 \right)^{-1} \cdot P_0^T A r_2 \\
\nu_{12} &:= - \frac{p_1^T A r_2 + p_1^T A P_0 u_{02}}{p_1^T A p_1} \\
p_2 &:= r_2 + P_0 u_{02} + p_1 \nu_{12}
\end{align*}

\subsubsection{Postprocessing}

For this example, we use the heuristics explained in Section \ref{sec:postprocessing} for the elimination of common subexpressions. The first expression that is identified is $A p_1$. The auxiliary variable $t_1 := A p_1$ is introduced and all occurrences of $A p_1$ are replaced with $t_1$. The resulting assignments are shown below.
\begin{align*}
t_1 &:= A p_1 \\
\delta_{11} &:= \frac{r_1^T r_1}{r_1^T t_1} \\
r_2 &:= r_1 - t_1 \delta_{11} \\
x_2 &:= x_1 - p_1 \delta_{11} \\
u_{02} &:= \left(- P_0^T A P_0 \right)^{-1} \cdot P_0^T A r_2 \\
\nu_{12} &:= - \frac{p_1^T A r_2 + p_1^T A P_0 u_{02}}{p_1^T t_1} \\
p_2 &:= r_2 + P_0 u_{02} + p_1 \nu_{12}
\end{align*}
Further auxiliary variables are introduced for $A r_2$ and $P_0 u_{02}$. The update that is obtained at the end of this step is shown in the filled out worksheet in Figure \ref{fig:ws:nonsymCG}. Some of the partitionings and the loop invariant are omitted in the interest of legibility.
\begin{figure}
\centering
\begin{minipage}[t]{4.5in}
	\resetsteps
	\setboolean{BlockedAlgQ}{false}
	
	\renewcommand{\ALGroutinename}{ nonsymmetric CG}
	
	\renewcommand{\WSprecondition}{}
	
	\renewcommand{\WSpartition}{
        $
        R \rightarrow
		\myFlaOneByThreeI{R_L}{r_R}{r_+}
	$
        }

	\renewcommand{\WSpartitionsizes}{
		$ R_{L} $ is $ n \times 0 $
	}
	
	\renewcommand{\WSpreprocessing}{
		$p_R := r_R$
	}
	
	\renewcommand{\WSguard}{ $\| r_R \| \geq \varepsilon$ }

	\renewcommand{\WSrepartition}{
	$\myFlaOneByThreeI{R_L}{r_R}{r_+} \rightarrow \myFlaOneByThreeI{R_0}{r_1}{r_2}$
	}
	
	\renewcommand{\WSrepartitionsizes}{
	$A$ is $k \times k$
	}

	\renewcommand{\WSupdate}{
		$\begin{aligned}
		t_1 &:= A p_1 \\
		\delta_{11} &:= \frac{r_1^T r_1}{r_1^T t_1} \\
		r_2 &:= r_1 - t_1 \delta_{11} \\
		x_2 &:= x_1 - p_1 \delta_{11} \\
		t_2 &:= A r_2 \\
		u_{02} &:= \left(- P_0^T A P_0 \right)^{-1} \cdot P_0^T t_2 \\
		t_3 &:= P_0 u_{02} \\
		\nu_{12} &:= - \frac{p_1^T t_2 + p_1^T A t_3}{p_1^T t_1} \\
		p_2 &:= r_2 + t_3 + p_1 \nu_{12}
		\end{aligned}$
	}
	
	\renewcommand{\WSmoveboundary}{
	$\myFlaOneByThreeI{R_L}{r_R}{r_+} \leftarrow \myFlaOneByFourTFF{R_0}{r_1}{r_2}{r_3}$
	}
	
	\renewcommand{\WSpostcondition}{
	
	}
	
	{
	\FlaAlgorithmIter
	}
\end{minipage}
\caption{Worksheet for a nonsymmetric CG algorithm.}
\label{fig:ws:nonsymCG}
\end{figure}
\section{Example: BiCG}
\label{sec:BiCGexample}

As a second example, we show the derivation of two algorithms for the biconjugate gradient method (BiCG). Since we already showed a full example in the previous section, we now proceed at a slightly higher pace. Pre- and postcondition are shown below.
\begin{align*}
P_\text{pre}: \{ &\text{\ttfamily Input}[A] \land \text{\ttfamily Matrix}[A] \land \text{\ttfamily NonSingular}[A] \land \\
		&\text{\ttfamily Output}[P] \land \text{\ttfamily Matrix}[P] \land \\
		&\text{\ttfamily Output}[\tilde{P}] \land \text{\ttfamily Matrix}[\tilde{P}] \land \\
		&\text{\ttfamily Output}[D] \land \text{\ttfamily Matrix}[D] \land \text{\ttfamily Diagonal}[D] \land \\
		&\text{\ttfamily FirstColumnInput}[R] \land \text{\ttfamily Matrix}[R] \land \\
		&\text{\ttfamily FirstColumnInput}[\underline{R}] \land \text{\ttfamily Matrix}[\underline{R}]  \land \\
		&\text{\ttfamily FirstColumnInput}[\tilde{R}] \land \text{\ttfamily Matrix}[\tilde{R}] \land  \\
		&\text{\ttfamily FirstColumnInput}[\underline{\tilde{R}}] \land \text{\ttfamily Matrix}[\underline{\tilde{R}}] \land \\
		&\text{\ttfamily Diagonal}[R^T \tilde{R}] \land \text{\ttfamily Diagonal}[\underline{R}^T \underline{\tilde{R}}] \land \\
		&\text{\ttfamily DiagonalR}[\underline{R}^T \tilde{R}] \land \text{\ttfamily DiagonalR}[R^T \underline{\tilde{R}}] \land \\
		&\text{\ttfamily FirstColumnInput}[X] \land \text{\ttfamily Matrix}[X] \land \\
		&\text{\ttfamily Output}[U] \land \text{\ttfamily Matrix}[U] \land \text{\ttfamily UpperDiagonal}[U] \land \\
		& \text{\ttfamily Matrix}[\underline{I} - \underline{J}] \land \text{\ttfamily LowerTrapezoidal}[\underline{I} - \underline{J}] \}
\end{align*}
\begin{align*}
P_\text{post}:	\{ 	&APD = R \left( \underline{I} - \underline{J} \right) \\
				&A^T \tilde{P}D = \tilde{R} \left( \underline{I} - \underline{J} \right) \\
				&P \left( I - U \right) = \underline{R} \\
				&\tilde{P} \left( I - U \right) = \underline{\tilde{R}} \\
				&PD = X \left( \underline{I} - \underline{J} \right) \\
				&\| R e_r^T \| < \varepsilon \}
\end{align*}
Following from the precondition, the function representing the operation is
\begin{align*}
\left\{ R, \tilde{R}, U, P, \tilde{P}, D, X \right\} &:= \text{BiCG} \left(A, R e_0, X e_0 \right) \text{.}
\end{align*}

\subsubsection{Derivation of Properties}

For BiCG, neither $R$ nor $\tilde{R}$ is orthogonal. Instead, they are mutually orthogonal, which means that $R^T \tilde{R}$ is diagonal. Except for that, the derivation of properties is mostly the same to the one shown in the example in Section \ref{sec:derivationOfPropertiesExample}. For this reason, it will not be shown here. Instead. we give a short overview of the relevant properties: 
\begin{itemize}
\item[-] $\underline{\tilde{R}}^T A P$ and $\underline{R}^T A^T \tilde{P}$ are lower triangular.
\item[-] $\tilde{P}^T A P$ (and thus $P^T A^T \tilde{P}$) is diagonal.
\item[-] $P^T \tilde{R}$ and $\tilde{P}^T R$ are lower triangular and rectangular.
\end{itemize}
%
%
%
\subsubsection{PME Generation}

As a first step towards the PME, the operands are partitioned. Initially, only the fist columns of $R$ and $\tilde{R}$ are known, so they are partitioned into \smash{$\myFlaOneByThree{R_L}{r_R}{r_+}$} and \smash{$\myFlaOneByThree{\tilde{R}_L}{\tilde{r}_R}{\tilde{r}_+}$}, respectively. The remaining operands are partitioned accordingly:
\begin{align*}
A \myFlaOneByTwo{P_L}{p_R} \myFlaTwoByTwo{D_{TL}}{0}{0}{\delta_{BR}} &= \myFlaOneByThree{R_L}{r_R}{r_+} \myFlaThreeByTwo{I - J}{0}{-e_r^T}{1}{0}{-1} \\
A^T \myFlaOneByTwo{\tilde{P}_L}{\tilde{p}_R} \myFlaTwoByTwo{D_{TL}}{0}{0}{\delta_{BR}} &= \myFlaOneByThree{\tilde{R}_L}{\tilde{r}_R}{\tilde{r}_+} \myFlaThreeByTwo{I - J}{0}{-e_r^T}{1}{0}{-1} \\
\myFlaOneByTwo{P_L}{p_R}\myFlaTwoByTwo{I- U_{TL}}{- u_{TR}}{0}{1} &= \myFlaOneByTwo{R_L}{r_R} \\
\myFlaOneByTwo{\tilde{P}_L}{\tilde{p}_R}\myFlaTwoByTwo{I- U_{TL}}{- u_{TR}}{0}{1} &= \myFlaOneByTwo{\tilde{R}_L}{\tilde{r}_R} \\
\myFlaOneByTwo{P_L}{p_R} \myFlaTwoByTwo{D_{TL}}{0}{0}{\delta_{BR}} &= \myFlaOneByThree{X_L}{x_R}{x_+} \myFlaThreeByTwo{I - J}{0}{-e_r^T}{1}{0}{-1}
\end{align*}
After the execution of the Matrix Arithmetic step, those expressions become
\begin{gather*}
\myFlaOneByTwo{A P_L D_{TL} = R_L \left( I - J \right) - r_R e_r^T}{ A p_R \delta_{BR} = r_R - r_+ } \\
\myFlaOneByTwo{A^T \tilde{P}_L D_{TL} = \tilde{R}_L \left( I - J \right) - \tilde{r}_R e_r^T}{ A \tilde{p}_R \delta_{BR} = \tilde{r}_R - \tilde{r}_+ } \\
\myFlaOneByTwo{P_L \left( I - U_{TL} \right) = R_L}{ - P_L u_{TR} + p_R = r_R} \\
\myFlaOneByTwo{\tilde{P}_L \left( I - U_{TL} \right) = \tilde{R}_L}{ - \tilde{P}_L u_{TR} + \tilde{p}_R = \tilde{r}_R} \\
\myFlaOneByTwo{ P_L D_{TL} = X_L \left( I - J \right) - x_R e_r^T}{ p_R \delta_{BR} = x_R - x_+ }\text{.}
\end{gather*}
The left-hand sides of those expressions can now be replaced by the BiCG function:
\begin{align*}
\left\{ \myFlaOneByTwo{R_L}{r_R}, \myFlaOneByTwo{\tilde{R}_L}{\tilde{R}_R}, U_{TL}, P_L, \tilde{P}_L, D_{TL}, \myFlaOneByTwo{X_L}{x_R} \right\} &:= \text{BiCG} \left(A, R_L e_0, X_L e_0 \right)
\end{align*}
Updates for the remaining unknown quantities are obtained by solving the equations on the right-hand using the derived properties. For this example, to find an assignment for $u_{TR}$, we use that $\tilde{P}^T A P$ is diagonal. Multiplying $\tilde{P}_L^T A$ from the left to both sides of $- P_L u_{TR} + p_R = r_R$ results in $- \tilde{P}_L^T A P_L u_{TR} = \tilde{P}_L^T A r_R$. Since $\tilde{P}_L^T A P_L$ is diagonal, and thus is invertible, the following assignment is obtained:
$$u_{TR} := - \left( \tilde{P}_L^T A P_L \right)^{-1} \tilde{P}_L^T A r_R$$
Note that it would have been also possible to use $- \tilde{P}_L u_{TR} + \tilde{p}_R = \tilde{r}_R$ to get to an assignment for $u_{TR}$.  Now, the following assignments are obtained for $p_R$ and $\tilde{p}_R$:
\begin{align*}
p_R &:= r_R +  P_L u_{TR} \\
\tilde{p}_R &:= \tilde{r}_R + \tilde{P}_L u_{TR}
\end{align*}
Similarly to $u_{TR}$, $\delta_{BR}$ can be computed in several different ways. Here, we use that $\tilde{P}^T R$ is lower triangular and rectangular, in combination with the equation $A p_R \delta_{BR} = r_R - r_+$. Multiplying $\tilde{p}_R^T$ from the left and solving to $\delta_{BR}$ yields
$$\delta_{BR} := \frac{\tilde{p}_R^T r_R}{\tilde{p}_R^T A p_R}$$
Finally, the following assignments are derived for $r_+$, $\tilde{r}_+$ and $x_+$:
\begin{align*}
r_+ &:= r_R -  A p_R \delta_{BR} \\
\tilde{r}_+ &:= \tilde{r}_R - A \tilde{p}_R \delta_{BR} \\
x_+ & := x_R -  p_R \delta_{BR}
\end{align*}
The complete PME is shown below (already in form of a list, in anticipation of the construction of the dependency graph):
%

\begin{enumerate}
\item $\left\{ \myFlaOneByTwo{R_L}{r_R}, \myFlaOneByTwo{\tilde{R}_L}{\tilde{R}_R}, U_{TL}, P_L, \tilde{P}_L, D_{TL}, \myFlaOneByTwo{X_L}{x_R} \right\} := \text{BiCG} \left(A, R_L e_0, X_L e_0 \right)$
\item $u_{TR} := - \left( \tilde{P}_L^T A P_L \right)^{-1} \tilde{P}_L^T A r_R$
\item $p_R := r_R +  P_L u_{TR}$
\item $\tilde{p}_R := \tilde{r}_R + \tilde{P}_L u_{TR}$
\item $\delta_{BR} := \frac{\tilde{p}_R^T r_R}{\tilde{p}_R^T A p_R}$
\item $r_+ := r_R -  A p_R \delta_{BR}$
\item $\tilde{r}_+ := \tilde{r}_R - A \tilde{p}_R \delta_{BR}$
\item $x_+  := x_R -  p_R \delta_{BR}$
\end{enumerate}

\subsubsection{Loop Invariant Identification}

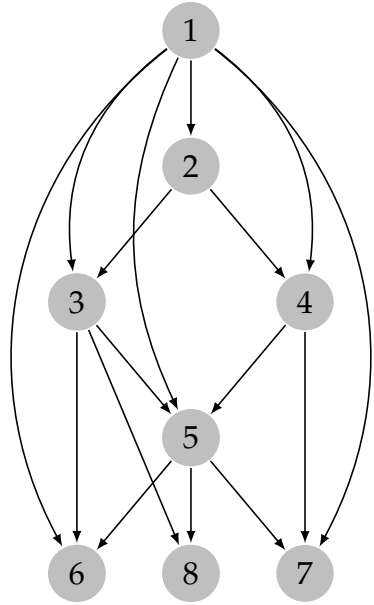
\begin{figure}[]
\centering
\begin{center}
\begin{tikzpicture}

\node[default]	(R1)				{1};

\node[default]	(U1)	[below=of R1, yshift=-0.3cm]	{2};

\node[default]	(P1)	[below left=of U1, yshift=-0.3cm]	{3};

\node[default]	(PT1)	[below right=of U1, yshift=-0.3cm]	{4};

\node[default]	(D1)	[below right=of P1, yshift=-0.3cm]	{5};

\node[default]	(R2)	[below left=of D1, yshift=-0.3cm, xshift=-0cm]	{6};

\node[default]	(RT2)	[below right=of D1, yshift=-0.3cm, xshift=0cm]	{7};

\node[default]	(X2)	[below=of D1, yshift=-0.3cm, xshift=0cm]	{8};

\path[->]	(R1)		edge							(U1)
		(R1)		edge		[bend right=30]			(P1)
		(R1)		edge		[bend left=30]			(PT1)
		(R1)		edge		[bend right=25]			(D1)
		(R1)		edge		[bend right=40]			(R2)
		(R1)		edge		[bend left=40]			(RT2);
		
\path[->]	(U1)		edge							(P1)
		(U1)		edge							(PT1);

\path[->]	(P1)		edge							(D1)
		(P1)		edge		[bend left=0]			(R2)
		(P1)		edge		[bend right=0]			(X2);

\path[->]	(PT1)	edge							(D1)
		(PT1)	edge							(RT2);

\path[->]	(D1)		edge							(R2)
		(D1)		edge							(RT2)
		(D1)		edge							(X2);

\end{tikzpicture}
\caption{Dependency graph for BiCG.}
\label{fig:dg:BiCG}
\end{center}
\end{figure}

The dependency graph is shown in Figure \ref{fig:dg:BiCG}.
%
%
%
%
%
%
%
%
%
%
%
%
%
%
%
%
%
%
Of all the subsets of this graph that respect the dependencies, only the empty one fails to satisfy the conditions for the feasibility of loop invariant. The remaining 13 subsets are feasible loop invariants. We will continue the derivation in this example with the loop invariant that corresponds to the full set.
%

Because of the additional predicate in the postcondition, $\| R e_r^T \| < \varepsilon$, and the property $\text{\ttfamily FirstColumnInput}[R]$ in the precondition, the loop guard $\| r_R \| < \varepsilon$ is determined.

\subsubsection{Algorithm Construction}

The first part of this step consists of finding the preprocessing operations. This is done by taking all assignments of the loop invariant except for the one from the first node, and eliminating all expression that are empty in the initial partitioning. The relevant assignments are shown below:
\begin{align*}
u_{TR} &:= - \left( \tilde{P}_L^T A P_L \right)^{-1} \tilde{P}_L^T A r_R \\
p_R &:= r_R +  P_L u_{TR} \\
\tilde{p}_R &:= \tilde{r}_R + \tilde{P}_L u_{TR} \\
\delta_{BR} &:= \frac{\tilde{p}_R^T r_R}{\tilde{p}_R^T A p_R}
\end{align*}
\begin{align*}
r_+ &:= r_R -  A p_R \delta_{BR} \\
\tilde{r}_+ &:= \tilde{r}_R - A \tilde{p}_R \delta_{BR} \\
x_+ & := x_R -  p_R \delta_{BR}
\end{align*}
Initially, $P_L$ and $\tilde{P}_L$ have the size $n \times 0$ and $u_{TR}$ is of size $0 \times 1$. As a result, the first assignment is eliminated, while the second and third are reduced to
\begin{align*}
p_R &:= r_R \\
\tilde{p}_R &:= \tilde{r}_R \text{.}
\end{align*}
The remaining four assignments do not change.

To determine the update, the predicates $P_\text{before}$ and $P_\text{after}$ have to be constructed by repartitioning the loop invariant. The ``Repartition'' rules are shown below:
\begin{align*}
\myFlaOneByThreeI{R_L}{r_R}{r_+} &\rightarrow \myFlaOneByThreeI{R_0}{r_1}{r_2} &
\myFlaOneByThreeI{\tilde{R}_L}{\tilde{r}_R}{\tilde{r}_+} &\rightarrow \myFlaOneByThreeI{\tilde{R}_0}{\tilde{r}_1}{\tilde{r}_2}\\
\myFlaTwoByTwoI{U_{TL}}{u_{TR}}{0}{0} &\rightarrow \myFlaTwoByTwoI{U_{00}}{u_{01}}{0}{0} &
\myFlaTwoByTwoI{D_{TL}}{0}{0}{\delta_{BR}} &\rightarrow \myFlaTwoByTwoI{D_{00}}{0}{0}{\delta_{11}} \\
\myFlaOneByTwoI{P_L}{p_R} &\rightarrow \myFlaOneByTwoI{P_0}{p_1} & \myFlaOneByTwoI{\tilde{P}_L}{\tilde{p}_R} &\rightarrow \myFlaOneByTwoI{\tilde{P}_0}{\tilde{p}_1} \\ 
\myFlaOneByThreeI{X_L}{x_R}{x_+} &\rightarrow \myFlaOneByThreeI{X_0}{x_1}{x_2} &&
\end{align*}
Applying that to the loop invariant yields the following equations which constitute $P_\text{before}$:
\begin{align*}
\left\{ \myFlaOneByTwo{R_0}{r_1}, \myFlaOneByTwo{\tilde{R}_0}{\tilde{R}_1}, U_{00}, P_0, \tilde{P}_0, D_{00}, \myFlaOneByTwo{X_0}{x_1} \right\} &:= \text{BiCG} \left(A, R_0 e_0, X_0 e_0 \right) \\
u_{01} &:= - \left( \tilde{P}_0^T A P_0 \right)^{-1} \tilde{P}_0^T A r_1 \\
p_1 &:= r_1 +  P_0 u_{01} \\
\tilde{p}_1 &:= \tilde{r}_1 + \tilde{P}_0 u_{01} \\
\delta_{11} &:= \frac{\tilde{p}_1^T r_1}{\tilde{p}_1^T A p_1} \\
r_2 &:= r_1 -  A p_1 \delta_{11} \\
\tilde{r}_2 &:= \tilde{r}_1 - A \tilde{p}_1 \delta_{11} \\
x_2 & := x_1 -  p_1 \delta_{11} 
\end{align*}
For the ``Continue with'' repartitioning it is important to note that $U$ is upper diagonal. Thus, $u_{02}$ is zero:
\begin{align*}
\myFlaOneByThreeI{R_L}{r_R}{r_+} &\leftarrow \myFlaOneByFourTFF{R_0}{r_1}{r_2}{r_3} \\
\myFlaOneByThreeI{\tilde{R}_L}{\tilde{r}_R}{\tilde{r}_+} &\leftarrow \myFlaOneByFourTFF{\tilde{R}_0}{\tilde{r}_1}{\tilde{r}_2}{\tilde{r}_3} \\
\myFlaTwoByTwoI{U_{TL}}{u_{TR}}{0}{0} &\leftarrow
\myFlaThreeByThreeTLI	{U_{00}}	{u_{01}}	{0}
					{0}		{0}		{\nu_{12}}
					{0}		{0}		{0} \\
\myFlaTwoByTwoI{D_{TL}}{0}{0}{\delta_{BR}} &\leftarrow
\myFlaThreeByThreeTLI	{D_{00}}	{0}			{0}
					{0}		{\delta_{11}}	{0}
					{0}		{0}			{\delta_{22}}\\
\myFlaOneByTwoI{P_L}{p_R} &\leftarrow \myFlaOneByThreeLI{P_0}{p_1}{p_2} \\
\myFlaOneByTwoI{\tilde{P}_L}{\tilde{p}_R} &\leftarrow \myFlaOneByThreeLI{\tilde{P}_0}{\tilde{p}_1}{\tilde{p}_2} \\
\myFlaOneByThreeI{X_L}{x_R}{x_+} &\leftarrow \myFlaOneByFourTFF{X_0}{x_1}{x_2}{x_3} \\
\end{align*}
This repartitioning transforms
$$\left\{ \myFlaOneByTwo{R_L}{r_R}, \myFlaOneByTwo{\tilde{R}_L}{\tilde{R}_R}, U_{TL}, P_L, \tilde{P}_L, D_{TL}, \myFlaOneByTwo{X_L}{x_R} \right\} := \text{BiCG} \left(A, R_L e_0, X_L e_0 \right)$$
into those assignments that are also contained in the predicate $P_\text{before}$. Applying it to 
$$u_{TR} := - \left( \tilde{P}_L^T A P_L \right)^{-1} \tilde{P}_L^T A r_R \text{,}$$
yields
$$\myFlaTwoByOne{0}{\nu_{12}} := \myFlaTwoByTwo{\tilde{P}_0^T A P_0}{0}{0}{\tilde{p}_1^T A p_1}^{-1} \myFlaTwoByOne{\tilde{P}_0^T A r_2}{\tilde{p}_1^T A r_2} \text{.}$$
Flattening this expression results in the assignment
$$\nu_{12} := \frac{\tilde{p}_1^T A r_2}{\tilde{p}_1^T A p_1}$$
for $\nu_{12}$. Because of 
\begin{align*}
p_R &:= r_R +  P_L u_{TR} \\
\tilde{p}_R &:= \tilde{r}_R + \tilde{P}_L u_{TR} \text{,}
\end{align*}
the following two expressions are obtained and added to $P_\text{after}$:
\begin{align*}
p_2 &:= r_2 +  p_1 \nu_{12} \\
\tilde{p}_2 &:= \tilde{r}_2 + \tilde{p}_1 \nu_{12}
\end{align*}
For the remaining assignments of the loop invariant, only the indices change:
\begin{align*}
\delta_{22} &:= \frac{\tilde{p}_2^T r_2}{\tilde{p}_2^T A p_2} \\
r_3 &:= r_2 -  A p_2 \delta_{22} \\
\tilde{r}_3 &:= \tilde{r}_2 - A \tilde{p}_2 \delta_{22} \\
x_3 & := x_2 -  p_2 \delta_{22}
\end{align*}
Now that the $P_\text{before}$ and $P_\text{after}$ are completely determined, the update is found by identifying those assignments that are contained in $P_\text{after}$, but not in $P_\text{before}$. For this example, no common subexpressions are replaced. The worksheet for this algorithm, together with a second one for the algorithm obtained from the set that only contains the first node, is shown in Figure \ref{fig:ws:BiCG}.

\begin{figure}
\centering
\scalebox{0.9}{
\begin{tabular}{cc}
\begin{minipage}[t]{0.5\textwidth}
	\resetsteps
	\setboolean{BlockedAlgQ}{false}
	
	\renewcommand{\ALGroutinename}{ BiCG}
	
	\renewcommand{\WSprecondition}{}
	
	\renewcommand{\WSpartition}{
        $
        R \rightarrow
		\myFlaOneByThreeI{R_L}{r_R}{r_+}
	$
        }

	\renewcommand{\WSpartitionsizes}{
		$ R_{L} $ is $ n \times 0 $
	}
	
	\renewcommand{\WSpreprocessing}{
		$\begin{aligned}
		p_R &:= r_R \\
		\tilde{p}_R &:= \tilde{r}_R \\
		\delta_{BR} &:= \frac{\tilde{p}_R^T r_R}{\tilde{p}_R^T A p_R} \\
		r_+ &:= r_R -  A p_R \delta_{BR} \\
		\tilde{r}_+ &:= \tilde{r}_R - A \tilde{p}_R \delta_{BR} \\
		x_+ & := x_R -  p_R \delta_{BR}
		\end{aligned}$
	}
	
	\renewcommand{\WSguard}{ $\| r_R \| \geq \varepsilon$ }

	\renewcommand{\WSrepartition}{
	$\myFlaOneByThreeI{R_L}{r_R}{r_+} \rightarrow \myFlaOneByThreeI{R_0}{r_1}{r_2}$
	}
	
	\renewcommand{\WSrepartitionsizes}{
	$A$ is $k \times k$
	}

	\renewcommand{\WSupdate}{
		$\begin{aligned}
		\nu_{12} &:= \frac{\tilde{p}_1^T A r_2}{\tilde{p}_1^T A p_1} \\
		p_2 &:= r_2 +  p_1 \nu_{12} \\
		\tilde{p}_2 &:= \tilde{r}_2 + \tilde{p}_1 \nu_{12} \\
		\delta_{22} &:= \frac{\tilde{p}_2^T r_2}{\tilde{p}_2^T A p_2} \\
		r_3 &:= r_2 -  A p_2 \delta_{22} \\
		\tilde{r}_3 &:= \tilde{r}_2 - A \tilde{p}_2 \delta_{22} \\
		x_3 & := x_2 -  p_2 \delta_{22}
		\end{aligned}$
	}
	
	\renewcommand{\WSmoveboundary}{
	$\myFlaOneByThreeI{R_L}{r_R}{r_+} \leftarrow \myFlaOneByFourTFF{R_0}{r_1}{r_2}{r_3}$
	}
	
	\renewcommand{\WSpostcondition}{
	
	}
	
	{
	\FlaAlgorithmIter
	}
\end{minipage}
\vspace{0.05cm} & \vspace{0.05cm}
\begin{minipage}[t]{0.5\textwidth}
	\resetsteps
	\setboolean{BlockedAlgQ}{false}
	
	\renewcommand{\ALGroutinename}{ BiCG}
	
	\renewcommand{\WSprecondition}{}
	
	\renewcommand{\WSpartition}{
        $
        R \rightarrow
		\myFlaOneByThreeI{R_L}{r_R}{r_+}
	$
        }

	\renewcommand{\WSpartitionsizes}{
		$ R_{L} $ is $ n \times 0 $
	}
	
	\renewcommand{\WSpreprocessing}{
	}
	
	\renewcommand{\WSguard}{ $\| r_R \| \geq \varepsilon$ }

	\renewcommand{\WSrepartition}{
	$\myFlaOneByThreeI{R_L}{r_R}{r_+} \rightarrow \myFlaOneByThreeI{R_0}{r_1}{r_2}$
	}
	
	\renewcommand{\WSrepartitionsizes}{
	$A$ is $k \times k$
	}

	\renewcommand{\WSupdate}{
		$\begin{aligned}
		u_{01} &:= - \left( \tilde{P}_0^T A P_0 \right)^{-1} \tilde{P}_0^T A r_1 \\
		p_1 &:= r_1 +  P_0 u_{01} \\
		\tilde{p}_1 &:= \tilde{r}_1 + \tilde{P}_0 u_{01} \\
		\delta_{11} &:= \frac{\tilde{p}_1^T r_1}{\tilde{p}_1^T A p_1} \\
		r_2 &:= r_1 -  A p_1 \delta_{11} \\
		\tilde{r}_2 &:= \tilde{r}_1 - A \tilde{p}_1 \delta_{11} \\
		x_2 & := x_1 -  p_1 \delta_{11} 
		\end{aligned}$
	}
	
	\renewcommand{\WSmoveboundary}{
	$\myFlaOneByThreeI{R_L}{r_R}{r_+} \leftarrow \myFlaOneByFourTFF{R_0}{r_1}{r_2}{r_3}$
	}
	
	\renewcommand{\WSpostcondition}{
	
	}
	
	{
	\FlaAlgorithmIter
	}
\end{minipage}
\end{tabular}
} 
\caption{Worksheets for two BiCG algorithms. The one on the left is obtained from the loop invariant that consists of the full set of nodes. The right one is derived from the set that only contains the first node.}
\label{fig:ws:BiCG}
\end{figure}

\subsection{Remark on the Equivalence of Loop Invariants and Algorithms}
\label{sec:RemarkEquivalence}

Comparing those algorithms, one notes that the updates are very similar. In fact, according to the criteria established in Section \ref{sec:equivalence}, those algorithms are considered equivalent. With a suitable replacement, the loop invariant for one can be transformed into the loop invariant of the other. The differences in the shape of some assignments stem from the fact that $u_{01}$ is zero except for the last position. In the algorithm on the left in Figure \ref{fig:ws:BiCG}, this is revealed by the repartitioning, while this is not the case on the right.

In general, one observes that the presented approach always produces two equivalent loop invariants, namely the subset of the dependency graph that only contains the first node, and the full set. For direct methods, the coarsest possible, that is, the standard $2 \times 2$ partitioning, never results in equivalent loop invariants. Finer partitionings result in new loop invariants than can not be found with a coarser one.

For the iterative methods covered in this thesis, this is different. The partitioning that is used for the derivation is the coarsest possible one, as a coarser one would only partition those matrices that are initially partially known. All other matrices would not be partitioned at all, so it would not be possible to derive a PME. On the other hand, for quite a few of the methods presented in Appendix \ref{chap:appendixRepresentations}, a finer partitioning does not result in new loop invariants.

This is the case for those methods where either there is no matrix $U$ or it is upper diagonal. The reason is that for those methods, finer partitionings do not partition quantities that are computed in one iteration in multiple parts. All they do is expose additional quantities that are fully computed in different iterations.

%
%
%

\section{Scope and Limitations}
\label{sec:scopeLimitations}

The presented method extends to a lot more iterative methods than those shown as examples throughout the thesis. Matrix representations of further methods are shown in Appendix \ref{chap:appendixRepresentations}. Note that this includes stationary iterative methods (\ref{sec:stationaryIterative}). For those, the derivation is even simpler because no properties have to be derived.

With the presented approach, it is possible to derive a large number of algorithms for most iterative methods. In case of CG for example, by using the derived properties to solve equations in different ways, four PMEs are found. From each PME, seven loop invariants are obtained, resulting in $28$ algorithms. Another considerable factor is added by the elimination of common subexpressions. Similar numbers can be expected for other, comparably complex iterative methods.
%
%
%
%

There are, however, some limitations. They are explained in the following.

%
%

\subsubsection{Rewriting of Updates}

In some cases, the presented method is not able to generate those assignments that are commonly found in literature. This is for example the case for symmetric CG. The derived updates for $\nu_{12}$ always have a shape like this:
\begin{align}
\nu_{12} :=  \frac{p_1^T A r_2}{p_1^T A p_1} \label{eq:nuUpdate}
\end{align}
Usually, the following formula is used \cite{barrett:templates, saad2000iterative, vanderVorst:book}:
$$\nu_{12} := \frac{r_2^T r_2}{r_1^T r_1}$$
The advantage is that the matrix-vector product $p_1^T A$ is eliminated. It is obtained as follows. We begin with rewriting $r_2 = r_1 - A p_1 \delta_{11}$, which is the update for $r_2$, as $A p_1 = (r_1 - r_2) \delta_{11}^{-1}$. Then, both sides of this equation are transposed, resulting in $p_1^T A = \delta_{11}^{-1} (r_1 - r_2)^T$. Now, this equation is used to replace $p_1^T A$ in the numerator of equation (\ref{eq:nuUpdate}):
\begin{align*}
\nu_{12} :=  \frac{\delta_{11}^{-1} (r_1 - r_2)^T r_2}{p_1^T A p_1}
\end{align*}
Because of the orthogonality or $R$, it simplifies to 
\begin{align*}
\nu_{12} :=  \frac{\delta_{11}^{-1} r_2^T r_2}{p_1^T A p_1} \text{.}
\end{align*}
$\delta_{11}$ is then replaced with
$$\delta_{11} = \frac{r_1^T r_1}{p_1^T A p_1} \text{,}$$
finally resulting in
\begin{align*}
\nu_{12} :=  \frac{p_1^T A p_1}{r_1^T r_1} \cdot \frac{r_2^T r_2}{p_1^T A p_1} = \frac{r_2^T r_2}{r_1^T r_1} \text{.}
\end{align*}

By itself, this transformation is not particularly difficult. Since all steps follow well defined algebraic rules, it is not even difficult to design a system that is able to perform the individual steps of this rewriting. The problem is that based on the initial equation (\ref{eq:nuUpdate}), there is no indication that it is possible to reduce the number of matrix-vector products. Even if that is known, there is no indication which steps to perform. Thus, any system that is supposed to rewrite the assignment has to perform some sort of exhaustive search. Unfortunately, since expressions are substituted, the search space is infinite. Heuristics have to be applied to guarantee the termination of this search. However, they must not limit the capability of the system to find such rewritings.

%
%
%

\subsubsection{Normalized Vectors}

Some iterative methods construct a set of orthonormal vectors. With the presented approach, it is not possible to derive algorithms for those methods. Consider the Arnoldi iteration as an example. The matrix representation is shown below. It is based on the description in \cite{saad2000iterative}:
\begin{align*}
\{ Q, H\}:= \text{AI} ( A, Q e_0) \equiv
\left\{
\begin{aligned}
P_\text{pre}: \{ &\text{\ttfamily Input}(A) \land \text{\ttfamily Matrix}(A) \\
		&\text{\ttfamily FirstColumnInput}(Q) \land \text{\ttfamily Matrix}[Q] \land  \\
		&\text{\ttfamily Orthonormal}[Q] \land \\
				&\text{\ttfamily Output}[H] \land \text{\ttfamily Matrix}[H] \land \text{\ttfamily UpperHessenberg}[H]  \} \\
P_\text{post}: \{ &A \underline{Q} = Q H \\
			&\text{size}(Q) = n \times m \}
\end{aligned}
\right.
\end{align*}
The postcondition is repartitioned as follows:
$$A \myFlaOneByTwo{Q_L}{q_R} = \myFlaOneByThree{Q_L}{q_R}{q_{+}} \myFlaThreeByTwo{H_{TL}}{h_{TR}}{h_{ML}}{\eta_{MR}}{0}{\eta_{BR}}$$
After that, the Matrix Arithmetic step yields the following expression:
$$\myFlaOneByTwo{AQ_L = Q_L H_{TL} + q_R h_{ML}}{A q_R = Q_L h_{TR} + q_R \eta_{MR} + q_+ \eta_{BR}}$$
As usual, the left-hand side is matched by the function:
$$\left\{ \myFlaOneByTwo{Q_L}{q_R} , \myFlaTwoByOne{H_{TL}}{h_{ML}} \right\}:= \text{AI} ( A, Q_L e_0)$$
Using that $Q$ is orthonormal, the following assignments are obtained for $h_{TR}$ and $\eta_{MR}$:
\begin{align*}
h_{TR} &:= Q_L^T A q_R \\
\eta_{MR} &:= q_R^T A q_R 
\end{align*}
The problem is now to compute the normalized vector $q_+$ and $\eta_{BR}$. From $A q_R = Q_L h_{TR} + q_R \eta_{MR} + q_+ \eta_{BR}$, the following equation can be obtained, where all quantities on the right-hand side are known:
$$q_+ \eta_{BR} = A q_R - Q_L h_{TR} - q_R \eta_{MR}$$
Since we know that $q_+$ is normalized, the scalar $\eta_{BR}$ has to be computed as
\begin{align}
\eta_{BR} :=  \| A q_R - Q_L h_{TR} - q_R \eta_{MR} \| \text{.}
\label{eq:ArnoldiProblem}
\end{align}
Then, $q_+$ is determined by the following assignment:
$$q_+ := \frac{A q_R - Q_L h_{TR} - q_R \eta_{MR}}{\eta_{BR}}$$
With the presented method of solving equations by applying properties, the assignment (\ref{eq:ArnoldiProblem}) can not be obtained. To ensure that algorithms for such methods can be derived, the presented approach for solving equations must be expanded.


\chapter{Conclusion}
\label{chap:conclusion}

This thesis introduces a methodology that allows the systematic derivation of algorithms for iterative methods; the starting point for this methodology is a formal description of an iterative method in matrix form. In addition, we presented an approach for deriving properties of matrices and matrix expressions from this representation; those properties are necessary for the derivation of algorithms.

The actual derivation of algorithms consists of four major steps. First, PMEs are generated by partitioning the operands of the matrix representation and applying the derived properties to solve equations. Then, from those PMEs, loop invariants are obtained. In the third step, from each loop invariant, one loop-based algorithm is constructed. Finally, common subexpressions are eliminated, generating an even larger number of algorithms.

One of the most important aspects, and indispensable for the automatic generation of libraries, is that the derived algorithms are provably correct. This is ensured by constructing them around a proof of correctness, based on the loop invariants generated in the second step.

A conscious effort was made to ensure that the entire process is systematic, that is, each step is performed according to well defined rules and no guidance by a human expert is required. This allows the approach to be implemented as a tool that automatically generates algorithms based on a formal description of the operation. We consider this to be another important step towards the automatic generation of linear algebra libraries as envisioned by the founders of the FLAME project.

As for future work, there are a number of ways to build on the results of this thesis:

\begin{description}
\item[Implementation] Executing the presented approach by hand is a laborious and thus error-prone task, not least because it was not designed to be executed by hand. To be used productively, the presented approach should be implemented as a computer program.

\item[Stability Analysis] To asses the usefulness of the derived algorithms in practice, a stability analysis is indispensable. In \cite{Bientinesi:thesis}, it was shown that the FLAME methodology can be combined with a systematic stability analysis. The presented method should be extended in a similar way.

\item[Performance Analysis] While it is desirable to derive a large number of algorithms to find new, potentially faster variants, the task of identifying them should not be left to the user. Thus, similar to a systematic stability analysis, the system should also be able to reason about the performance of the generated algorithms and select the best ones.

\item[Matrix Representations] In this thesis, only a small number of matrix representations for iterative methods is presented. Clearly, it would be desirable to find representations of many more methods. Additionally, it might be interesting to find out if this representation reveals new insights about different iterative methods and their relations to each other.
\end{description}

\appendix

\chapter{Matrix Properties}
\label{chap:appendixProperties}

In the following, we define those matrix properties used throughout the thesis that are not self-explanatory. Let $A \in \mathbb{R}^{n \times m}$ be a matrix. The elements of this matrix are denoted as $a_{ij}$ with $i \in \{0, \ldots, n-1\}$ and $j \in \{0, \ldots, m-1\}$.

%
%
%
%

\begin{itemize}
\item[-] Upper diagonal ($\text{\ttfamily UpperDiagonal}$): $a_{ij} = 0$ for $i + 1 \neq j$ with $n = m$. Consider the matrix below as an example.
$$
\left(
\begin{array}{cccc}
0 & a_{01} & 0 & 0 \\ 
0 & 0 & a_{12} & 0 \\ 
0 & 0 & 0 & a_{23} \\ 
0 & 0 & 0 & 0 \\ 
\end{array}
\right)
$$
\item[-] Lower diagonal ($\text{\ttfamily LowerDiagonal}$): $a_{ij} = 0$ for $i - 1 \neq j$ with $n = m$.
\item[-] Diagonal and rectangular ($\text{\ttfamily DiagonalR}$): $a_{ij} = 0$ for $i \neq j$ with $n \neq m$.
\item[-] Upper diagonal and rectangular ($\text{\ttfamily UpperDiagonalR}$): $a_{ij} = 0$ for $i + 1 \neq j$ with $n \neq m$.
\item[-] Lower diagonal and rectangular ($\text{\ttfamily LowerDiagonalR}$): $a_{ij} = 0$ for $i - 1 \neq j$ with $n \neq m$.
\item[-] Upper trapezoidal ($\text{\ttfamily UpperTrapezoidal}$): $a_{ij} = 0$ for $i > j$ with $n < m$. Thus, the following matrix is upper trapezoidal:
$$
\left(
\begin{array}{cccc}
a_{00} & a_{01} & a_{02} & a_{03} \\ 
0 & a_{11} & a_{12} & a_{13} \\ 
0 & 0 & a_{22} & a_{23} \\ 
\end{array}
\right)
$$
\item[-] Lower trapezoidal ($\text{\ttfamily LowerTrapezoidal}$): $a_{ij} = 0$ for $i < j$ with $n > m$.
\item[-] Upper triangular and rectangular ($\text{\ttfamily UpperTriangularR}$): $a_{ij} = 0$ for $i > j$ with $n > m$. An example of such a matrix is shown below.
$$
\left(
\begin{array}{ccc}
a_{00} & a_{01} & a_{02} \\ 
0 & a_{11} & a_{12} \\ 
0 & 0 & a_{22} \\ 
0 & 0 & 0 \\ 
\end{array}
\right)
$$
\item[-] Lower triangular and rectangular ($\text{\ttfamily LowerTriangularR}$): $a_{ij} = 0$ for $i < j$ with $n < m$.
\item[-] Elements on the diagonal are zero ($\text{\ttfamily ZeroDiagonal}$): $a_{ij} = 0$ for $i = j$ with $n = m$.
\end{itemize}


\chapter{Matrix Representations of Iterative Methods}
\label{chap:appendixRepresentations}

In this chapter, a collection of matrix representations of iterative methods is provided. The representations for the Krylov sequence, Steepest Descent, symmetric CG and nonsymmetric CG are modifications of the ones introduced in \cite{eijkhout:CGderivation}. The differences lie in the use of the underline. The representation for BiCG is based on the one for symmetric CG.
\section{Krylov Subspace Methods}
\subsection{Krylov Sequence}

\begin{align*}
P_\text{pre}: \{ &\text{\ttfamily Input}(A) \land \text{\ttfamily Matrix}(A) \land \\
		&\text{\ttfamily Matrix}[\underline{J}] \land \text{\ttfamily LowerDiagonalR}[\underline{J}] \land \\
		&\text{\ttfamily FirstColumnInput}(K) \land \text{\ttfamily Matrix}[K] \}
\end{align*}
\begin{align*}
P_\text{post}: \{ &A \underline{K} = K \underline{J} \\
			&\text{size}(K) = n \times m \}
\end{align*}

\subsection{Steepest Descent}

\begin{align*}
P_\text{pre}: \{ &\text{\ttfamily Input}[A] \land \text{\ttfamily Matrix}[A] \land \text{\ttfamily NonSingular}[A] \land \\
		&\text{\ttfamily Output}[D] \land \text{\ttfamily Matrix}[D] \land \text{\ttfamily Diagonal}[D] \land \\
		&\text{\ttfamily FirstColumnInput}[R] \land \text{\ttfamily Matrix}[R] \land \\
		&\text{\ttfamily ZeroDiagonal}[R^T R J] \land \text{\ttfamily ZeroDiagonal}[R^T R J^T] \land \\
		&\text{\ttfamily FirstColumnInput}[X] \land \text{\ttfamily Matrix}[X] \land \\
		& \text{\ttfamily Matrix}[\underline{I} - \underline{J}] \land \text{\ttfamily LowerTrapezoidal}[\underline{I} - \underline{J}] \}
\end{align*}
\begin{align*}
P_\text{post}:	\{ 	&ARD = R \left( \underline{I} - \underline{J} \right) \\
				&RD = X \left( \underline{I} - \underline{J} \right) \\
				&\| R e_r^T \| < \varepsilon \}
\end{align*}

\subsection{Conjugate Gradient (symmetric)}
\begin{align*}
P_\text{pre}: \{ &\text{\ttfamily Input}[A] \land \text{\ttfamily Matrix}[A] \land \text{\ttfamily NonSingular}[A] \land \text{\ttfamily Symmetric}[A] \land \\
		&\text{\ttfamily Output}[P] \land \text{\ttfamily Matrix}[P] \land \\
		&\text{\ttfamily Output}[D] \land \text{\ttfamily Matrix}[D] \land \text{\ttfamily Diagonal}[D] \land \\
		&\text{\ttfamily FirstColumnInput}[R] \land \text{\ttfamily Matrix}[R] \land \text{\ttfamily Orthogonal}[R] \land \\
		&\text{\ttfamily FirstColumnInput}[\underline{R}] \land \text{\ttfamily Matrix}[\underline{R}] \land \text{\ttfamily Orthogonal}[\underline{R}] \land \\
		&\text{\ttfamily DiagonalR}[R^T \underline{R}] \land \text{\ttfamily DiagonalR}[\underline{R}^T R] \land \\
		&\text{\ttfamily FirstColumnInput}[X] \land \text{\ttfamily Matrix}[X] \land \\
		&\text{\ttfamily Output}[U] \land \text{\ttfamily Matrix}[U] \land \text{\ttfamily UpperDiagonal}[U] \land \\
		& \text{\ttfamily Matrix}[\underline{I} - \underline{J}] \land \text{\ttfamily LowerTrapezoidal}[\underline{I} - \underline{J}] \}
\end{align*}
\begin{align*}
P_\text{post}:	\{ 	&APD = R \left( \underline{I} - \underline{J}  \right) \\
				&P \left( I - U \right) = \underline{R} \\
				&PD = X \left( \underline{I} - \underline{J} \right) \\
				&\| R e_r^T \| < \varepsilon \}
\end{align*}

\subsection{Conjugate Gradient (nonsymmetric)}

\begin{align*}
P_\text{pre}: \{ &\text{\ttfamily Input}[A] \land \text{\ttfamily Matrix}[A] \land \text{\ttfamily NonSingular}[A] \land \\
		&\text{\ttfamily Output}[P] \land \text{\ttfamily Matrix}[P] \land \\
		&\text{\ttfamily Output}[D] \land \text{\ttfamily Matrix}[D] \land \text{\ttfamily Diagonal}[D] \land \\
		&\text{\ttfamily FirstColumnInput}[R] \land \text{\ttfamily Matrix}[R] \land \text{\ttfamily Orthogonal}[R] \land \\
		&\text{\ttfamily FirstColumnInput}[\underline{R}] \land \text{\ttfamily Matrix}[\underline{R}] \land \text{\ttfamily Orthogonal}[\underline{R}] \land \\
		&\text{\ttfamily DiagonalR}[R^T \underline{R}] \land \text{\ttfamily DiagonalR}[\underline{R}^T R] \land \\
		&\text{\ttfamily FirstColumnInput}[X] \land \text{\ttfamily Matrix}[X] \land \\
		&\text{\ttfamily Output}[U] \land \text{\ttfamily Matrix}[U] \land \text{\ttfamily StrictlyUpperTriangular}[U] \land \\
		& \text{\ttfamily Matrix}[\underline{I} - \underline{J}] \land \text{\ttfamily LowerTrapezoidal}[\underline{I} - \underline{J}] \}
\end{align*}
\begin{align*}
P_\text{post}:	\{ 	&APD = R \left( \underline{I} - \underline{J} \right) \\
				&P \left( I - U \right) = \underline{R} \\
				&PD = X \left( \underline{I} - \underline{J} \right) \\
				&\| R e_r^T \| < \varepsilon \}
\end{align*}

\subsection{BiCG}

\begin{align*}
P_\text{pre}: \{ &\text{\ttfamily Input}[A] \land \text{\ttfamily Matrix}[A] \land \text{\ttfamily NonSingular}[A] \land \\
		&\text{\ttfamily Output}[P] \land \text{\ttfamily Matrix}[P] \land \\
		&\text{\ttfamily Output}[\tilde{P}] \land \text{\ttfamily Matrix}[\tilde{P}] \land \\
		&\text{\ttfamily Output}[D] \land \text{\ttfamily Matrix}[D] \land \text{\ttfamily Diagonal}[D] \land \\
		&\text{\ttfamily FirstColumnInput}[R] \land \text{\ttfamily Matrix}[R] \land \\
		&\text{\ttfamily FirstColumnInput}[\underline{R}] \land \text{\ttfamily Matrix}[\underline{R}]  \land \\
		&\text{\ttfamily FirstColumnInput}[\tilde{R}] \land \text{\ttfamily Matrix}[\tilde{R}] \land  \\
		&\text{\ttfamily FirstColumnInput}[\underline{\tilde{R}}] \land \text{\ttfamily Matrix}[\underline{\tilde{R}}] \land \\
		&\text{\ttfamily Diagonal}[R^T \tilde{R}] \land \text{\ttfamily Diagonal}[\underline{R}^T \underline{\tilde{R}}] \land \\
		&\text{\ttfamily DiagonalR}[\underline{R}^T \tilde{R}] \land \text{\ttfamily DiagonalR}[R^T \underline{\tilde{R}}] \land \\
		&\text{\ttfamily FirstColumnInput}[X] \land \text{\ttfamily Matrix}[X] \land \\
		&\text{\ttfamily Output}[U] \land \text{\ttfamily Matrix}[U] \land \text{\ttfamily UpperDiagonal}[U] \land \\
		& \text{\ttfamily Matrix}[\underline{I} - \underline{J}] \land \text{\ttfamily LowerTrapezoidal}[\underline{I} - \underline{J}] \}
\end{align*}
\begin{align*}
P_\text{post}:	\{ 	&APD = R \left( \underline{I} - \underline{J} \right) 
				&&A^T \tilde{P}D = \tilde{R} \left( \underline{I} - \underline{J} \right) \\
				&P \left( I - U \right) = \underline{R}
				&&\tilde{P} \left( I - U \right) = \underline{\tilde{R}} \\
				&PD = X \left( \underline{I} - \underline{J} \right) 
				&&\| R e_r^T \| < \varepsilon \}
\end{align*}

\section{Stationary Iterative Methods}
\label{sec:stationaryIterative}

Let $Ax = b$ be a linear system. $e$ is a column vector where all entries are one. We write $A$ as $A = D - L - U$, where $D$ contains the entries on the main diagonal of $A$, $L$ contains the entries below the main diagonal and $U$ the ones above the main diagonal. The representations for the Gauss-Seidel, Jacobi and Successive Overrelaxation method are based on the descriptions (in indexed notation) in \cite{barrett:templates}. The one for the Richardson iteration is based on the the description (in indexed notation) in \cite{eijkhout:CGderivation}.
\subsection{Gauss-Seidel Method}
\begin{align*}
P_\text{pre}: \{ &\text{\ttfamily Input}[D] \land \text{\ttfamily Matrix}[D] \land \text{\ttfamily Diagonal}[D] \land \\
		&\text{\ttfamily Input}[L] \land \text{\ttfamily Matrix}[L] \land \text{\ttfamily LowerTriangular}[L] \land \\
		&\text{\ttfamily Input}[U] \land \text{\ttfamily Matrix}[U] \land \text{\ttfamily UpperTriangular}[U] \land \\
		&\text{\ttfamily Input}[b] \land \text{\ttfamily Vector}[b] \land \\
		&\text{\ttfamily FirstColumnInput}[X] \land \text{\ttfamily Matrix}[X] \}
\end{align*}
\begin{align*}
P_\text{post}:	\{ 	&(D - L) X \underline{J} = U \underline{X} + be^T \\
				& \| X e_r^T - X e_{r-1}^T \| < \varepsilon \}
\end{align*}
\subsection{Jacobi Method}
The precondition of the Jacobi method is identical to the one for the Gauss-Seidel method.
\begin{align*}
P_\text{post}:	\{ 	&D X \underline{J} = ( L + U ) \underline{X} + be^T \\
				& \| X e_r^T - X e_{r-1}^T \| < \varepsilon \}
\end{align*}
\subsection{Successive Overrelaxation Method}
\begin{align*}
P_\text{pre}: \{ &\text{\ttfamily Input}[D] \land \text{\ttfamily Matrix}[D] \land \text{\ttfamily Diagonal}[D] \land \\
		&\text{\ttfamily Input}[L] \land \text{\ttfamily Matrix}[L] \land \text{\ttfamily LowerTriangular}[L] \land \\
		&\text{\ttfamily Input}[U] \land \text{\ttfamily Matrix}[U] \land \text{\ttfamily UpperTriangular}[U] \land \\
		&\text{\ttfamily Input}[b] \land \text{\ttfamily Vector}[b] \land \\
		&\text{\ttfamily Input}[\omega] \land \text{\ttfamily Scalar}[\omega] \\
		&\text{\ttfamily FirstColumnInput}[X] \land \text{\ttfamily Matrix}[X] \}
\end{align*}
\begin{align*}
P_\text{post}:	\{ &(D - \omega L) X \underline{J} = ( \omega U + (1 - \omega) D ) \underline{X} + \omega b e^T \\
			& \| X e_r^T - X e_{r-1}^T \| < \varepsilon\}
\end{align*}
%
%

\subsection{Richardson Iteration}
\begin{align*}
P_\text{pre}: \{ &\text{\ttfamily Input}[A] \land \text{\ttfamily Matrix}[A] \land \text{\ttfamily NonSingular}[A] \land \\
		&\text{\ttfamily Input}[b] \land \text{\ttfamily Vector}[b] \land \\
		&\text{\ttfamily Input}[\alpha] \land \text{\ttfamily Scalar}[\alpha] \\
		&\text{\ttfamily FirstColumnInput}[X] \land \text{\ttfamily Matrix}[X] \land \\
		&\text{\ttfamily Output}[R] \land \text{\ttfamily Matrix}[R] \}
\end{align*}
\begin{align*}
P_\text{post}:	\{ & \alpha R = X \left( \underline{I} - \underline{J} \right) \\
			& R = A \underline{X} - be^T \\
			& \| R e_r^T \| < \varepsilon \}
\end{align*}


\cleardoublepage

\addcontentsline{toc}{chapter}{Bibliography}
\bibliographystyle{plain}
\bibliography{MS_thesis_bibliography}

\end{document}